\documentclass[11pt]{amsart}
\usepackage{verbatim} 
\usepackage[utf8]{inputenc}
\usepackage{mathbbol,mathrsfs,amsfonts,amssymb,amsthm,stmaryrd}
\usepackage{amsmath}
\usepackage[all]{xypic}
\usepackage{mathpazo}
\usepackage[mathcal]{euscript}
\usepackage{graphics,supertabular,color}
\usepackage{bm}
\usepackage[colorlinks,hyperindex]{hyperref}
\hypersetup{linkcolor=blue, citecolor=cyan}

\newcommand{\R}{\mathbb{R}}

\newcommand{\Z}{\mathbb{Z}}
\newcommand{\T}{\mathbb{T}}
\newcommand{\U}{{\rm U}}
\newcommand{\PU}{{\rm PU}}
\newcommand{\HH}{\mathcal{H}}
\newcommand{\K}{\mathcal{K}}
\newtheoremstyle{ansgarthmstyle}{}{}{\it}{}{\bf}{}{0.4cm}{}
\theoremstyle{ansgarthmstyle}
\newtheorem{exa}{Example}[section]
\newtheorem{defi}{Definition}[section]
\newtheorem{prop}{Proposition}[section]
\newtheorem{rem}{Remark}[section]
\newtheorem{lem}{Lemma}[section]
\newtheorem{thm}{Theorem}[section]
\newtheorem{cor}{Corollary}[section]
\newenvironment{pf}{\begin{proof}[\textrm{{\bf Proof\ :\!}}]} {\end{proof}}

\begin{document}

\title{Equivariant T-Duality of Locally Compact Abelian Groups}

\author{Ansgar Schneider}

\thanks{ansgar.schneider@mathematik.uni-regensburg.de\hspace{12cm} }

\maketitle

\begin{center}
{\sc Abstract}
\end{center}
\noindent
Equivariant T-duality triples of locally compact abelian 
groups are considered. The motivating example dealing with the group 
$\R^n$ containing a lattice $\Z^n$ comes with
an isomorphism in twisted equivariant K-theory. 

\tableofcontents

\section{Introduction}

\noindent
In the mathematical study of T-duality
the notion of T-duality triples has been 
introduced in \cite{BS}. A key feature of 
T-duality triples is that they come 
with an isomorphism in twisted K-theory.
In \cite{BS2}  T-duality over orbifolds has been studied, and a 
corresponding isomorphism in Borel K-theory was proven.
However, they did not prove the corresponding isomorphism 
in K-theory. At least partially, this work fills that gap.

In \cite{Sch} it has been shown that the approach to T-duality 
of \cite{BS} using T-duality triples is equivalent to the $C^*$-algebraic approach 
of \cite{MR}. In the latter  Connes' Thom isomorphism \cite{Con}
turns out to be  the key tool to identify the K-theories of T-dual objects. 
The application in K-theory given in this work is also 
based on Connes' Thom isomorphism, however, 
the techniques used herein are quite different 
from \cite{MR}.

In this work we study T-duality triples (we will call them topological triples)
over the singleton space equipped with an action of a finite
group $\Gamma$. Let us start with its precise definition.
A topological triple is a  $\Gamma$-equivariant
commutative diagram of trivialisable principal fibre bundles    
\begin{equation*}
\xymatrix{
& 
& P \times \widehat E
\ar[rd]\ar[ld]&
&E\times\widehat P\ar[rd]\ar[ld]\ar[ll]_\cong
&\\
&P \ar[rd]_{\PU}&&
E\times \widehat E
\ar[rd]\ar[ld]&\quad\qquad& \widehat P\ar[ld]^{\PU}\\
&&E\ar[rd]_{G/N}&&\widehat E\ar[ld]^{\widehat G/N^\perp}&\\
&&&\ast,&&
}
\end{equation*}
where the structure groups of the bundles are as indicated.
Here $\PU$ is the projective unitary group of some separable 
Hilbert space, $G$ is a locally compact abelian group with 
discrete and cocompact subgroup $N$, and $\widehat G$ is the dual group 
of $G$ containing $N^\perp$, the annihilator of $N$.
Moreover, the  isomorphism on top of the diagram must 
satisfy a certain condition:
By trivialisability of the diagram, for a trivialisation the top isomorphism 
defines a function $G/N\times \widehat G/N^\perp\to \PU$, 
thus a class in $\check H^1(G/N\times \widehat G/N^\perp,\underline{\U(1)})$.
The requirement is that there exists a trivialisation such that
this class is contained in the subgroup of 
$\check H^1(G/N\times \widehat G/N^\perp,\underline{\U(1)})$
which is generated by the Poincaré class. 
The Poincaré class is the class of the canonical line bundle 
over $G/N\times \widehat G/N^\perp.$
If for an integer $L$ the class of the top isomorphism is 
$L$ times the Poincaré class, then we call $L$ the order
of the topological triple.

To a topological triple we can associate two $C^*$-algebras 
$C^*(E,P)$ and $C^*(\widehat E,\widehat P)$ which are 
the crossed product $C^*$-algebras of the group $\Gamma$
with the $C^*$-algebras of sections of the associated 
$C^*$-bundles $P\times_\PU\K\to E$ and 
$\widehat P\times_\PU\K\to \widehat E$. Here $\K$ is
the $C^*$-algebra of compact operators.
It is the principal aim of this paper to understand how the two $C^*$-algebras 
$C^*(E,P)$ and $C^*(\widehat E,\widehat P)$ are related to each other.
{The basic idea is that one can construct out of the data of
a topological triple a $G$-action  on $P$ which commutes with the $\Gamma$-action
and which covers the $G/N$ action on  $E$. 
Then there is an induced $G$-action on $C^*(E,P)$, and we 
show that the crossed product of $G$ with $C^*(E,P)$ is Morita equivalent 
to $C^*(\widehat E,\widehat P)$.
For $\Gamma$ being the trivial group this statement reduces to a well-known 
Result proven by Rieffel, but for a non-trivial group $\Gamma$ the situation is more involved 
as we are faced with the given $\Gamma$-actions of the topological triple.
}
In particular, {for the groups $G=\R^n$ and $N=\Z^n$,}
if the order of the topological triple is $L=1$, 
we establish an isomorphism in  K-theory between
$K^*_{\Gamma,P}(E):= K_*(C^*(E,P){)}$ and 
$K^*_{\Gamma,\widehat P}(\widehat E):= K_*(C^*(\widehat E,\widehat P){)}$.
These groups are the equivariant twisted K-theories of $E$ and $\widehat E$
with twists given by $P$ and $\widehat P$, respectively.
Alternatively, { a $\Gamma$-equivariant $\PU$-bundle
$P\to E$ is the same amount of data as a $\PU$-bundle  over
the transformation groupoid $\Gamma\ltimes E$, and thus, as shown in \cite[Sec. 4]{TXC},}
these are the twisted K-theories of the transformation
groupoids $\Gamma\ltimes E$ and $\Gamma\ltimes \widehat E$
with twists $P$ and $\widehat P$, respectively.
Turning from a groupoid to the stack presented by it, 
one may also take these groups as the twisted K-theories of the quotient stacks 
$[E/\Gamma]$ and $[\widehat E/\Gamma]$
with twists $P$ and $\widehat P$, respectively.
\\
\\
We give a short overview of this work.

We introduce the notion of pairs in section \ref{subsecTheCategoryOfPairs}.
A pair consists of  principal fibre bundle over the singelton space 
with structure group $G/N$ and a projective unitary
trivialisable principal bundle on its total space, where a finite
group $\Gamma$ acts on all spaces by bundle automorphisms.

If on a pair there is an additional action of $G$, 
we call these data a dynamical triple. A dynamical triple defines in 
a natural way a class in $H^2(N,\U(1))$ called the Mackey obstruction. 
If the Mackey obstruction of a dynamical triple
vanishes we call the triple dualisable (section \ref{SecDynT}), and
we prove in section \ref{SecTCODDT} a classification theorem
for dualisable dynamical triples. This is the key tool for our  
first important result, namely, the existence of a natural duality 
theory (section \ref{SecDualityTheoryODT}) of dualisable 
dynamical triples (on the level of equivalence classes).
The dual objects are given by replacing everything by its dual
in the sense of the duality theory of abelian groups 
(section \ref{SecDualPairsAndTriples}).

In section \ref{SecTDATCAODT} we show that the 
duality theory of section \ref{SecDualityTheoryODT}
is the same duality theory which one obtains if one concerns
the associated $C^*$-dynamical systems and uses the duality
theory of abelian crossed products.

In section \ref{SecTopologicalTriples} we introduce topological 
triples the objects of our main interest.
The second important result of this work, the 
classification theorem of topological triples
(section \ref{SecTheClassifOfET}), 
enables us to point out the relation between 
dynamical triples and topological triples:
there is a natural bijection (section \ref{SecFromDToTT}) 
between (the equivalence classes of)
dualisable dynamical triples and 
(the equivalence classes of) regular
topological triples of order $1$. 
Regular topological triples (section \ref{SecRegET}) are a 
subclass of all topological {triples} defined by 
means of their classification.

The main result of this work is
then stated in section \ref{SecAAppliEquKT},
where all partial results are puzzled together. 
If we fix the  groups $G=\R^n$ and $N=\Z^n$, then we show first 
that all topological triples are regular.
Thus, there is a natural bijection between topological triples
of order $1$ and dualisable dynamical triples.
By the duality theory of section \ref{SecDualityTheoryODT} 
and Connes' Thom isomorphism, we then obtain the 
mentioned isomorphism in equivariant twisted K-theory. 
\\
\\
\\

\hspace{4,23cm}{\bf Acknowledgement}
\nopagebreak
\\
\nopagebreak
\\
\nopagebreak
\noindent
The author is grateful to Ulrich Bunke
for  several discussions and for drawing 
his attention to the problem.

\section{Notation}

\noindent
{\bf 1.}
By  $\Gamma$ we will always denote a finite group. 
By $G$ will always denote a second countable,
Hausdorff, locally compact abelian group which
has a discrete, cocompact subgroup $N$, i.e.
the quotient $G/N$ is compact. The dual group
${\rm Hom}(G,\U(1))$ of $G$ is denoted by $\widehat G$.
It contains the cocompact subgroup 
$N^\perp:=\{\hat g\in\widehat G\ |\  \hat g|_N=1\}$.
Recall the classical isomorphisms of 
Pontrjagin duality \cite{Ru}
\begin{eqnarray*}
G\cong \widehat {\widehat G},\quad \widehat N\cong {\widehat{\textcolor{black}{G}}}/N^{{\perp}},\quad \widehat{G/N}\cong N^\perp.
\end{eqnarray*}
For all pairings between an abelian group $H$ and its dual $\widehat H$ 
we are going to use the bracket notation 
\begin{eqnarray*}
\langle\ , \ \rangle:\widehat H\times H\to \U(1).
\end{eqnarray*}
As an example, take $g\in G$ which is mapped 
to $z\in G/N$ by the quotient map, and let $n\in N^\perp\subset\widehat G $. 
Then we have the identity
\begin{eqnarray*}
\langle  n,g\rangle=\langle n,z\rangle,
\end{eqnarray*}
where the first pairing is between $\widehat G$ and $G$
and the second is between $N^\perp\cong \widehat{G/N}$ and $G/N$.
\\
\\
{\bf 2.}
Once and for all we fix Borel sections $\sigma:G/N\to G$
and $\hat\sigma:\widehat G/N^\perp\to \widehat G$ of 
the quotient maps $G\to G/N$ and $ \widehat G\to \widehat G/N^\perp$.
\\
\\
{\bf 3.}
If $M$ is an abelian group and module over the finite group 
$\Gamma$, we have the group cohomological chain
complex 
\begin{eqnarray*}
\cdots\stackrel{\delta}{\longrightarrow} C^n(\Gamma,M)
\stackrel{\delta}{\longrightarrow} C^{n+1}(\Gamma,M)
\stackrel{\delta}{\longrightarrow} \cdots
\end{eqnarray*}
where $C^n(\Gamma,M)$ is the set of all maps from $\Gamma^n$ to $M$.
We use the following (non standard) 
convention for the boundary operator $\delta$
\begin{eqnarray*}
(\delta m)(a_0,\dots,a_n)&:=& (a_0)^{-1}\cdot m(a_1,\dots,a_n)
-m(a_1a_0,a_2,\dots,a_n)\\
&&+ m(a_0,a_2a_1,\dots,a_n) -\cdots+\cdots\\
&&+(-1)^{{n}} m(a_0,\dots, a_n a_{n-1})
+(-1)^{{n+1}}m(a_0,\dots,a_{n-1})
\end{eqnarray*}
which fits better in our formulas than the usual one.
We will write the group law in $M$ multiplicative{ly} whenever 
the modul{e} is  $\U(1)$ (trivial module structure)
or a space of $\U(1)$-valued functions on $G/N$ or $\widehat G/N^\perp$.
In the latter case the action of $a\in \Gamma$ on a function $f$ on $G/N$
is of the form\footnote{
Observe that the minus sign of the action cancels 
the minus sign of $(a_0)^{-1}$ in the boundary operator.
}
$(a\cdot f)(z):=f(-\chi(a)+z)$, where
$\chi:\Gamma\to G/N$ is a homomorphism.
Dually, $\Gamma$ acts on a space of functions
on $\widehat G/N^\perp$ by a homomorphism
$\hat\chi:\Gamma\to \widehat G/N^\perp$.
The corresponding boundary operators are denoted 
$\delta_\chi$ or $\delta_{\hat\chi}$.
\\
\\
{\bf 4.}
We will also consider modules over the groups $G$ and $\widehat G$,
then we use the notation $d$ and $\hat d$ for the boundary operator
which is given by the same algebraic formula as $\delta$ above.
On a space of functions on $G/N$ the group $G$ acts in obvious way
in the arguments of the functions by (the negative of) the canonical homomorphism 
$G\to G/N$. Similarly, $\widehat G$ acts on functions on $\widehat G/N^\perp$.
\\
\\
{\bf 5.}
$\U(\HH)$ is the unitary 
group of some  separable Hilbert space
$\HH$, it is equipped with the strong operator topology.
The quotient by its center $\U(1)$ is the projective 
unitary group $\PU(\HH)$.
We denote the quotient map by 
\begin{eqnarray*}
{\rm Ad}:\U(\HH)\to \PU(\HH).
\end{eqnarray*}
The action of $\PU(\HH)$ on $\U(\HH)$ by conjugation
is denoted by
\begin{eqnarray*}
\PU(\HH)\times\U(\HH)&\to& \U(\HH)\\
(v,U)&\mapsto&v[U]:=V U V^{-1},
\end{eqnarray*}
where $V\in {\rm Ad}^{-1}(v)$ is some pre-image
of $v$. 
\\
\\
{\bf 6.}
The  action of $\PU(\HH)$ on the $C^*$-algebra
of compact operators $\K(\HH)$ by conjugation 
is also denoted with squared 
brackets: $v[K]$, for $K\in\K(\HH), v\in\PU(\HH)$. 
This defines an isomorphism between $\PU(\HH)$ and the 
$C^*$-automorphism group of $\K(\HH)$.
\\
\\
{\bf 7.}
A $\U(\HH)$-valued or $\U(1)$-valued Borel function 
on  any of the groups $\Gamma, G, G/N,\dots$ gives rise 
to a unitary multiplication operator on the corresponding $L^2$-space. 
By abuse of notation, we denote  the function and its multiplication
operator by same symbol. 
E.g. consider $G/N\ni z\mapsto\langle g, \sigma(z)\rangle\in \U(1)$, 
then 
\begin{eqnarray*}
\langle g,\sigma(\_)\rangle: L^2(G/N)&\to& L^2(G/N)\\
F&\mapsto& (z\mapsto \langle g,\sigma(z)\rangle F(z)).
\end{eqnarray*}
To make the reader more familiar with that notation 
we note here that by composition with ${\rm Ad}$ we obtain 
a map
\begin{eqnarray*}
G&\to& \PU(L^2(G/N))\\
g&\mapsto& {\rm Ad}(\langle g,\sigma(\_)\rangle).
\end{eqnarray*}

\section{Pairs}
\label{subsecTheCategoryOfPairs}

\noindent
Let $\Gamma$ be any finite group, 
and  let $\HH$ be any separable 
 Hilbert space.

\begin{defi}\label{DefiPairOverB}
A  {\bf  pair} $(P,E)$ 
(over the one point space $\ast$) 
with underlying Hilbert space $\HH$
is a $\Gamma$-equivariant sequence 
\begin{eqnarray*}
\begin{array}{c}
\xymatrix{
\ar@(dl,ul)[]^\Gamma 
}
\end{array}
\hspace{-0.4cm}
\left(
\begin{array}{c}
\xymatrix{ 
P\ar[d]&\hspace{-0,9cm}\rotatebox {90}{$\circlearrowleft$}\ \PU(\HH) \\ 
E\ar[d]&\hspace{-1,2cm}\rotatebox {90}{$\circlearrowleft$}\ G/N \\ 
\ast&
}
\end{array}
\right)
\end{eqnarray*}
of trivialisable principal fibre bundles, 
where  $\Gamma$ acts (from the left) by bundle automorphisms.
\end{defi}

Due to the triviality condition on $P$ 
we have a diagram of bundle isomorphisms
\begin{equation}\label{DiagOfTrivPairOverPt} 
\xymatrix{
G/N\times\PU(\HH)\ar[r]\ar[d]&P\ar[d]\\
G/N\ar[d]\ar[r]&E\ar[d]\\
\ast\ar[r]^=&\ast .
}
\end{equation}
We call any choice of such bundle isomorphisms  a {\bf chart} for the pair.
Pullback of the $\Gamma$-actions along 
the chart induces $\Gamma$-actions on 
$G/N$ and  $G/N\times\PU(\HH)$
which are given
by a homomorphism 
\begin{equation}\label{EqTheHomoPropOfChi}
\chi:\Gamma\to G/N
\end{equation}
and a group cohomological 1-cocycle 
$\lambda:\Gamma\to C(G/N,\PU(\HH))$,
i.e.
\begin{equation}\label{EqTheCocyclePropertyOfLambda}
\lambda(ba,z)=\lambda(b,z+\chi(a))\ \lambda(a,z),
\end{equation}
for $a,b\in\Gamma$ and $z\in G/N$.
Equivalently, one can consider  $(\chi,\lambda)$ 
as a single homomorphism
$$
\chi\times\lambda: \Gamma\to G/N\ltimes C(G/N,\PU(\HH)),
$$
where $G/N\ltimes C(G/N,\PU(\HH))$ is the semi-direct product.

It is obvious that any  tuple $(\chi,\lambda)$ with these properties 
defines the structure of a pair on $G/N\times\PU(\HH)\to G/N\to\ast$, and in this sense 
we also call such a tuple $(\chi,\lambda)$ a {\bf pair}.

A {\bf morphism} $(\varphi,\vartheta,\theta)$  from a   pair
$P{\to} E\to *$ with underlying Hilbert space $\HH$ 
to a pair  $P'{\to} E'\to *$ with underlying Hilbert space $\HH'$  is
a commutative diagram of $\Gamma$-equivariant bundle isomorphisms
\begin{equation}\label{DiagramMorphismOfPairs}
\xymatrix{
P\ar[r]^\vartheta        \ar[d]  & \varphi^*P'\ar[d]&\\
E\ar[r]^{\theta} \ar[d]  & E'\ar[d]& \\
{\ast}\ar[r]^{=}             &\ast,      &
}
\end{equation}
where $\varphi$ is the
mod-$\U(1)$ class of a unitary isomorphism $\HH\to \HH'$, and
$\varphi^*P'$ is the $\PU(\HH)$-bundle 
$P'\times_{\PU(\HH')}\PU(\HH)$,
where
$\PU(\HH')$ acts on $\PU(\HH)$ 
via 
$\varphi^*:\PU(\HH')\to \PU(\HH)$.
Pairs 
and their morphisms form a category;
composition of morphisms $(\varphi,\vartheta,\theta)$ and $(\varphi',\vartheta',\theta')$ 
is just component-wise composition 
$(\varphi'\circ \varphi,\vartheta'\circ \vartheta,\theta'\circ \theta).$
The resulting category is  a groupoid, 
i.e. every  morphism is an isomorphism. 

If $\HH_0$  is a separable Hilbert space and
$P$ is a principal $\PU(\HH)$-bundle, we use the notation 
\begin{eqnarray}
P_{\HH_0} &:=&\PU(\HH_0)\otimes P\nonumber\\
&:=& P\times_{\PU(\HH)}\PU(\HH_0\otimes \HH).
\label{EqStablilisedPUBundleP}
\end{eqnarray}
for the associated (stabilised) bundle.
We call  
two pairs 
$(P,E)$ and $(P',E')$ with underlying Hilbert 
spaces $\HH$ and $\HH'$  
{\bf stably isomorphic} if there exists a separable Hilbert
space $\HH_0$ such that
the pairs $(P_{\HH_0},E)$ and $(P'_{\HH_0},E')$
are isomorphic.
 
We call two pairs 
$(P,E)$ and $(P',E')$ {\bf outer conjugate}
if they are isomorphic up to a unitary cocycle,
i.e. if there are morphisms of pairs 
\begin{equation}
\xymatrix{
P\ar[r]        \ar[d]  & G/N\times\PU(\HH)\ar[d]
&
\hspace{-1.3cm} \rotatebox {90}{$\circlearrowleft$}\ \lambda&
P'\ar[r]\ar[d]&G/N\times\PU(\HH)\ar[d]&
\hspace{-1.3cm} \rotatebox {90}{$\circlearrowleft$}\ \lambda'
\\
E\ar[r] \ar[d]  & G/N\ar[d]
& 
\hspace{-2.8cm} \rotatebox {90}{$\circlearrowleft$}\ \chi
&
E'\ar[r]\ar[d]&G/N\ar[d]&
\hspace{-2.8cm} \rotatebox {90}{$\circlearrowleft$}\ \chi'
\\
{\ast}\ar[r]^{=}             &\ast ,  &    &\ast\ar[r]^{=}&\ast,&
}
\end{equation}
then the induced  pairs $(\chi,\lambda)$ and $(\chi',\lambda')$
satisfy $\chi=\chi'$ and 
$\lambda'=\lambda\ {\rm Ad}(l),
$
where  
$l:\Gamma\times G/N\to\U(\HH)$
is a continuous function
which satisfies the unitary 
cocycle condition
\begin{eqnarray}\label{EqOaCOET}
l(ba,z)=\lambda(a,z)^{-1}[l(b,z+\chi(a))]\ l(a,z).
\end{eqnarray}

We call two  pairs 
$(P,E)$ and $(P',E')$ {\bf stably outer conjugate}
if there exists a separable Hilbert
space $\HH_0$ such that
the pairs $(P_{\HH_0},E)$ and $(P'_{\HH_0},E')$
are outer conjugate.

It is elementary to check that the notions of
isomorphism, stable isomorphism, outer conjugation and
stable outer conjugation are equivalence relations which 
may be arranged in a diagram of implications
\begin{eqnarray*}
\xymatrix{
{\begin{array}{c} 
\textrm{ isomorphism} \\ 
\textrm{of pairs} 
\end{array}}
\ar@{=>}[r]\ar@{=>}[d]&
{\begin{array}{c} 
\textrm{outer conjugation} \\ 
\textrm{of pairs} 
\end{array}}
\ar@{=>}[d]\\
{\begin{array}{c} 
\textrm{stable isomorphism} \\ 
\textrm{of pairs} 
\end{array}}
\ar@{=>}[r]&
{\begin{array}{c} 
\textrm{stable outer conjugation} \\ 
\textrm{of pairs} 
\end{array}}.
}
\end{eqnarray*}

\begin{exa}\label{ExSOC}
Let $(\chi,\lambda)$ be any pair with 
underlying Hilbert space $\HH$, 
and let $\rho:\Gamma\to\U(L^2(\Gamma))$
be the right regular representation, then
$( \chi,\lambda\otimes{\rm Ad}(\rho))$
is a pair with underlying Hilbert space
$\HH\otimes L^2(\Gamma)$ and
stably outer conjugate to 
$(\chi,\lambda)$.
\end{exa}
\begin{pf}
Firstly,  $(\chi,\lambda)$ and 
$(\chi,\lambda\otimes{\rm Ad}(\Eins_{L^2(\Gamma)}) )$
are stably isomorphic, for
let $\HH_0$ be another infinite dimensional Hilbert space, 
then we have an isomorphism of hilbert spaces
$ L^2(\Gamma)\otimes\HH_0\to\HH_0$
which induces an isomorphism
of pairs
\begin{eqnarray*}
\Big(\chi,\big(\lambda\otimes{\rm Ad}(\Eins_{L^2(\Gamma)})\big)
\otimes{\rm Ad}(\Eins_{\HH_0}) \Big)\to 
\Big(\chi,\lambda\otimes{\rm Ad}(\Eins_{\HH_0}) \Big)
\end{eqnarray*}
 with underlying Hilbert spaces 
 $ L^2(\Gamma)\otimes\HH_0$ and $\HH_0$.

Secondly, $(\chi,\lambda\otimes{\rm Ad}(\Eins_{L^2(\Gamma)}) )$
and 
$(\chi,\lambda\otimes{\rm Ad}(\rho) )$
are outer conjugate, for we can define 
\begin{eqnarray*}
l(a,z):=\Eins_\HH\otimes\rho(a)
\end{eqnarray*}
which clearly satisfies (\ref{EqOaCOET}).
\end{pf}

Let us denote by ${\rm Par}$ the set valued contravariant functor
that sends a finite group $\Gamma$ to the set of stable outer conjugation classes 
of pairs, i.e.
\begin{eqnarray}
{\rm Par}(\Gamma):= \{\textrm{ pairs } \}\big/_{\rm st. out.conj.}\ ,
\end{eqnarray}
and if $f:\Gamma'\to \Gamma$ is a homomorphism of finite groups,
then pullback defines a map
$f^*:{\rm Par}(\Gamma)\to {\rm Par}(\Gamma')$.

Let $(P,E)$ be a pair.
We choose a chart 
to obtain an induced  
$(\chi,\lambda)$ pair.
Note that the homomorphism 
$\chi:\Gamma\to G/N$
is independent of the choice of
the chart.
So we obtain another set valued functor 
\begin{eqnarray}
{\rm Par}(\Gamma,\chi):= \{\textrm{  pairs with fixed } \chi\}\big/_{{\rm 
st.out.conj.} }\ ,
\end{eqnarray}
defined 
on the category of homomorphisms
$\chi:\Gamma\to G/N$ whose  morphisms 
are commutative diagrams of group homomorphisms
$$
\xymatrix{
\Gamma'\ar[r]\ar[d]^{\chi'}&\Gamma\ar[d]^\chi\\
G/N\ar@{=}[r]&G/N.
}
$$ 
The set ${\rm Par}(\Gamma,\chi)$ becomes an abelian group
by tensor product and complex conjugation of projective unitary 
bundles.

Let $L\in\Z$ be any integer. If $(L\chi,\lambda)$
is a pair, then it is immediate that also $(\chi,L^*\lambda)$ is  a pair, 
where $(L^*\lambda)(a,z):=\lambda(a,Lz)$.
This leads to a natural map
\begin{eqnarray}\label{EqTheLMapOfPairs}
L^*:{\rm Par}(\Gamma,L\chi)\to {\rm Par}(\Gamma,\chi)
\end{eqnarray}
which we  define by sending 
 $(L\chi,\lambda)$ to 
$(\chi,L^*\lambda)$.

\section{Pairs, their $C^*$-Algebras and K-Theory}

\noindent
For any pair $(P,E)$ we can consider the $C^*$-algebra
of continuous sections $\Gamma(E, F)$
for the associated $C^*$-bundle 
$F:=P\times_{\PU(\HH)}\K(\HH)$
with the compacts as fibre.
Naturally, this $C^*$-algebra together with its inherited action $\alpha$ 
of $\Gamma$ defines a $C^*$-dynamical system 
$(\Gamma,\alpha,\Gamma(E,F))$.
We denote by $C^*(E,P):=\Gamma\ltimes\Gamma(E,F)$
the corresponding crossed product $C^*$-algebra \cite{Pe}.
As $\Gamma$ is finite, we know by \cite{Ju}
that the equivariant K-theory of $\Gamma(E,F)$
is given by the K-theory of its crossed product,
\begin{eqnarray*}
K_*^\Gamma(\Gamma(E,F))\cong K_*(\Gamma \ltimes \Gamma(E,F)).
\end{eqnarray*}
In this paper we take this last expression as 
definition for twisted equivariant K-theory:  
\begin{defi}
The {\bf twisted equivariant K-theory} of a pair
$(P,E)$ is 
\begin{eqnarray*}
K^*_{\Gamma,P}(E):=K_*(C^*(E,F)).
\end{eqnarray*}
\end{defi}
We also use the notation $C^*(\chi,\lambda)$ and $K^*_{\chi,\lambda}(G/N)$
for the crossed product $C^*$-algebra and the 
twisted equivariant K-theory of the pair $(\chi,\lambda)$, respectively.  
 
Now, consider two pairs $(P,E)$ and $(P',E')$
such that their classes match in ${\rm Par}(\Gamma)$, 
i.e. they are stably outer conjugate. 
Then the $C^*$-dynamical system
$(\Gamma,\alpha',\Gamma(E',F'))$
is 
Morita equivalent to 
$(\Gamma,\alpha,\Gamma(E,F))$
\cite{Com}, and therefore we know in particular 
that there is an isomorphism in twisted equivariant K-theory
\begin{eqnarray*}
K^*_{\Gamma,P}(E) \cong
K^*_{\Gamma,P'}(E').
\end{eqnarray*}

The sections  $\Gamma(G/N,G/N\times \K(\HH))$ of a 
pair $(\chi,\lambda)$ can be 
identified with the continuous functions $C(G/N,\K(\HH))$,
and for a function
 $f:G/N\to\K(\HH)$
the induced action of $a\in\Gamma$ is  given by
\begin{eqnarray*}
(\alpha_{(\chi,\lambda)}(a) f )(z)=\lambda(a,z)^{-1}[f(z+\chi(a))].
\end{eqnarray*}
Now, let  $L\in\Z$ be any integer, and let $(L\chi,\lambda)$ be a pair.
Then we  have a morphism 
of $C^*$-dynamical systems
\begin{eqnarray}\label{EqDochNochNeNumma}
\Big(\Gamma,\alpha_{(L\chi,\lambda)},C(G/N,\K(\HH))\Big)&\to&
\Big(\Gamma,\alpha_{(\chi,L^*\lambda)}, C(G/N,\K(\HH))\Big)\\
f&\mapsto& \big( z\mapsto f(Lz)\big) \nonumber
\end{eqnarray}
This map induces maps between
the  crossed products
$L^*:C^*(L\chi,\lambda)\to C^*(\chi,L^*\lambda)$
and in K-theory 
$L^\#:K^*_{L\chi,\lambda}(G/N)\to K^*_{\chi,L^*\lambda}(G/N)$.
For later purpose we note that  $L^\#$ 
fits into the six term exact sequence 
\begin{eqnarray}\label{DiagMyFirstSixPack}
\xymatrix{
K^1_{\chi,L^*\lambda}(G/N)\ar[r]&K_0(C(L^*))\ar[r]&K^0_{L\chi,\lambda}(G/N)\ar[d]^{L^\#}\\
K^1_{L\chi,\lambda}(G/N)\ar[u]^{L^\#}&K_1(C(L^*))\ar[l]&K^0_{\chi,L^*\lambda}(G/N)\ar[l]
}
\end{eqnarray}
which is induced by the exact sequence 
\begin{eqnarray*}
\Eins\to C_0( (0,1),C^*(\chi,L^*\lambda))\to C(L^*)\to 
C^*(L\chi,\lambda)\to\Eins,
\end{eqnarray*}
where $C(L^*)$ is the mapping cone of $L^*:
C^*(L\chi,\lambda)\to C^*(\chi,L^*\lambda)$.

\section{Dynamical Triples}
\label{SecDynT}

\noindent
Let $(P,E)$ be a pair.
The quotient map $G \ni g\mapsto gN\in G/N$  
induces a right $G$-action on $E$.

\begin{defi}
A {\bf decker} is just a continuous (right)
action $\rho$ of $G$ on $P$ by bundle 
automorphisms that
lifts the induced $G$-action on $E$
and commutes with the $\Gamma$-action of $P$.
\end{defi}

The existence of deckers can be a very 
restrictive condition  on the bundle $P\to E$. 
(See e.g. Proposition \ref{PropDeckerForG/NExIff} below.)
In fact, they  need not exist and need not to be unique
in general, but they play a central rôle in what follows, 
therefore we introduced an extra name.

\begin{defi}
A {\bf dynamical  triple} $(\rho,P,E)$ (over
the point) is a pair $(P,E)$ together  with
a decker $\rho$.
 
Equivalently, a dynamical triple is a 
$\Gamma\times G$-equivariant sequence
\begin{eqnarray*}
\begin{array}{c}
\xymatrix{
\ar@(dl,ul)[]^{\Gamma\times G} 
}
\end{array}
\hspace{-0.4cm}
\left(
\begin{array}{c}
\xymatrix{ P\ar[d]&\hspace{-0,9cm}\rotatebox {90}{$\circlearrowleft$}\ \PU(\HH) \\ 
E\ar[d]&\hspace{-1,2cm}\rotatebox {90}{$\circlearrowleft$}\ G/N \\ 
\ast}
\end{array}
\right)
\end{eqnarray*}
of trivialisable principal fibre bundles 
where $\Gamma\times G$ acts from the 
left\footnote{$G$ is commutative so any right action is a left action.}
by bundle automorphisms such that the restricted action of $1\times G$ 
on $E$ is the induced action.
\end{defi}

Let $(\rho,P,E)$ be a dynamical tripel.
As soon as  we have chosen a chart for 
the pair, we see that the decker is given by a continuous 1-cocycle 
$\mu:G\to C(G/N,\PU(\HH))$, i.e.
for $g,h\in G$  and $z\in G/N$ we have
\begin{equation}\label{EqTheCocyclePropOfMu}
\mu(g+h,z)= \mu(g,z+hN)\ \mu(h,z).
\end{equation}
As the two actions of $\Gamma$ and $G$ commute, 
we directly obtain the relation 
\begin{eqnarray}\label{EqLocStucOfDeckers}
\lambda(a,z+gN)\ \mu(g,z)=\mu(g,z+\chi(a))\ \lambda(a,z)
\end{eqnarray}
between $\mu$ and $(\chi,\lambda)$. 

Conversely, if we are given a  tuple $(\chi,\lambda,\mu)$
satisfying (\ref{EqTheCocyclePropOfMu}), (\ref{EqTheHomoPropOfChi}), (\ref
{EqTheCocyclePropertyOfLambda})
and (\ref{EqLocStucOfDeckers}), we obtain the structure
of a dynamical triple on
 $G/N\times\PU(\HH)\to G/N\to *$.
We therefore call such tuples also {\bf dynamical triples}.

Equation (\ref{EqLocStucOfDeckers}) 
can be used to give a complete answer to the existence of deckers 
in the case of $N=0$, i.e. $G=G/N$.

\begin{prop}\label{PropDeckerForG/NExIff}
Assume  $N=0$. Let $P{\to} E\stackrel{p}{\to}*$ be a pair.
Then a decker exists if and only if there is a $\Gamma$-equivariant 
bundle isomorphism $P\cong p^*P'$ for 
a $\Gamma$-equivariant $\PU(\HH)$-bundle  $P'\to *$.
\end{prop}

\begin{pf}
If $P\cong p^*P'$ for some 
$\PU(\HH)$-bundle $P'\to *$ we obtain a decker 
by the $G$-action on the first entry of the fibered product 
$p^*P'=E\times P'.$

Conversely, we fix a chart by two isomorphisms
$E\cong G/N$ and $P\cong G/N\times \PU(\HH)$
to get an induced triple $(\mu,\chi,\lambda).$
Then we define
$\lambda': \Gamma\to\PU(\HH)$ by
$\lambda'(a):=\mu(\chi(a),0)^{-1}\ \lambda(a,0)$,
for $a\in\Gamma$.
This is well-defined since $G=G/N.$
\\
Claim 1 : $\lambda'$ is a homomorphism, i.e.
a projective unitary representation.\\
Proof : Let $a,b\in\Gamma$, then
\begin{eqnarray*}
&&\lambda'(b)\ \lambda'(a)\\
&=&
\mu(\chi(b),0)^{-1}\ \lambda(b,0)\ \mu(\chi(a),0)^{-1}\ \lambda(a,0)\\
&\stackrel{(\ref{EqLocStucOfDeckers})}{=}&
\underbrace{\mu(\chi(b),0)^{-1}\  \mu(\chi(a),\chi(b)+0)^{-1}}\
\underbrace{ \lambda(b,0+\chi(a))\ \lambda(a,0)}\\
&&\qquad\qquad\qquad\qquad
\stackrel{(\ref{EqTheCocyclePropOfMu})}{=}\mu(\chi(a)+\chi(b),0)^{-1}\qquad
\stackrel{(\ref{EqTheCocyclePropertyOfLambda})}{=}\lambda(ba,0)\\
&=&\lambda'(ba).
\end{eqnarray*} 
Thus, $P':=\PU(\HH)\to*$ becomes $\Gamma$-equivariant.
\\
Claim 2 : $P\cong p^*P'$.\\
Proof :  The $\Gamma$-action on 
$p^*P'=E\times \PU(\HH)\cong G/N\times \PU(\HH)$ 
is given by the cocycle 
$$
\Gamma \ni a\mapsto \lambda'(a)\in\PU(\HH)\subset C(G/N,\PU(\HH)).
$$
We define an  isomorphism $f:P\to p^*P'$ 
by use of the chart 
by
\begin{eqnarray*}
P  &\stackrel{f}
{\longrightarrow}& p^*P'\\
\cong\downarrow&&\downarrow\cong\\
 G/N\times \PU(\HH)
&\longrightarrow&  G/N\times \PU(\HH)\\
(z,U)&\longmapsto&(z,\mu(z,0)^{-1}U).
\end{eqnarray*}
This is in fact a $\Gamma$-equivariant isomorphism
since 
from  (\ref{EqLocStucOfDeckers})
it follows for $G=G/N$ and $a\in\Gamma$
\begin{eqnarray*}
\mu(z,0)^{-1}&=&\lambda(a,0)^{-1}\ \mu(z,\chi(a))^{-1}\ \lambda(a,z)\\
&=&\lambda(a,0)^{-1}\ \mu(\chi(a),0)\ \mu(z+\chi(a),0)^{-1}\ \lambda(a,z)\\
&=&\lambda'(a)^{-1}\ \mu(z+\chi(a),0)^{-1}\ \lambda(a,z).
\end{eqnarray*}
Thus the  isomorphism $f$ commutes with the action of $\Gamma$.
\end{pf}

We introduce the notions of (iso-)morphism, 
stable isomorphism, outer conjugation and
stable outer conjugation in the same way as
we did for pairs, we just have to replace the
group $\Gamma$ and its actions by 
$\Gamma\times G$ everywhere.
We then have a similar diagram of implications 
$$
\xymatrix{
{\begin{array}{c} 
\textrm{ isomorphism} \\ 
\textrm{of dyn. triples} 
\end{array}}
\ar@{=>}[r]\ar@{=>}[d]&
{\begin{array}{c} 
\textrm{outer conjugation} \\ 
\textrm{of dyn. triples} 
\end{array}}
\ar@{=>}[d]\\
{\begin{array}{c} 
\textrm{stable isomorphism} \\ 
\textrm{of dyn. triples} 
\end{array}}
\ar@{=>}[r]&
{\begin{array}{c} 
\textrm{stable outer conjugation} \\ 
\textrm{of dyn. triples} 
\end{array}}
}
$$

We make the notion of outer conjugation
more explicit in the next lemma.
Let $(\chi,\lambda,\mu)$ be a dynamical 
triple.
The joint action of $\Gamma\times G$
is described by the cocycle 
$\nu(a,g,z):=\lambda(a,z+gN)\mu(g,z)$,
and outer conjugate triples are given by
$\nu'(a,g,z):=\nu(a,g,z)\ {\rm Ad}(n(a,g,z))$,
wherein $n:\Gamma\times G\times G/N\to\U(\HH)$ 
is continuous and satisfies 
\begin{eqnarray*}
n( yx,z)= \nu(x,z)^{-1}[n(y,z\cdot x)]\ n(x,z),
\end{eqnarray*}
for $x,y\in \Gamma\times G$, 
and $z\cdot x:=z+\chi(a)+gN$ if $x=(a,g)$.
We can characterise outer conjugate triples in terms of 
the cocycles $\lambda$ and $\mu$:

\begin{lem}\label{LemSOCODT}
Two dynamical triples  $(\chi,\lambda,\mu)$ and
$(\chi',\lambda',\mu')$ (with same underlying Hilbert space $\HH$) 
are outer conjugate if and only if (up to isomorphism)
\begin{eqnarray*}
\chi=\chi',\quad \lambda'=\lambda\ {\rm Ad}(l),\quad  \mu'=\mu\ {\rm Ad}(m),
\end{eqnarray*}
where $l:\Gamma\times G/N\to\U(\HH)$
and $m:G\times G/N\to\U(\HH)$ are continuous, unitary functions
such that the (joint) cocycle condition
\begin{eqnarray*}
l(ba,z)&=&\lambda(b,z)^{-1}[ l(b,z+\chi(a))]\ l(a,z)\\
m(g+h,z)&=&\mu(g,z)^{-1}[ m(h,z+gN)]\ m(g,z)\\
\lambda(a,z)^{-1}[m(g,z+\chi(a)]\ l(a,z) &=& \mu(g,z)^{-1}[l(a,z+gN)]\ m(g,z)
\end{eqnarray*}
is satisfied.
\end{lem}
\begin{pf}
To obtain $l$ and $m$ from $n$ one can define
\begin{eqnarray*}
l(a,z):=n(a,0,z)\quad\text{ and }\quad m(g,z):=n(1,g,z).
\end{eqnarray*}
Conversly, if $l$ and $m$ are given,  one can define
$$
n(a,g,z):= \mu(g,z)^{-1}[l(a,z+gN)]\ m(g,z).
$$
It is then easily checked that the functions defined in
this way satisfy the corresponding cocycle conditions.
\end{pf}

We illustrate our notions by two examples.

\begin{exa}
Let $(\chi,\lambda,\mu)$ be a dynamical triple, and\label{PageTheExampleOfStableEquivalence}
let $\lambda_{_{G}}:G\to\U(L^2(G))$ be the left regular
representation of $G$, 
and let $\rho:\Gamma\to \U(L^2(\Gamma))$ be 
the right regular representation of $\Gamma$.

Then the two triples
$(\chi ,\lambda,\mu)$ and $(\chi, \lambda\otimes {\rm Ad}( \Eins_{L^2(G)})\otimes {\rm Ad}(\rho),\mu\otimes{\rm Ad}(\lambda_{_{G}})\otimes {\rm Ad}(\Eins_{L^2(\Gamma)}))$ are stably outer conjugate.
\end{exa}

\begin{pf}
The argument is the same as in Example \ref{ExSOC}.
The triple $(\chi,\lambda.\mu)$ and its stabilisation  
$(\chi,\lambda\otimes\Eins\otimes\Eins,\mu\otimes\Eins\otimes\Eins)$
are stably isomorphic, and  the triples
$(\chi,\lambda\otimes\Eins,\mu\otimes\Eins)$
and $(\chi, \lambda\otimes {\rm Ad}( \Eins_{L^2(G)})\otimes {\rm Ad}(\rho),\mu\otimes{\rm Ad}(\lambda_{_{G}})\otimes {\rm Ad}(\Eins_{L^2(\Gamma)}))$ 
are outer conjugate by 
$$
l(a,z):=\Eins_{L^2(G)}\otimes \rho(a),\quad m(g,z):=
\lambda_{_{G}}(g)\otimes\Eins_{L^2(\Gamma)}
$$
which satisfy the conditions of Lemma \ref{LemSOCODT}.
\end{pf}

\begin{exa}\label{ExaOfExteriorEquivalence}
Let $(\chi,\lambda,\mu)$ be a dynamical triple
and let $(\chi,\lambda {\rm Ad}(l))$
be an outer conjugate pair to $(\chi,\lambda)$ which is 
also isomorphic to $(\chi,\lambda)$
by a morphism
\begin{eqnarray*}
\xymatrix{
\lambda {\rm Ad}(l)\ \rotatebox{90}{\rotatebox{90}{\rotatebox{90}{$\circlearrowright$}}}
\hspace{-0.9cm}& 
G/N\times\PU(\HH)\ar[d]\ar[r]^v
&G/N\times\PU(\HH)\ar[d]&
\hspace{-0.9cm} \rotatebox {90}{$\circlearrowleft$}\ \lambda,\ \mu\\
&G/N\ar[r]^=\ar[d]&G/N\ar[d]&\\
&\ast\ar[r]^=&\ast&
}
\end{eqnarray*}
which has vanishing 
class $[v]=0\in\check H^1_\Gamma(G/N,\underline{\U(1)})$
(see Appendix \ref{SecTECeCOBI12} for the definition of $[v]$).

Then the triples $(\chi,\lambda,\mu)$ and
$(\chi,\lambda{\rm Ad}(l),v^*\mu)$ are
outer conjugate.
\end{exa}

\begin{pf}
Because the class of $v$ vanishes, 
we can assume without restriction that it is 
implemented by  continuous, unitary functions
 $u:G/N\to\U(\HH)$ such that
$\lambda(a,z)^{-1}[u(z+\chi(a))]=u(z)\ l(a,z)^{-1}$.
The pulled back cocycle $v^*\mu$ is given by
\begin{eqnarray*}
(v^*\mu)(g,z)&=&{\rm Ad}(u(z+gN)^{-1}) \mu(g,z)\ {\rm Ad}(u(z))\\
&=&\mu(g,z)\ {\rm Ad}\big(\underbrace{\mu(g,z)^{-1}[u(z+gN)^{-1}] u(z)} \big).\\
&&\hspace{4cm}=: m(g,z)
\end{eqnarray*}
We have to check that $l$ and $m$ satisfy the
joint cocycle condition
of Lemma \ref{LemSOCODT}.
The first equality just involves $l$, 
and therefore holds by assumption. 
We compute the other two
\begin{eqnarray*}
&&m(g+h,z)\\
&=&\mu(g+h,z)^{-1}\big[u(z+{gN+hN})^{-1}\big]\ u(z)\\
&=&
\mu(g,z)^{-1}\Big[ \mu(h,z+gN)^{-1}\big[u(z+{gN+hN})^{-1}\big]\Big]\ u(z)\\
&=&
\mu(g,z)^{-1}\Big[ \mu(h,z+{gN})^{-1}\big[u(z+{gN+hN})^{-1}\big]\ u(z+{gN})\Big]\\
&&
\mu(g,z)^{-1}\big[u(z+{gN})^{-1}\big]\ u(z)\\
&=&\mu(g,z)^{-1}[ m(h,z+{gN}]\ m(g,z)
\end{eqnarray*}
which proves the second equality;
and
\begin{eqnarray*}
&&\lambda(a,z)^{-1}\big[m(g,z+\chi(a))\big]l(a,z)\\
&=&\lambda(a,z)^{-1}\Big[
\mu(g,z+\chi(a))^{-1}\big[
u(z+{gN+\chi(a)})^{-1}\big]\ u(z+\chi(a))\Big]l(a,z)\\
&\stackrel{(\ref{EqLocStucOfDeckers})}{=}&
\mu(g,z)^{-1}
\lambda(a,z+{gN})^{-1}\big[
u(z+{gN+\chi(a)})^{-1}\big]\ \lambda(a,z)^{-1}\big[u(z+\chi(a))\big]l(a,z)\\
&\stackrel{[\nu]=0}{=}&
\mu(g,z)^{-1}\big[
l(a,z+{gN}) u(z+gN)^{-1}\big]
 u(z)\\
&=&\mu(g,z)^{-1}\big[
l(a,z+{gN}) \big]
 m(g,z)
\end{eqnarray*}
proves the last equality. 
\end{pf}

By ${\rm Dyn}$ \label{PageOfDyn} we denote the set valued functor
that sends a finite group $\Gamma$ to the set of 
equivalence classes of stably outer conjugate dynamical triples, 
i.e.
\begin{eqnarray}
{\rm Dyn}(\Gamma):=\{\textrm{ dynamical triples } \}\big/_\textrm{st.out.conj.}.
\end{eqnarray}
In the same manner as we did for the functor ${\rm Par}$
we can define  a subfunctor of ${\rm Dyn}$
by 
\begin{eqnarray}
{\rm Dyn}(\Gamma,\chi):=\{\textrm{ dynamical triples } \textrm{ with fixed }
\chi\}
\big/_{\rm st.out.conj.}.
\end{eqnarray}
This yields a decomposition 
${\rm Dyn}(\Gamma)= \coprod_{\chi} {\rm Dyn}(\Gamma,\chi)$.
For fixed $\chi$, each  ${\rm Dyn}(\Gamma,\chi)$ has 
the structure of an abelian group given by tensor product 
and complex conjugation of projective unitary bundles, and
by Lemma \ref{LemSOCODT}  there
are  well-defined forgetful maps
\begin{eqnarray*}
{\rm Dyn}(\Gamma,\chi)\to{\rm Par}(\Gamma,\chi), \quad 
{\rm Dyn}(\Gamma)\to{\rm Par}(\Gamma)
\end{eqnarray*}
which extend to natural transformations
\begin{eqnarray*}
{\rm Dyn}\to{\rm Par}.
\end{eqnarray*}

\section{Dualisable Dynamical Triples}

\noindent
Fix a homomorphism
\begin{eqnarray*}
\chi:\Gamma\to G/N.
\end{eqnarray*}
We consider $L_\chi:=L^\infty(G/N,\U(1))$
as a topological $\Gamma$-$G$-bimodule
in a natural fashion: The action is given by shifting the 
argument of a function by $-\chi(a)$ or $-gN$ for
$a\in\Gamma$ or $g\in G$.
The topology on $L_\chi$ is the 
strong (or weak) operator topology, i.e. we consider
$L_\chi$ as a a set of unitary multiplication operators
on the Hilbert space $L^2(G/N,\HH)$. 

Let 
$
B^{k,l}(\Gamma,G,L_\chi)
$ 
be the group of all Borel functions
$\Gamma^k\times G^l\to L_\chi$,
then we find a double complex
$$
\xymatrix{
\vdots&\vdots&\vdots& \\
 B^{2,0}(\Gamma,G,L_\chi)\ar[r]^{d}\ar[u]^{\delta_\chi}&
 B^{2,1}(\Gamma,G,L_\chi)\ar[r]^{d}\ar[u]^{\delta_\chi}&
 B^{2,2}(\Gamma,G,L_\chi)\ar[r]^{\quad d}\ar[u]^{\delta_\chi}&\cdots\\
 B^{1,0}(\Gamma,G,L_\chi)\ar[r]^{d}\ar[u]^{\delta_\chi}&
 B^{1,1}(\Gamma,G,L_\chi)\ar[r]^{d}\ar[u]^{\delta_\chi}&
 B^{1,2}(\Gamma,G,L_\chi)\ar[r]^{\quad d}\ar[u]^{\delta_\chi}&\cdots\\
 B^{0,0}(\Gamma,G,L_\chi)\ar[r]^{d}\ar[u]^{\delta_\chi}&
 B^{0,1}(\Gamma,G,L_\chi)\ar[r]^{d}\ar[u]^{\delta_\chi}&
 B^{0,2}(\Gamma,G,L_\chi)\ar[r]^{\quad d}\ar[u]^{\delta_\chi}&\cdots,
}
$$
wherein $d$ for the group $G$ and $\delta_\chi$ 
for $\Gamma$ are  the  group cohomological  
boundary operators.  
We denote the resulting total complex
by $\big( B^\bullet(\Gamma,G,L_\chi),\partial_\chi\big)$, i.e.
$$
B^p(\Gamma,G,L_\chi):=\bigoplus_{p=k+l}  B^{k,l}(\Gamma,G,L_\chi)
$$
and (if we write the group law in $L$ additive for a moment)
$$
\partial_\chi|_{B^{k,l}(\Gamma,G,L_\chi)}:=d+(-1)^l\delta_\chi.
$$ 
We denote by $H^\bullet_{\rm Bor}(\Gamma,G,L_\chi)$ the corresponding  
cohomology groups. 
\\

The algebraic structure of a dynamical triple gives rise to 
a 2-cohomology class as  we explain now.
Let $(\chi,\lambda,\mu)$
be a dynamical triple.
One should  realise at this point that when we
suppress the non-commutativity of $\PU(\HH)$ for a moment 
the  cocycle conditions for  $\lambda$ and $\mu$ together 
with (\ref{EqLocStucOfDeckers}) are equivalent 
to $\partial_\chi(\lambda,\mu)=\Eins$. 
Now, we lift the cocycles to unitary valued, Borel functions
which will define a 2-cocycle. In detail:
Choose for both, $\mu$ and $\lambda$, unitary lifts to Borel functions 
$\overline\mu:G\to L^\infty(G/N,\U(\HH))$ 
and $\overline \lambda:\Gamma\to L^\infty(G/N,\U(\HH))$
such that ${\rm Ad}(\overline\mu(g,z))=\mu(g,z)$ and
${\rm Ad}(\overline\lambda(a,z))=\lambda(a,z)$.
Then we define $(\psi,\phi,\omega)$ by
\begin{eqnarray}
\overline\lambda(b,z+\chi(a))\ 
\overline\lambda(a,z)& =&\overline\lambda(ba,z)\ \psi(a,b,z),\nonumber\\
\overline\lambda(a,z+gN)\ \overline\mu(g,z) 
&=& \overline\mu(g,z+\chi(a))\ \overline\lambda(a,z)\ \phi(a,g,z),\nonumber\\
\overline\mu(h,z+gN)\ \overline\mu(g,z) &=& \overline\mu(g+h,z)\ \omega(g,h,z).
	\label{EqTheDefOfPsPhOm}
\end{eqnarray}
Due to (\ref{EqTheCocyclePropertyOfLambda}), (\ref{EqTheCocyclePropOfMu}) and (\ref
{EqLocStucOfDeckers}),
these three functions $\psi,\phi$ and $\omega$
are $\U(1)$-valued and satisfy the algebraic relations
(written multiplicative)\label{PageOfDefiOfTotalCocycle}
\begin{eqnarray}
\delta_\chi\psi&=&1,\nonumber\\
\delta_\chi\phi&=& d\psi,\nonumber\\
d\phi&=&\delta_\chi\omega^{-1},\nonumber\\
d\omega&=&1
	\label{EqTheAlgRelOfPsPhOm}
\end{eqnarray}
which is (again in additive notation) equivalent to  
\begin{equation}
\partial_\chi(\psi,\phi,\omega)=0\in B^3(\Gamma,G,L_\chi),
\end{equation}
i.e. $(\psi,\phi,\omega)$ is a 2-cocycle. 
Of course, one can verify this by direct computation,
but indeed it is implicitly clear, because, informally\footnote{
i.e. up to the non-comutativity of $\U(\HH)$},
we have defined
$(\psi,\phi,\omega):=\partial_\chi
(\overline\lambda,\overline\mu)$.

\begin{prop}\label{PropWouldLikeToClassifyTriples}
The assignment 
$(\chi,\lambda,\mu)\mapsto (\psi,\phi,\omega)$
constructed above defines a homomorphism of groups
$$
{\rm Dyn}(\Gamma,\chi)\to H^2_{\rm Bor}(\Gamma,G,L_\chi).
$$
\end{prop}
\begin{pf}
We must check that the defined total cohomology class
is independent of all choices.
This is simple to verify for the choice of 
the chart, and the choice of the lifts of 
the transition functions and cocycles.
It is also clear that stably isomorphic triples define the same
total cohomology class.
By Lemma \ref{LemSOCODT}
we can directly calculate 
the cocycle $(\psi',\phi',\omega')$ 
for an outer conjugate triple. This calculation 
is straight forward, and, in fact, for
$\overline\mu':=\overline\mu m$ and
$\overline\lambda':=\overline\lambda l$
we find $(\psi',\phi',\omega')=(\psi,\phi,\omega)$.
\end{pf}

Note that we do not claim anything about the
injectivity or surjectivity of the map in Proposition
\ref{PropWouldLikeToClassifyTriples}.
We do better in the next section, where we consider 
dualisable triples only. For the definition of these 
consider the forgetful map
\begin{eqnarray}
H^2_{\rm Bor}(\Gamma,G,L_\chi)&\to&
H^2_{\rm Bor}(G,L^\infty(G/N,\U(1)))\cong H^2(N,\U(1)).\\
\ [(\psi,\phi,\omega)]&\mapsto& [\omega]\nonumber
\end{eqnarray}
The corresponding class
in $H^2(N,\U(1))$ is called the 
{\bf Mackey obstruction} of the dynamical triple.
We will have our focus on those triples which 
have a vanishing Mackey obstruction.

\begin{defi}
A dynamical triple is called {\bf dualisable} if its 
Mackey obstruction vanishes.
\end{defi}

There is the contravariant set valued sub-functor 
of ${\rm Dyn}$ which 
sends  a finite group  to the 
stable outer conjugation classes of 
dualisable dynamical triples. It is denoted by
${\rm Dyn}^\dag$, so
\begin{eqnarray}
{\rm Dyn}^\dag(\Gamma)\subset{\rm Dyn}(\Gamma), \quad 
{\rm Dyn}^\dag(\Gamma,\chi)
\subset{\rm Dyn}(\Gamma,\chi).
\end{eqnarray}

\section{The Classification of Dualisable Dynamical Triples}
\label{SecTCODDT}

\noindent
We start similar to the previous section, but we 
stick to the continuous setting rather than to the
Borel setting.
For  a homomorphism
\begin{eqnarray*}
\chi:\Gamma\to G/N
\end{eqnarray*}
we consider 
the topological $\Gamma$-$G$-bimodule
$M_\chi:=C(G/N,\U(1))$,
where $(a,g)\in\Gamma\times G$
acts by shift with $-\chi(a)-gN$ in the 
arguments of the functions.
The topology on $M_\chi$ is the 
compact-open  topology.
Let 
$
C^{k,l}(\Gamma,G,M_\chi)
$ 
be the group of all continuous  functions
$\Gamma^k\times G^l\to M_\chi$.
Again we have a double complex
$$
\xymatrix{
\vdots&\vdots&\vdots& \\
 C^{2,0}(\Gamma,G,M_\chi)\ar[r]^{d}\ar[u]^{\delta_\chi}&
 C^{2,1}(\Gamma,G,M_\chi)\ar[r]^{d}\ar[u]^{\delta_\chi}&
 C^{2,2}(\Gamma,G,M_\chi)\ar[r]^{\quad d}\ar[u]^{\delta_\chi}&\cdots\\
 C^{1,0}(\Gamma,G,M_\chi)\ar[r]^{d}\ar[u]^{\delta_\chi}&
 C^{1,1}(\Gamma,G,M_\chi)\ar[r]^{d}\ar[u]^{\delta_\chi}&
 C^{1,2}(\Gamma,G,M_\chi)\ar[r]^{\quad d}\ar[u]^{\delta_\chi}&\cdots\\
 C^{0,0}(\Gamma,G,M_\chi)\ar[r]^{d}\ar[u]^{\delta_\chi}&
 C^{0,1}(\Gamma,G,M_\chi)\ar[r]^{d}\ar[u]^{\delta_\chi}&
 C^{0,2}(\Gamma,G,M_\chi)\ar[r]^{\quad d}\ar[u]^{\delta_\chi}&\cdots,
}
$$
We denote the resulting total cohomology groups
by $H^\bullet_{\rm cont}(\Gamma,G,M_\chi)$. 
For the classification of dualisable dynamical triples 
we have to consider the kernel of the forgetful map
\begin{eqnarray*}
H^2_{\rm cont}(\Gamma,G,M_\chi)\to H^2_{\rm cont}(G,C(G/N,\U(1)))
\end{eqnarray*}
which we denote  by
$H^{2,\dag}_{\rm cont}(\Gamma,G,M_\chi)$.

Now, we can state the classification theorem of 
dualisable dynamical triples.

\begin{thm}\label{ThmTheCFODDT}
There is a natural isomorphism of groups
\begin{eqnarray*}
{\rm Dyn}^\dag(\Gamma,\chi)\stackrel{\cong}{\to}
 H^{2,\dag}_{\rm cont}(\Gamma,G,M_\chi).
\end{eqnarray*}
\end{thm}

\noindent
The proof of Theorem \ref{ThmTheCFODDT} is the content
of the remainder of this section.
\\

Consider a dualisable dynamical triple $(\chi,\lambda,\mu)$.
As its Mackey obstruction vanishes, we can 
find a unitary Borel lift  $\overline \mu:G\to L^\infty(G/N,\U(\HH))$
of $\mu$ which satisfies again the cocycle condition
\begin{eqnarray}\label{EqTheUniCCOfMu}
\overline\mu(g+h,z)=\overline\mu(h,z+gN)\ \overline\mu(g,z).
\end{eqnarray}
The function $\overline\mu(g,\_)\in L^\infty(G/N,\U(\HH))$ is 
a multiplication operator which acts
by $(\overline\mu(g,\_)F)(z):=\overline\mu(g,z)F(z)$ for 
$F\in L^2(G/N)\otimes\HH$.
We also may consider 
$\overline\mu(\ .\ , \_)$ or 
$\overline \mu(\ .\ , z)$ as multiplication
operators on $L^2(G)\otimes L^2(G/N)\otimes\HH$
or $L^2(G)\otimes \HH$ respectively.

As a matter of fact, the cocycle condition (\ref{EqTheUniCCOfMu})
implies some useful continuity properties of the functions mentioned.
This is the content of the next technical lemmata.

\begin{lem}\label{LemTechLFMu11}
The map 
\begin{eqnarray*}
G&\to& \U(L^2(G/N)\otimes\HH)\\
g&\mapsto& \overline\mu(g,\_)
\end{eqnarray*}
is continous.
\end{lem}

\begin{pf}
We consider 
$
\Eins\otimes\overline\mu(g,\_)\in \U(L^2(G)\otimes L^2(G/N)\otimes\HH).
$
Using the cocycle condition we have
\begin{eqnarray*}
\Eins\otimes\overline\mu(g,\_)&=&\overline\mu(\ .\ ,\_+gN)^{-1}\ \overline\mu(g+\ .\ ,\_)\\
&=&\lambda_{_{G/N}}(-gN) \overline\mu(\ .\ ,\_)^{-1}\lambda_{_{G/N}}(gN)
\ \lambda_{_{G}}(-g)\overline\mu(\ .\ ,\_)\lambda_{_{G}}(g),
\end{eqnarray*}
where $\lambda_{_{G/N}}$  and $\lambda_{_{G}}$ are the left regular 
representations of the groups $G/N$ and $G$ respectively.
As the left regular representations are (strongly) continuous 
the assertion follows.
\end{pf}

In the next lemma $\sigma :G/N\to G$ is a Borel 
section of the quotient map $G\to G/N$.

\begin{lem}\label{LemTechLFMu22}
The maps 
\begin{eqnarray*}
 G/N&\to& \PU(L^2(G)\otimes\HH)\\
  z&\mapsto& {\rm Ad}(\overline\mu(\ .\ ,z))
\end{eqnarray*}
and
\begin{eqnarray*}
 G/N&\to& \PU(L^2(G/N)\otimes\HH)\\
  z&\mapsto& {\rm Ad}(\overline\mu(\sigma(\_) ,z))
\end{eqnarray*}
are continuous.
\end{lem}

\begin{pf}
Let $z_\alpha\to z$ be a converging net. Let $x_\alpha:=z_\alpha-z$
and choose $g_\alpha\to 0\in G$ such that\footnote{
Such $g_\alpha$
exist -- take a local section of the quotient $G\to G/N$. 
}
$g_\alpha N=x_\alpha$.
Then 
\begin{eqnarray*}
{\rm Ad}(\overline\mu(\ .\ ,z_\alpha) )
&=&{\rm Ad}(\overline\mu(\ .\ ,z+x_\alpha) )\\
&=&{\rm Ad}(\overline\mu(\ .\ +g_\alpha,z) )
{\rm Ad}(\overline\mu(g_\alpha,z)^{-1} )\\
&=&{\rm Ad}(\overline\mu(\ .\ +g_\alpha,z) )
\ \mu(g_\alpha,z)^{-1} \\
&=&{\rm Ad}(\lambda_{_{G}}(-g_\alpha)\overline\mu(\ .\ ,z)\lambda_{_{G}}(g_\alpha) )
\ \mu(g_\alpha,z)^{-1}\\ 
&\to&{\rm Ad}(\overline\mu(\ .\ ,z) )
\ \mu(0,z)^{-1}. 
\end{eqnarray*}
As $\mu(0,z)=\Eins$, the continuity of the first map follows.

To prove the continuity of the second map 
we show first that 
\begin{eqnarray*}
\sigma^*:L^\infty(G,\U(\HH))&\to L^\infty(G/N,\U(\HH))\\
\nu(\ .\ )\mapsto \nu(\sigma(\_))
\end{eqnarray*}
is  continuous.
Indeed, let $f\in L^2(G/N)$ and let
$\varphi$ be the characteristic function of $\sigma(G/N)\subset G$,
so $g\mapsto f(gN)\varphi(g)$ is a function in $L^2(G)$.
Let $\nu_n(\ .\ )\to \nu(\ .\ )\in L^\infty(G,\U(\HH))$ be a converging sequence.
Then
\begin{eqnarray*}
\|\nu_n(\sigma(\_))f(\_)-\nu(\sigma(\_))f(\_)\|^2
&=&\int_{G/N}|\nu_n(\sigma(z))f(z)-\nu(\sigma(z))f(z)|^2\ dz\\
&=&\int_{G}|\nu_n(g)f(gN)\varphi(g)-\nu(g)f(gN)\varphi(g)|^2\ dg\\
&\to& 0,\ {\rm for\ } n\to \infty.
\end{eqnarray*}
Because ${\rm Ad}\circ\sigma^*$ is constant along 
the orbits of $\U(1)$ the dotted arrow
in
\begin{eqnarray*}
\xymatrix{
L^\infty(G,\U(\HH))\ar[d]^{\rm Ad}\ar[r]^{\sigma^*}&L^\infty(G/N,\U(\HH))\ar[d]^{\rm Ad}\\
P L^\infty(G,\U(\HH))\ar@{.>}[r]^{\sigma^*}&P L^\infty(G/N,\U(\HH))
}
\end{eqnarray*}
is well-defined and continuous by the universal property 
of the quotient map.

This proves the Lemma.
\end{pf}

With these two lemmata at hand we can proof 
the following proposition.

\begin{prop}\label{PropForDefOfClMap}
For each dualisable dynamical triple $(\chi,\lambda,\mu)$
with underlying Hilbert space $\HH$
there exists a stably outer conjugate triple
$(\chi,\lambda',\mu')$ with underlying Hilbert space
$\HH':= L^2(G/N)\otimes \HH$ such that
the cocycles $\lambda'$ and $\mu'$ permit 
continuous lifts
\begin{eqnarray*}
\begin{array}{c}
\xymatrix{
&\U(\HH')\ar[d]\\
\Gamma\times G/N\ar[r]^{\lambda'}\ar[ru]^{\overline{\lambda'}}&\PU(\HH')
}
\end{array}
\quad and
\quad
\begin{array}{c}
\xymatrix{
&\U(\HH')\ar[d]\\
G\times  G/N\ar[r]^{\mu'}\ar[ru]^{\overline{\mu'}}&\PU(\HH')
}
\end{array},
\end{eqnarray*}
where $\overline{\mu'}$ still satisfies the cocycle condition (\ref{EqTheUniCCOfMu}). 
\end{prop}

\begin{pf}
The triple  $(\chi,\Eins\otimes\lambda,\Eins\otimes\mu)$
with underlying Hilbert space $L^2(G/N)\otimes\HH$ is 
stably isomorphic to $(\chi,\lambda,\mu)$.
Let 
\begin{eqnarray*}
\xymatrix{
& 
G/N\times\PU(\HH')\ar[d]\ar[r]^\theta
&G/N\times\PU(\HH')\ar[d]&
\hspace{-0.9cm} \rotatebox {90}{$\circlearrowleft$}\ \Eins\otimes\lambda,\ \Eins\otimes\mu\\
&G/N\ar[r]^=\ar[d]&G/N\ar[d]&\\
&\ast\ar[r]^=&\ast&
}
\end{eqnarray*}
be the bundle isomorphism which is 
given by  $\theta(z,U):=(z, {\rm Ad}(\overline\mu(\sigma(\_) ,z))^{-1} U)$.
By Lemma \ref{LemTechLFMu22} this is well-defined.
We define 
\begin{eqnarray*}
(\chi,\lambda',\mu'):=\theta^*(\chi,\Eins\otimes\lambda,\Eins\otimes\mu),
\end{eqnarray*}
i.e. 
\begin{eqnarray*}
\mu'(g,z)&:=& {\rm Ad}\big(\overline\mu(\sigma(\_) ,z+gN)\big)\ \mu(g,z)\ 
{\rm Ad}\big(\overline\mu(\sigma(\_) ,z)\big)^{-1}\\
&=&{\rm Ad}\big( \overline\mu(\sigma(\_ ),z+gN)\ 
\overline\mu(g,z)\ \overline\mu(\sigma(\_) ,z)^{-1} \big)\\
&=&{\rm Ad}\big( \overline\mu(g+\sigma(\_) ,z+gN)\ \ \overline\mu(\sigma(\_) ,z)^{-1} \big)\\
&=&{\rm Ad}\big( \overline\mu(g ,z+ \_  )  \big)\\
&=&{\rm Ad}\big( \underbrace{\lambda_{_{G/N}}(-z)\overline\mu(g ,\_)
\lambda_{_{G/N}}(z) } \big)\\
&&\hspace{2,5cm} =:\overline{\mu'}(g,z).
\end{eqnarray*}
By Lemma \ref{LemTechLFMu11} $\overline{\mu'}$ is continous,
it is also clear that it satisfies (\ref{EqTheUniCCOfMu}).

It remains to show that $\lambda'$ also posesses a 
continous lift.
But this follows from equation (\ref{EqLocStucOfDeckers})
which implies for $z=0$ that
\begin{eqnarray*}
\lambda'(a,gN)=\mu'(g,\chi(a))\ \lambda(a,0) \ \mu'(g,0)^{-1};
\end{eqnarray*}
as $\Gamma$ is discrete,  
the right-hand side has a continuous lift, hence 
the left-hand side has.
\end{pf}

Proposition \ref{PropForDefOfClMap} enables us 
to define the classification map of Theorem \ref{ThmTheCFODDT}.
In fact, we define $(\psi,\phi,1)$ as in equation (\ref{EqTheDefOfPsPhOm}),
but with with the additional assumption, that 
$\overline\lambda$ and $\overline\mu$ are continuous.
This defines 
\begin{eqnarray}\label{EqGRROFOFOO}
{\rm Dyn}^\dag(\Gamma,\chi)\to H^{2,\dag}_{\rm cont}(\Gamma,G,M_\chi).
\end{eqnarray}
To check that our definition is well-defined requires 
the same arguments as in Proposition \ref{PropWouldLikeToClassifyTriples}.

To show that (\ref{EqGRROFOFOO}) is an isomorphism we
construct a map in opposite direction.
To do so, let $(\psi,\phi,1)$ represent an element of 
$ H^{2,\dag}_{\rm cont}(\Gamma,G,M_\chi)$.
We regard  $\psi(a, : ,z)$ and 
$\phi( : ,g,z)$ as unitary multiplication operators on
$L^2(\Gamma)$.
Let
$\rho:\Gamma\to \U(L^2(\Gamma))$ be the right regular
representation, then 
we define
\begin{eqnarray*}
\lambda_{\psi}(a,z)&:=&{\rm Ad}(\psi(a, : , z) \rho(a))\\
\mu_\phi(g,z)&:=&{\rm Ad}(\phi( : ,g,z)).
\end{eqnarray*}
The cocycle condition $\partial_\chi(\psi,\phi,1)=1$
is equivalent to 
\begin{eqnarray*}
\psi(b,:,z+\chi(a)\rho(b)\ \psi(a,:,z)\rho(a)&=& \psi(ba,:,z)\rho(ba)\ \underbrace{\psi(a,b,z)},\\
&&\hspace{3cm} \in\U(1)\\
\psi(a,:,z+gN)\rho(a)\ \phi(:,g,z)&=&
\phi(:,g,z+\chi(a))\ \psi(a,:,z)\rho(a) \underbrace{\phi(a,g,z)},\\
&&\hspace{5cm} \in\U(1)\\
\phi(:,h,z+gN)\ \phi(:,g,z)&=&\phi(:,g+h,z).
\end{eqnarray*}
By taking ${\rm Ad}:\U(L^2(\Gamma))\to\PU(L^2(\Gamma))$ on both sides 
of the equalities above, it follows 
that $(\chi,\lambda_\psi,\mu_\phi)$ is a dynamical triple 
which is mapped to $[\psi,\phi,1]$ under (\ref{EqGRROFOFOO}).
In this way we obtain a map
\begin{eqnarray}\label{EqGRROFOFOOInv}
 H^{2,\dag}_{\rm cont}(\Gamma,G,M_\chi)\to
{\rm Dyn}^\dag(\Gamma,\chi)
\end{eqnarray}
such that the composition with (\ref{EqGRROFOFOO})
is the identity on $H^{2,\dag}_{\rm cont}(\Gamma,G,M_\chi)$.

To conclude that (\ref{EqGRROFOFOO}) is an isomorphism 
we show that the composition of (\ref{EqGRROFOFOO}) with 
(\ref{EqGRROFOFOOInv}) is the identity on ${\rm Dyn}^\dag(\Gamma,\chi)$.
So take a triple $(\chi,\lambda,\mu)$ (with underlying Hilbert space $\HH$) 
which permits continuous lifts
of $\lambda$ and $\mu$, then define $(\psi,\phi,1)$ 
as in (\ref{EqTheDefOfPsPhOm})
and $\lambda_{\psi},\mu_\phi$ (with underlying Hilbert space $L^2(\Gamma)$) as above.
We claim that  $(\chi,\lambda,\mu)$ and
$(\chi,\lambda_\psi,\mu_\phi)$ are stably
outer conjugate. In fact, 
$(\chi,\lambda,\mu)$ and $(\chi,\lambda{\rm Ad}(\rho),\mu\otimes\Eins)$
(with underlying Hilbert space $\HH\otimes L^2(\Gamma)$) are 
stably outer conjugate. 
Now, note that (\ref{EqTheDefOfPsPhOm}) implies
\begin{eqnarray*}
\lambda(a,z){\rm Ad}(\rho(a))&=&
 {\rm Ad}\Big(\overline\lambda(:,z+\chi(a))^{-1}\  
\overline\lambda(:a,z)\ \psi(a,:,z)\ \rho(a)\Big)\\
&=& {\rm Ad}\big(\overline\lambda(:,z+\chi(a))\big)^{-1}\ \lambda_\psi(a,z)\ 
{\rm Ad}\big(\overline \lambda(:,z)\big)
\end{eqnarray*}
and
\begin{eqnarray*}
\mu(g,z)&=&
{\rm Ad}\Big(\overline\lambda(:,z+gN)^{-1}\  
\overline\mu(g,z+\chi(:))\ \overline\lambda(:a,z)\ \phi(:,g,z)\Big)\\
&=&
{\rm Ad}\big(\overline\lambda(:,z+gN)\big)^{-1}\  
\mu_\phi(g,z)\ {\rm Ad}\big( \underbrace{\overline\mu(g,z+\chi(:))}\big)    {\rm Ad}\big(\overline\lambda(:a,z)\big).\\
&&\hspace{6,1cm}=:m(g,z)
\end{eqnarray*}
Thus,  we see that there is an isomorphism
of dynamical triples 
\begin{eqnarray*}
\xymatrix{
\Eins\otimes\lambda_\psi,\mu_\phi{\rm Ad}(m)\ 
\rotatebox {90}{\rotatebox {90}{\rotatebox {90}{$\circlearrowright$}}}
\hspace{-0,9cm}
& 
G/N\times\PU(\HH_\Gamma)\ar[d]\ar[r]^\theta
&G/N\times\PU(\HH_\Gamma)\ar[d]&
\hspace{-0.9cm} \rotatebox {90}{$\circlearrowleft$}\ \lambda{\rm Ad}(\rho),\ \mu\otimes\Eins\\
&G/N\ar[r]^=\ar[d]&G/N\ar[d]&\\
&\ast\ar[r]^=&\ast,&
}
\end{eqnarray*}
where $\HH_\Gamma:=\HH\otimes L^2(\Gamma)$ and
  $\theta(z,U):=(z, {\rm Ad}(\overline\lambda( :,z)) U)$.
But the triple 
$(\chi,\Eins\otimes\lambda_\psi,\mu_\phi{\rm Ad}(m))$
is stably outer conjugate to 
$(\chi,\lambda_\psi,\mu_\phi)$, because 
$m$ satisfies the cocycle condition of 
Lemma \ref{LemSOCODT}.
This shows that the composition  of (\ref{EqGRROFOFOO}) with 
(\ref{EqGRROFOFOOInv}) is the identity on ${\rm Dyn}^\dag(\Gamma,\chi)$.
\\

We have just proven Theorem \ref{ThmTheCFODDT}.

\section{Dual Pairs and Triples}
\label{SecDualPairsAndTriples}

\noindent
Of course, in the whole discussion we 
made so far the group $G$ and its subgroup $N$ 
 can be replaced by 
 the dual group  $\widehat G:={\rm Hom}(G,\U(1))$ and 
the anihilator  $N^\perp:=\{\hat g|\hat g|_N=1\} \subset \widehat G$ of $N$
everywhere.
This is meaningful as 
$N^\perp\cong\widehat{G/N}$ is discrete, 
and $\widehat G/N^\perp\cong \widehat N$ is compact.
\\

Let $\Gamma$ be a finite group and $\HH$ be a separable 
Hilbert space.

\begin{defi}
\begin{enumerate}
\item
A {\bf dual  pair} $(\widehat P,\widehat E)$ with underlying 
Hilbert space $\HH$ is 
a $\Gamma$-equivariant sequence $\widehat P\to\widehat E\to *$, wherein 
$\widehat E\to *$ and  
$\widehat P\to \widehat E$ are trivialisable principal 
fibre bundles with structure groups $\widehat G/N^\perp$,
 $\PU(\HH)$ respectively,
and $\Gamma$ acts by bundle automorphisms.
 
\item \label{PageOfAllTheDuals}
A {\bf dual decker } $\hat\rho$ is an action 
by bundle automorphisms of 
$\widehat G$ on $\widehat P$  
which lifts the induced 
$\widehat G$-action on $\widehat E$ and commutes with
the given $\Gamma$-action of the pair.

\item
A {\bf dual dynamical triple} $(\hat\rho,\widehat P,\widehat E)$
is a dual pair $(\widehat P,\widehat E)$ together  with
a dual decker  $\hat\rho$

\item 
A {\bf dualisable dual dynamical triple}
is a dual dynamical triple whose Mackey 
obstruction vanishes.
\end{enumerate}
\end{defi}

It is clear how to define   
$\widehat{\rm Par}(\Gamma), 
\widehat{\rm Dyn}(\Gamma)$ and
$\widehat{\rm Dyn}^\dag(\Gamma)$,
and for a homomorphism
\begin{eqnarray*}
\hat\chi:\Gamma\to \widehat G/N^\perp
\end{eqnarray*}
we also define
$\widehat{\rm Par}(\Gamma,\hat\chi),$
\label{PageOfDualFunctors}
$\widehat{\rm Dyn}(\Gamma,\hat\chi)$
and
$\widehat{\rm Dyn}^\dag(\Gamma,\hat\chi)$
as before.

All statements 
we have achieved so far
translate to the dual setting
in the obvious way.
In particular we stress the classification 
theorem for dualisable dual dynamical triples
(cp. Theorem \ref{ThmTheCFODDT}):
\begin{thm}\label{ThmTheCFODDTDualer}
There is a natural isomorphism of groups
\begin{eqnarray*}
\widehat{\rm Dyn}^\dag(\Gamma,\hat\chi)\stackrel{\cong}{\to} 
H^{2,\dag}_{\rm cont}(\Gamma,\widehat G,M_{\hat\chi}).
\end{eqnarray*}
\end{thm}
Here $M_{\hat\chi}:=C(\widehat G/N^\perp,\U(1))$
is the dual counterpart of $M_{\chi}$ introduced in 
section \ref{SecTCODDT}, i.e. 
it has the $\Gamma$-$\widehat G$-bimodule structure
given by shift with $-\hat\chi(a)-\hat gN^\perp$ in the 
arguments of the functions, for $a\in\Gamma, \hat g\in \widehat G$.

\section{The Duality Theory of Dynamical Triples}
\label{SecDualityTheoryODT}

\noindent
In this section we show that the two functors 
${\rm Dyn}^\dag$ and $\widehat{\rm Dyn}^\dag$
are naturally isomorphic functors. To do so
we construct an assignment 
$(\chi,\psi,\phi)\mapsto(\hat\chi,\hat\psi,\hat\phi)$
and then
use the classification theorems of dynamical
and dual dynamical triples,
Theorem \ref{ThmTheCFODDT} and 
Theorem \ref{ThmTheCFODDTDualer}.

Let $(\chi,\lambda,\mu)$ be a dualisable dynamical triple
represented by $(\psi,\phi,1)$.
As $d\phi=1$ we have
\begin{eqnarray*}
\phi(a,g,{hN})\ \phi(a,h,0)&=&\phi(a,g+h,0)\\
&=&\phi(a,h+g,0)\\
&=&\phi(a,h,{gN})\ \phi(a,g,0)
\end{eqnarray*}
which implies for $g=n\in N$ 
that
$\phi(a,n,{hN})=\phi(a,n,0)$ holds, hence
$z\mapsto\phi(a,n,z)$ is a constant.
Further,  $N\ni n\mapsto\phi(a,n,0)$ is a continuous homomorphism. Thus
there exists a map
\label{PageThePageOfHatChi} 
$\hat \chi :\Gamma\to \widehat G/N^\perp (\cong\widehat N)$
such that
\begin{equation}\label{EqTheDualTorusCoycle123}
\langle \hat \chi(a),n\rangle=\phi(a,n,0)^{-1} 
\end{equation}

\begin{prop}
 $\hat \chi :\Gamma \to \widehat G/N^\perp$  is a
homomorphism of groups.
\end{prop}

\begin{pf}
We have $\delta_\chi\phi=d\psi$, explicitely this reads
for $g\in G, a,b\in\Gamma$
$$
\phi(b,g,{\chi(a)})
\phi(ba,g,0)^{-1}
\phi(a,g,0)
=\psi(a,b,{gN})
\psi(a,b,0)^{-1},
$$
and for $g=n\in N$ the right-hand side vanishes.
\end{pf}

We choose a lift $\overline\chi$ of $\chi$ in
\begin{eqnarray*}
\xymatrix{
&G\ar[d]\\
\Gamma\ar[ru]^{\overline\chi}\ar[r]^\chi&G/N
},
\end{eqnarray*}
and then we define
\begin{eqnarray}\label{EqTheDeOTDPsi}
\hat\psi(a,b,\hat z) &:=&
\psi(a,b,0)\ \phi(b,\overline\chi(a),0)^{-1}\nonumber\\
&&
\langle \hat\chi(ba)+\hat z, \overline\chi(ba)-\overline\chi(a)-\overline\chi(b)\rangle
\end{eqnarray}
and
\begin{eqnarray}\label{EqTDOTDPhi}
\hat\phi(a,\hat g,\hat z):=
\langle \hat g,\overline\chi(a)\rangle^{-1},
\end{eqnarray}
wherein the scalar product in the definition of $\hat\psi$ is
the paring of $\widehat G/N^\perp$ with $N$, and paring
in the definition of $\hat\phi$ is between $\widehat G$ and
$G$.

\begin{prop}\label{PropTheDualTCCC}
(In multiplicative notation)  the cocycle condition 
$$
\partial_{\hat\chi}(\hat\psi,\hat\phi,1)=1
$$
is satisfied, and the class
\begin{eqnarray*}
[\hat\psi,\hat\phi,1]\in H^{2,\dag}_{\rm cont}(\Gamma,\widehat G, M_{\hat\chi})
\end{eqnarray*}
is independent of the choices of $\overline\chi$
and the representative $(\hat\psi,\hat\phi,1)$.
\end{prop}

\begin{pf}
The proof is a straight forward calculation.

Firstly, the equality $\hat d\hat\phi=1$ is satisfied
as $\widehat G\ni g\mapsto\hat\phi(a,\hat g,\hat z)$
is a homomorphism.

Secondly, we have $\hat d\hat\psi=\delta_{\hat\chi}\hat\phi$ as
\begin{eqnarray*}
(\hat d\hat\psi)(a,b,\hat g, \hat z)&=&
\langle gN,\overline\chi(ba)-\overline\chi(a)-\overline\chi(b)\rangle\\
&=&
(\delta_{\hat\chi}\hat\phi)(a,b,\hat g,\hat z).
\end{eqnarray*}

Thirdly, we show that  $\delta_{\hat\chi}\hat\psi=1$ which is the most 
lengthy equation to check. The proof requires the cocycle condition 
$\delta_{\chi}\psi=1, d\psi=\delta_{\chi}\phi, d\phi=1$.
In fact, we have
\begin{eqnarray*}
&&(\delta_{\hat\chi}\hat\psi)(a,b,c,z)\\
&=&
\hat\psi(b,c,\hat z+\hat\chi(a))\ \psi(ba,c,\hat z)^{-1}\ 
\hat\psi(a,cb,\hat z)\ \hat\psi(a,b,\hat z)^{-1}\\
&=&
\psi(b,c,0)\ \psi(ba,c,0)^{-1}\ \psi(a,cb,0)\ \psi(a,b,0)^{-1}\\
&&
\phi(c,\overline\chi(b),0)^{-1}\
\phi(c,\overline\chi(ba),0)\
\phi(cb,\overline\chi(a),0)^{-1}\
 \phi(b,\overline\chi(a),0)\\
&&
\langle \hat\chi(cb)+\hat\chi(a)+\hat z, \overline\chi(cb)-\overline\chi(b)-\overline\chi(c)\rangle\\
&&\langle \hat\chi(cba)+\hat z, \overline\chi(cba)-\overline\chi(ba)-\overline\chi(c)\rangle^{-1}
\\
&&
\langle \hat\chi(cba)+\hat z, \overline\chi(cba)-\overline\chi(a)-\overline\chi(cb)\rangle
\\
&&
\langle \hat\chi(ba)+\hat z, \overline\chi(ba)-\overline\chi(a)-\overline\chi(b)\rangle^{-1}
\\
&=&
\phi(c,\overline\chi(b),0)^{-1}\
\phi(c,\overline\chi(ba),0)\
\phi(c,\overline\chi(a),\chi(b))^{-1}\\
&&\langle \hat\chi(c), \overline\chi(ba)-\overline\chi(a)-\overline\chi(b)\rangle\\
&=&
\phi(c,\overline\chi(b),0)^{-1}\
\phi(c,\overline\chi(ba),0)\
\phi(c,\overline\chi(a),\chi(b))^{-1}\\
&& \phi(c, \overline\chi(ba)-\overline\chi(a)-\overline\chi(b),0)^{-1}\\
&=&1.
\end{eqnarray*}
Thus, we have shown that $\partial_{\hat\chi}(\hat\psi,\hat\phi,1)=1$.

To check that all choices involved do not change the class 
$[\hat\psi,\hat\phi,1]$ is tedious and left to the reader.
\end{pf}

It follows that the assignment 
$(\chi,\psi,\phi)\mapsto (\hat\chi,\hat\psi,\hat\phi)$
defines a duality map
\begin{eqnarray}\label{EqTheDualityMap}
\widehat{}\ :{\rm Dyn}^\dag(\Gamma)\to\widehat{\rm Dyn}^\dag(\Gamma).
\end{eqnarray}

\begin{prop}
The duality map (\ref{EqTheDualityMap})
is an isomorphism, and its inverse is 
given by replacing everything by its 
dual counterpart.
\end{prop}

\begin{pf}
We compute the double dual
$(\hat{\hat\chi},\hat{\hat\psi},\hat{\hat\phi})$,
show that $\hat{\hat\chi}=\chi$, and that
$(\hat{\hat\psi},\hat{\hat\phi},1)$ defines the
the same class
in $H^{2,\dag}_{\rm cont}(\Gamma,G,M_\chi)$ as
$({\psi},{\phi},1)$.

Firstly, $\hat{\hat\chi}=\chi$ as according to (\ref {EqTheDualTorusCoycle123})
\begin{eqnarray*}
\langle\hat{\hat\chi}(a),n^\perp\rangle:=\hat\phi(a,n^\perp,z)^{-1}
=\langle n^\perp,\overline\chi(a)\rangle
=\langle n^\perp,\chi(a)\rangle,
\end{eqnarray*}
for $n^\perp\in N^\perp.$

Secondly, we compute $\hat{\hat\phi}$.
So let $\overline{\hat\chi}:\Gamma\to \widehat G$ be a 
lift of $\hat\chi :\Gamma\to \widehat G/N^\perp$,
then according to (\ref{EqTDOTDPhi})
\begin{eqnarray*}
\hat{\hat\phi}(a,g,z)&:=&\langle \overline{\hat\chi}(a),g\rangle^{-1}\\
&=& \langle \overline{\hat\chi}(a),g\rangle^{-1}\ \phi(a,g,z)^{-1}\ \phi(a,g,z)\\
&=& (d\nu)(a,g,z)\ \phi(a,g,z)\\
\end{eqnarray*}
where $\nu:\Gamma\times G/N\to \U(1)$ 
is such that 
$(d\nu)(a,g,z):=\nu(a,gN+z)\nu(a,z)^{-1}
= \langle \overline{\hat\chi}(a),g\rangle^{-1} \phi(a,g,z)^{-1}.$
Such $\nu$ exists as $g\mapsto \langle \overline{\hat\chi}(a),g\rangle^{-1}\ \phi(a,g,z)^{-1}$ factors over $G/N$ by definition of $\hat\chi$, and hence
$(g,z)\mapsto \langle \overline{\hat\chi}(a),g\rangle^{-1}\ \phi(a,g,z)^{-1}$
is a boundary.

Thridly, we have to compute $\hat{\hat\psi}$ according to (\ref{EqTheDeOTDPsi}):
\begin{eqnarray*}
\hat{\hat\psi}(a,b,z)&:=&
\hat\psi(a,b,0)\ \hat\phi(b,\overline{\hat\chi}(a),0)^{-1}\\
&&
\langle \overline{\hat\chi}(ba)-\overline{\hat\chi}(a)-\overline{\hat\chi}(b),\chi(ba)+ z\rangle.
\end{eqnarray*}
The calculation is rather lengthy, but straight forward. 
It only uses the cocycle condition for $(\psi,\phi,1)$ 
and the definitions of $\hat\chi,\hat\phi,\hat\psi$.
The details are left to the reader, we just state the result,
this is
\begin{eqnarray*}
\hat{\hat\psi}(a,b,z)= \psi(a,b,z)\  (\delta_\chi\nu_1)(a,b,z),
\end{eqnarray*}
where $\nu_1:\Gamma\times G/N\to\U(1)$ is defined by
\begin{eqnarray*}
\nu_1(a,z):= \nu(a,z)\ \nu(a,-\chi(a))^{-1} \phi(a,-\overline\chi(a),0)^{-1}.
\end{eqnarray*}
Note that $d\nu_1=d\nu$, therefore
we have shown that (written multiplicative)
\begin{eqnarray*}
(\hat{\hat\psi},\hat{\hat\phi},1)=(\psi,\psi,1)\cdot \partial_\chi(\nu_1,1),
\end{eqnarray*}
and thus the proposition follows.
\end{pf}

It is easily checked now that the duality map (\ref{EqTheDualityMap}) 
is natural in $\Gamma$, and therefore we even have constructed a
natural equivalence of set-valued functors
\begin{eqnarray*}
{\rm Dyn}^\dag \cong \widehat{\rm Dyn}^\dag.
\end{eqnarray*}

\section{The Duality Theory, Crossed Products and K-Theory}
\label{SecTDATCAODT}

\noindent
In this section we show that the duality 
theory of dynamical triples coincides with 
the $C^*$-algebraic duality theory using
crossed product $C^*$-algebras.
\\
\\
Let $(\chi,\lambda,\mu)$ be a dualisable dynamical 
triple. As in section \ref{SecTCODDT} we have 
continuous lifts $\overline\lambda, \overline\mu$,
where $\overline\mu$ still satisfies the cocycle condition.
Canonically, we identify the $C^*$-algebra of sections of the associated 
$C^*$-bundle 
$$(G/N\times\PU(\HH))\times_{\PU(\HH)}\K(\HH)\to G/N$$ 
with the functions $C(G/N,\K(\HH))$. 
It has an action of $G$ which is given by
\begin{eqnarray*}
(g\cdot f)(z):= \mu(g,z)^{-1}[f(z+gN)].
\end{eqnarray*}
For this action we denote by $G\ltimes C(G/N,\K(\HH))$
the corresponding crossed product $C^*$-algebra \cite{Pe}.
In \cite[Thm. 3.7]{Sch} an isomorphism of $C^*$-algebras 
\begin{eqnarray*}
T_{\overline\mu}:   G\ltimes C(G/N,\K(\HH))
\to C(\widehat G/N^\perp,\K(L^2(N^\perp)\otimes\HH))
\end{eqnarray*}
has been introduced which  on the dense subspace of compactly supported
functions $f\in C_c(G\times G/N, \K(\HH))\subset G\ltimes C(G/N,\K(\HH))$
is given by
\begin{eqnarray*}
(T_{\overline\mu}f)(\hat g N^\perp)
:=\mathcal F_{G/N}\circ \langle -\hat g, \sigma(\_)\rangle\circ\mathcal F_{G/N}^{-1}
\circ f^{\overline\mu}(\hat g)\circ
\mathcal F_{G/N}\circ \langle \hat g, \sigma(\_)\rangle\circ\mathcal F_{G/N}^{-1},
\end{eqnarray*}
where $\mathcal F_{G/N}:L^2(G/N)\to L^2(N^\perp)$ is the Fourier transform
and $f^{\overline\mu}(\hat g)\in \K(L^2(N^\perp)\otimes\HH)$ is the Hilbert-Schmidt 
operator whose integral kernel is 
$$
N^\perp\times N^\perp\ni (m,n)\mapsto 
\mathcal F_{G\times G/N}(f\cdot\overline\mu)(\hat g+n,n-m)
\in \K(\HH).
$$
The crossed product $G\ltimes C(G/N,\K(\HH))$ 
comes along with its natural action of $\widehat G$, 
and   because  the actions of $G$ and $\Gamma$
on $C(G/N,\K(\HH))$ commute, the crossed product
inherits  also a $\Gamma$-action commuting with $\widehat G$.
For $f\in C_c(G\times G/N, \K(\HH))$ these actions are
\begin{eqnarray*}
(\hat g \cdot f)(g,z):= \langle \hat g,g\rangle f(g,z), \quad
(a\cdot f)(g,z):= \lambda(a,z)^{-1}[f(g,z+\chi(a))].
\end{eqnarray*}
So, the isomorphism $T_{\overline\mu}$ induces 
two commuting actions of $\widehat G$ and $\Gamma$ on 
$C(\widehat G/N^\perp,\K(L^2(N^\perp)\otimes\HH))$.
It is lengthy but straight forward to compute these induced 
actions, namely, the following identities holds
\begin{eqnarray*}
T_{\overline\mu}(\hat g\cdot f)(\hat z)
=\hat\mu(\hat g,\hat z)^{-1}[ T_{\overline\mu}(f)(\hat z+\hat gN^\perp)],
\end{eqnarray*}
where $\hat\mu(\hat g,\hat z):={\rm Ad}( \mathcal F_{_{G/N}}\circ
\langle \hat g,-\sigma(\_)\rangle\circ\mathcal F_{_{G/N}}^{-1})$.
For the action of $\Gamma$ we find
\begin{eqnarray*}
T_{\overline\mu}(a\cdot f)(\hat z)
=\hat\lambda(a,\hat z)^{-1}[ T_{\overline\mu}(f)(\hat z+\hat \chi(a))],
\end{eqnarray*}
where
\begin{eqnarray*}
\hat\lambda(a,\hat z)&:=&{\rm Ad}
\Big( \mathcal F_{_{G/N}}\circ
\langle \hat\chi(a)+\hat z,\sigma(\_+\chi(a))-\sigma(\_)-\sigma(\chi(a))\rangle\\
&&\lambda_{_{G/N}}(-\chi(a))\overline\lambda(a,-\_)\phi(a,-\sigma(\_),0)^{-1}
\circ\mathcal F_{_{G/N}}^{-1}
\Big)
\end{eqnarray*}
and $\hat\chi$ is from (\ref{EqTheDualTorusCoycle123}).
The doubtful reader is advised to have a look at 
section 3.4 of \cite{Sch}, where essentially the same
calculations are done.
 
We see from the structure of the induced actions,
that $(\hat\chi,\hat\lambda,\hat\mu)$ is a dualisable
dual dynamical triple with underlying Hilbert space 
$L^2(N^\perp)\otimes\HH$.
Now, the following lemma links the $C^*$-algebraic duality
theory to the duality theory of the previous section.

\begin{lem}\label{LemTheLinkBetweenMeAndMrTakai}
The class of $(\hat\chi,\hat\lambda,\hat\mu)$
in $\widehat{\rm Dyn}^\dag(\Gamma,\hat\chi)\cong 
H^{2,\dag}_{\rm cont}(\Gamma,\widehat G,M_{\hat\chi})$ coincides 
with the dual of the class of $(\chi,\lambda,\mu)$ given
by the duality map (\ref{EqTheDualityMap}).
\end{lem} 

\begin{pf}
We have to compute the cocycle $(\hat\psi,\hat\phi,1)$
of  $(\hat\chi,\hat\lambda,\hat\mu)$ according to
(\ref{EqTheDefOfPsPhOm}) for continuous lifts
$\overline{\hat\lambda},\overline{\hat\mu}$.
In the definitions of $\hat\lambda$ and $\hat\mu$
the arguments of ${\rm Ad}$ are continuous 
functions and we define them to be  
$\overline{\hat\lambda}$ and $\overline{\hat\mu}$.
Then the calculation is again straight forward and leads to 
\begin{eqnarray*}
\overline{\hat\lambda}(a,\hat z+\hat gN^\perp)\overline{\hat\mu}(\hat g,\hat z)
=\overline\mu(\hat g,\hat z+\hat\chi(a))\overline{\hat\lambda}(a,\hat z)
\ \underbrace{\langle \hat g ,\sigma(\chi(a)) \rangle^{-1}},\\
=:\hat\phi(a,\hat g, z)
\hspace{-0,5cm}
\end{eqnarray*}
and (somewhat more lengthy to compute)  we also find
\begin{eqnarray*}
&&\overline{\hat\lambda}(ba,\hat z^{-1}\overline{\hat\lambda}(b,\hat z+\hat\chi(a))\overline{\hat\lambda}(a,\hat z)\\
&&=
\underbrace{ \psi(a,b,0)\phi(b,\sigma(\chi(a)),0)^{-1}
\langle \hat\chi(ba)+\hat z,\sigma(\chi(ba))-\sigma(\chi(b))-\sigma(\chi(a))\rangle}.\\
&&\hspace{6cm}=:\hat\psi(a,b,z)
\end{eqnarray*}
A comparison with (\ref{EqTDOTDPhi}) and 
(\ref{EqTheDeOTDPsi})  proves the lemma.
\end{pf}

\begin{rem}
Combining Lemma \ref{LemTheLinkBetweenMeAndMrTakai} 
and the Takai duality theorem (see \cite{Pe})
we also obtain that the duality map (\ref{EqTheDualityMap})
is an isomorphism.
However, we obtained that result completely independent of any
theory about $C^*$-algebras.
\end{rem}

If we fix the groups $G:=\R^n$, $N:=\Z^n$
we have a K-theoretic application of 
Lemma \ref{LemTheLinkBetweenMeAndMrTakai}.
Namely, consider a dualisable dynamical triple $(\rho,E,P)$
and a dualisable dual dynamical triple $(\hat\rho,\widehat E,\widehat P)$
such that their classes match under the duality map
(\ref{EqTheDualityMap}).
Then there is an isomorphism in 
twisted equivariant K-theory
\begin{eqnarray}\label{EqAlterDiePapersKlebenNicht}
K^{*-n}_{\Gamma,P}(E)\cong K^{*}_{\Gamma,\widehat P}(\widehat E).
\end{eqnarray}
In fact, 
by  Connes' Thom isomorphism 
$K_*(\R^n\ltimes A)\cong K_{*-n}(A)$, for any
$\R^n$-$C^*$-algebra $A$ \cite{Con},
we have a chain of isomorphisms 
\begin{eqnarray*}
K^*_{\Gamma,P}(E)&\cong& 
K_*(\Gamma\ltimes \Gamma(E,P\times_{\PU(\HH)}\K(\HH)))\\
&\cong& 
K_{*+n}(\R^n\ltimes\Gamma\ltimes \Gamma(E,P\times_{\PU(\HH)}\K(\HH)))\\
&\cong& 
K_{*+n}(\Gamma\ltimes\R^n\ltimes \Gamma(E,P\times_{\PU(\HH)}\K(\HH)))\\
&\cong& 
K_{*+n}(\Gamma\ltimes \Gamma(\widehat E,\widehat P\times_{\PU(\widehat \HH)}\K(\widehat \HH)))\\
&\cong& 
K^{*+n}_{\Gamma,P}(E).
\end{eqnarray*}

\section{Topological Triples}
\label{SecTopologicalTriples}

\noindent
In this section we introduce the main objects of our interest 
-- ($\Gamma$-equivariant) topological triples.
In case of $G/N=\R^n/\Z^n=\widehat G/N^\perp$
their non-equivariant version has been introduced first 
in \cite{BS} under the name T-duality triples.
For general compact abelian groups the name topological 
triples has been used in \cite{Sch}. 
\\

There is a canonical
$\U(1)$-principal fibre bundle over
$G/N\times \widehat G/N^\perp$
which is called Poincaré bundle.
We shortly recall its definition. 
Let 
\begin{equation}\label{EqThePoincareBundle}
Q:= \big(G/N\times \widehat G\times \U(1)\big)/N^\perp,
\end{equation}
where the action of $N^\perp$ is defined  by
$(z,\hat g,t)\cdot n^\perp:=(z,\hat g+n^\perp,t\ \langle n^\perp,z\rangle^{-1})$.
Then the obvious map
$Q\to G/N\times \widehat G/N^\perp$ is a $\U(1)$-principal fibre bundle.
Dually, we may exchange $G$ and $\widehat G$ 
to find a second $\U(1)$-bundle 
$$
R:=(G\times \widehat G/N^\perp\times \U(1))/N\to G/N\times\widehat G/N^\perp.
$$ 

We denote the $\check{\rm C}$ech classes
in $\check H^1(G/N\times \widehat G/N^\perp,\underline{\U(1)})$
defined by these bundles by $[Q]$ and $[R]$.

\begin{defi}\label{DefiOfPoincareClass}
The class $\pi:=-[Q]$ constucted above is
called the {\bf Poincaré class} of $G/N\times\widehat G/N^\perp$.
\end{defi}

In fact, there is no proper choice  for this definition 
as the following is true.

\begin{lem}
$[Q]$ and $[R]$ are inverses of each other, i.e.
$$[Q]+[R]=0\in \check H^1(G/N\times \widehat G/N^\perp,\underline{\U(1)}).
$$
\end{lem}

\begin{pf}
See Lemma 2.2 of \cite{Sch}.
\end{pf}

We now turn to the definition of  topological 
triples.
Let $P\to E\to *$ be a  pair and 
let  $\widehat P\to\widehat E\to *$
be a dual pair with same underlying Hilbert space $\HH$.
So we have the following
diagram in the $\Gamma$-equivariant category
\begin{equation}\label{DiagAnAlmostTDualityDiagram}
\xymatrix{
&& P\times \widehat E\ar[rd]\ar[ld]&
&E\times\widehat P\ar[rd]\ar[ld]&\\
&P \ar[rd]&&
E\times \widehat E
\ar[rd]\ar[ld]&& \widehat P\ar[ld]\\
&&E\ar[rd]&&\widehat E\ar[ld]&\\
&&&\ast&.
}
\end{equation}
Assume that there is  a $\Gamma$-equivariant $\PU(\HH)$-bundle isomorphism
$\kappa: E\times \widehat P\to  P\times \widehat E
$
(being the identity on the  base $E\times \widehat E$)
which makes  diagram 
(\ref{DiagAnAlmostTDualityDiagram})
commute. 
We choose charts 
for the pair and the dual pair to
trivialise (\ref{DiagAnAlmostTDualityDiagram}).
This induces 
an automorphism  of the
trivial $\PU(\HH)$-bundle over \label{PageOfKappai}
$G/N\times \widehat G/N^\perp$,
$$
\resizebox*{!}{2,5cm}{
\xymatrix{
G/N\times \widehat G/N^\perp\times
\PU(\HH)\ar[rdd]&&G/N\times \widehat G/N^\perp\times
\PU(\HH)\ar[ll]\ar[ldd]\\
&&\\
&G/N\times\widehat G/N^\perp,&
}
}
$$
which is given by a function
(which we also denote by)
$$
\kappa:G/N\times \widehat G/N^\perp\to \PU(\HH),
$$
and this function defines a $\check{\rm C}$ech class
$[\kappa]\in\check H^1(G/N\times \widehat G/N^\perp,\underline{\U(1)})$.

\begin{defi}
Let $L\in\Z$ be any integer.
$\kappa$  satisfies the {\bf Poincaré condition of order $\boldsymbol L$}
if  the equality  $[\kappa]= L\cdot\pi+p_1^*x+ p_2^*y$ holds,
for the Poincaré class $\pi$ and some classes
$x\in \check H^1(G/N,{\underline{\U(1)}})$ and
$y\in \check H^1(\widehat G/N^\perp,{\underline{\U(1)}})$.
Here $p_1,p_2$ are the projections from $G/N\times\widehat G/N^\perp$
on the first and second factor.
\end{defi}
Note that in this definition the classes $x,y$ are just of minor importance.
They are manifestations of the freedom to choose another chart. 
In fact, one can always modify  the given charts
such that $x$ and $y$ vanish.

\begin{defi}\label{DefiTopTriples}
Let $L\in\Z$ be any integer.
A  ($\Gamma$-equivariant) {\bf topological  triple of order $\boldsymbol L$}
(over the point), 
$\big(\kappa,(P,E),(\widehat P,\widehat E)\big)$,
 is a commutative diagram
\begin{equation}\label{DiagTDualityDiagram}
\xymatrix{
& 
& P \times \widehat E
\ar[rd]\ar[ld]&
&E\times\widehat P\ar[rd]\ar[ld]\ar[ll]_\kappa
&\\
&P \ar[rd]&&
E\times \widehat E
\ar[rd]\ar[ld]&\quad\qquad& \widehat P\ar[ld]\\
&&E\ar[rd]&&\widehat E\ar[ld]&\\
&&&\ast&,&
}
\end{equation}
wherein $(P,E)$ is a pair, $(\widehat P,\widehat E)$ is
 a dual pair (with same underlying Hilbert space $\HH$),
 and $\kappa$ is an isomorphism that satisfies
the Poincaré condition of order $L$.
\end{defi}

We call two topological triples 
$\big(\kappa,(P,E),(\widehat P,\widehat E)\big)$ 
and
$\big(\kappa',(P',E'),(\widehat P',\widehat E')\big)$
(of same order, but with underlying Hilbert spaces $\HH,\HH'$ respectively)
{\bf isomorphic} if there is  a morphism of pairs $(\varphi,\vartheta,\theta)$ from
$(P,E)$ to $(P',E')$ and a morphism of 
dual pairs $(\hat\varphi,\hat \vartheta,\hat \theta)$
from $(\widehat P,\widehat E)$
to $(\widehat P',\widehat E')$ such that
the induced diagram
\begin{equation}\label{DiagMorphOfTopTriples}
\xymatrix{
P\times\widehat E\ar[d]^{\vartheta\times\hat\theta}
&&E\times\widehat P\ar[ll]_\kappa \ar[d]^{\theta\times_{_{}}\hat\vartheta}\\
\varphi^*P'\times\widehat E'&&E'\times\hat\varphi^*\widehat P'\ar[ll]_{\kappa'}
}
\end{equation}
is  strictly commutative.
The triples are called {\bf stably isomorphic}
if the stabilised triples
$\big(\Eins\otimes\kappa,(P_{\HH_1},E),(\widehat P_{\HH_1},\widehat E)\big)$ and
$\big(\Eins\otimes\kappa',(P'_{\HH_1},E'),(\widehat P'_{\HH_1},\widehat E')\big)$
are isomorphic  for some separable Hilbert space $\HH_1$.
The meaning 
of the index $\HH_1$ is stabilisation as in equation (\ref{EqStablilisedPUBundleP}).

We call two topological triples 
$\big(\kappa,(P,E),(\widehat P,\widehat E)\big)$ 
and
$\big(\kappa',(P',E'),(\widehat P',\widehat E')\big)$
(of same order, but with underlying Hilbert spaces $\HH,\HH'$ respectively)
{\bf equivalent} if the
pairs $(P,E)$ and  $(P',E')$ and the dual pairs 
$(\widehat P,\widehat E)$
and $(\widehat P',\widehat E')$ are outer conjugate
and the  joint equivariant $\check{\rm C}$ech class 
(see Appendix \ref{SecTECeCOBI}) of 
the bundle isomorphisms $\kappa$ and $\kappa'$
vanishes:
\begin{eqnarray*}
[\kappa,\kappa']=0\in \check 
H^1_\Gamma(G/N\times\widehat G/N^\perp,\underline{\U(1)}).
\end{eqnarray*}
The triples are called {\bf stably equivalent}
if for a separable Hilbert space $\HH_1$ the stabilised triples
$\big(\Eins\otimes\kappa,(P_{\HH_1},E),(\widehat P_{\HH_1},\widehat E)\big)$ and
$\big(\Eins\otimes\kappa',(P'_{\HH_1},E'),(\widehat P'_{\HH_1},\widehat E')\big)$
are equivalent.

These notions can be arranged in a diagram of implications
\begin{eqnarray*}
\xymatrix{
{\begin{array}{c} 
\textrm{ isomorphism} \\ 
\textrm{of top. triples} 
\end{array}}
\ar@{=>}[r]\ar@{=>}[d]&
{\begin{array}{c} 
\textrm{equivalence } \\ 
\textrm{of top. triples} 
\end{array}}
\ar@{=>}[d]\\
{\begin{array}{c} 
\textrm{stable isomorphism} \\ 
\textrm{of top. triples} 
\end{array}}
\ar@{=>}[r]&
{\begin{array}{c} 
\textrm{stable equivalence} \\ 
\textrm{of top. triples} 
\end{array}}.
}
\end{eqnarray*}
Stable equivalence will be the correct notion of equivalence for us, and
we introduce a set valued functor ${\rm Top}_L$ which \label{PageOfTopFunc}
associates to a finite group $\Gamma$ the set of stable equivalence classes
of topological triples of order $L$, i.e.
$$\label{PageTopEB}
{\rm Top}_L(\Gamma):=\{\text{ topological triples of order $L$} \}\big/_\textrm{st. 
equiv.}.
$$
If we fix a homomorphism $\chi:\Gamma\to G/N$-bundle,
we can consider all  topological
triples for which (after a chart is chosen)
the $\Gamma$-action on $E$
is given by $\chi$.
We define
$$
{\rm Top}_L(\Gamma,\chi):=\{\textrm{topological triples of order $L$ with fixed }\chi\}
\big/_{\textrm{st. equiv.  }},
$$
so there is a decomposition 
${\rm Top}_L(\Gamma)=\coprod_{\chi} {\rm Top}_L(\chi,\Gamma)$.
We also let 
$$
{\rm Top}_L(\Gamma,\chi,\hat\chi):=\{\textrm{topological triples of 
order $L$ with fixed }\chi\ {\rm and}\ \hat\chi\}
\big/_{\textrm{st. equiv.  }},
$$
for fixed  
 $\chi:\Gamma\to G/N$ 
and $\hat\chi:\Gamma\to \widehat G/N^\perp$ in the obvious way. 
This yields a decomposition
${\rm Top}_L(\Gamma)=\coprod_{\chi,\hat\chi}{\rm Top}_L(\Gamma,\chi,\hat\chi)$.

The notion of stable equivalence on topological triples
is made such that we have natural forgetful maps
\begin{eqnarray*}
\xymatrix{
&{\rm Top}_L(\Gamma)\ar[ld]\ar[dr]&\\
{\rm Par}(\Gamma)&&\widehat{\rm Par}(\Gamma).
}
\end{eqnarray*}

There are  also natural maps
\begin{eqnarray*}
{\rm Top}_L(\Gamma,\chi,\hat\chi)\times {\rm Top}_K(\Gamma,\chi,\hat\chi)
&\to& {\rm  Top}_{L+K}(\Gamma,\chi,\hat\chi),\\
{\rm Top}_L(\Gamma,\chi,\hat\chi)&\to& {\rm Top}_{-L}(\Gamma,\chi,\hat\chi)
\end{eqnarray*}  
induced by  tensor product
and complex conjugation of projective unitary
bundles.
These maps turn 
$
{\rm Top}_0(\Gamma,\chi,\hat\chi)$ into an abelian group,
where the class of  $(\Eins,(\chi,\Eins),(\hat\chi,\Eins))$
is the unit.
Each 
${\rm Top}_L(\Gamma,\chi,\hat\chi) $ becomes a  
${\rm Top}_0(\Gamma,\chi,\hat\chi)$-torsor, and
${\rm Top}_\bullet(\Gamma,\chi,\hat\chi):=\coprod_{L\in\Z} T_L(\Gamma,\chi,\hat\chi)$
becomes an abelian group which fits 
into the short exact sequence
\begin{eqnarray*}
\xymatrix{
0\ar[r]& {\rm Top}_0(\Gamma,\chi,\hat\chi)\ar[r]&{\rm Top}_\bullet(\Gamma,\chi,\hat\chi)\ar[r]&\Z\ar[r]& 0.
}
\end{eqnarray*}

\begin{rem}\label{RemOnTorsorStories}
Note that the empty set is a torsor for any group, 
and that the quotient of the empty set by any group
is the one point space.
Therefore the above short exact sequence is 
meaningful even if 
${\rm Top}_L(\Gamma,\chi,\hat\chi)=\emptyset$
for all $L\not= 0$.
${\rm Top}_0(\Gamma,\chi,\hat\chi)$  is never empty, 
as it contains at least its unit $[\Eins,(\chi,\Eins),(\hat\chi,\Eins)]$. 
\end{rem}

By its very definition, the Poincar\'e class $\pi$ has a geometric 
interpretation by the Poincar\'e bundle (\ref{EqThePoincareBundle}).
However, what we need is an analytical description of $\pi$ which is
more adequate for the analysis we are going to do.
So let $\sigma:G/N\to G$ and $\hat\sigma:\widehat G/N^\perp\to\widehat G$ be 
Borel sections of the quotient maps, then
for each fixed $(z,\hat z)\in G/N\times\widehat G/N^\perp$
there are multiplication operators
\begin{eqnarray*}
\overline\kappa^\sigma(z,\hat z)&:=&
\langle\hat\sigma(\hat z),\sigma(\_-z)-\sigma(\_)\rangle
\in L^\infty(G/N,\U(1))\subset\U(L^2(G/N)),\\
\overline{\hat\kappa}^{\hat\sigma}(z,\hat z)&:=&
\langle\hat\sigma(..-\hat z)-\hat\sigma(..),\sigma(z)\rangle
\in L^\infty(\widehat G/N^\perp,\U(1))\subset\U(L^2(\widehat G/N^\perp)),
\end{eqnarray*}
wherein for a $\U(\HH)$-valued function $f$ on $G/N$, 
or $g$ on $\widehat G/N^\perp$, we denoted
by $f(\_)$ and  $g(..)$, respectively, the corresponding 
multiplication operator; for example if 
$F\in L^2(G/N)$, then 
$(f(\_)F)(x):=f(x)F(x)$.
These operators give rise to the following  projective unitary functions
\begin{eqnarray}
\kappa^\sigma:
G/N\times \widehat G/N^\perp&\to&\PU(L^2(G/N))\label{GlDefGlfuerKappa},\\
(z,\hat z)&\mapsto & {\rm Ad}\big( \overline\kappa^\sigma(z,\hat z) \big)\nonumber\\
\nonumber\\
\hat\kappa^{\hat\sigma}:
G/N\times \widehat G/N^\perp&\to&\PU(L^2(\widehat G/N^\perp)).\\
(z,\hat z)&\mapsto&{\rm Ad}\big(\overline{\hat\kappa}^{\hat\sigma}(z,\hat z)\big)\nonumber
\end{eqnarray}

\begin{lem}\label{LemTheDefiOfThePoincareClass}
The maps $\kappa^\sigma$ and $\hat\kappa^{\hat\sigma}$ are continuous,
$\kappa^\sigma$ is independent of $\hat\sigma$ and $\hat\kappa^{\hat\sigma}$
is independent of $\sigma$.

Moreover, the $\check C$ech classes 
$
[\kappa^\sigma],[\hat\kappa^{\hat\sigma}]\in 
\check H^1(G/N\times\widehat G/N^\perp,\underline{\U(1)}),
$
of the two maps are independent of $\sigma$ and $\hat\sigma$ and 
$$
\pi=[\kappa^\sigma]=-[\hat\kappa^{\hat\sigma}].
$$
\end{lem} 
\begin{pf}
See Lemma 2.3 and Lemma 2.4 of \cite{Sch}.
\end{pf}

It is worth mentioning that $ \kappa^\sigma(0,\hat z)=\Eins=\kappa^\sigma(z,0)$.
The first equality is trivial, the second follows from $\hat\sigma(0)\in N^\perp$,
so 
$\langle\hat\sigma(0),\sigma(\_-z)-\sigma(\_)\rangle=\langle\hat\sigma(0),\sigma(-z)\rangle\in\U(1).
$
Clearly, the same is true for $\hat\kappa^{\hat\sigma}$.
By similar reasoning, we 
can prove  the following Lemma.

\begin{lem}\label{LemKLAPPTOPI}
For  integers $K,L\in \Z$  consider the map $(K,L):(z,\hat z)\mapsto (K z,L \hat z)$.
Its induced map 
\begin{eqnarray*}
(K,L)^*:\check H^1(G/N\times\widehat G/N^\perp,\underline{\U(1)})\to
\check H^1(G/N\times\widehat G/N^\perp,\underline{\U(1)})
\end{eqnarray*}
maps the Poincaré class $\pi$ to $K\cdot L\cdot \pi$.
\end{lem}

\begin{pf}
We show  that 
$\kappa^\sigma(Kz,L\hat z)$ equals $\kappa^\sigma(z,\hat z)^{KL}$ up
to a unitary function.

As $\langle \hat\sigma(L\hat z),\sigma(z)\rangle\in\U(1)$ and
$\sigma(\_-z)-\sigma(\_)+\sigma(z)\in N$, we have
first
\begin{eqnarray*}
\kappa^\sigma(z,L\hat z)&=&
{\rm Ad}(\langle \hat \sigma(L\hat z),\sigma(\_-z)-\sigma(\_)+\sigma(z)\rangle)\\
&=&{\rm Ad}(\langle L\hat z,\sigma(\_-z)-\sigma(\_)+\sigma(z)\rangle)\\
&=&{\rm Ad}(\langle \hat z,\sigma(\_-z)-\sigma(\_)+\sigma(z)\rangle)^L.\\
\end{eqnarray*}
Second, by use of 
$\sigma(\_-Kz)-K\sigma(\_-z)+(K-1)\sigma(\_)\in N$,
we have
\begin{eqnarray*}
\kappa^\sigma(Kz,\hat z)&=&
{\rm Ad}(\langle \hat \sigma(\hat z),\sigma(\_-Kz)-\sigma(\_)\rangle)\\
&=&{\rm Ad}(\langle \hat\sigma(\hat z),\sigma(\_-Kz)-K\sigma(\_-z)+(K-1)\sigma(\_)\rangle)\\
&&{\rm Ad}(\langle\hat \sigma(\hat z),K\sigma(\_-z)-K\sigma(\_)\rangle)\\
&=&{\rm Ad}(\langle \hat z,\sigma(\_-Kz)-K\sigma(\_-z)+(K-1)\sigma(\_)\rangle)\\
&&{\rm Ad}(\langle\hat \sigma(\hat z),\sigma(\_-z)-\sigma(\_)\rangle)^K.
\end{eqnarray*}
Finally,  note that
\begin{eqnarray*}
(z,\hat z)\mapsto \langle \hat z,\sigma(\_-Kz)-K\sigma(\_-z)+(K-1)\sigma(\_)\rangle\in \U(L^2(G/N))
\end{eqnarray*}
is continuous. In fact, the map 
\begin{eqnarray*}
G/N\times \widehat G\ni (z,\hat g)\mapsto \langle \hat g,\sigma(\_-Kz)-K\sigma(\_-z)+(K-1)\sigma(\_)\rangle\end{eqnarray*}
factors over $G/N\times \widehat G/N^\perp$ and 
$\langle \hat g,\sigma(\_-z)\rangle=\lambda_{_{G/N}}(z) \langle \hat g,\sigma(\_)\rangle
\lambda_{_{G/N}}(z)^{-1}$
is continuous by the same argument as in the proof of Lemma \ref{LemTechLFMu22}.
\end{pf}

By Lemma \ref{LemKLAPPTOPI} we can define for any 
triple of integers $K,L,M\in\Z$ a natural map
\begin{eqnarray}\label{EqThereAreNMF!TOL}
{\rm Top}_M(\Gamma,K\chi,L\hat\chi)\to{\rm Top}_{KLM}(\Gamma,\chi,\hat\chi),
\end{eqnarray}
namely, let $(\kappa,(K\chi, \lambda),(L \hat\chi,\hat\lambda))$
be a topological triple of order $M$, then the isomorphism $\kappa$
is given by a continuous function 
(also denoted by) 
$\kappa:G/N\times\widehat G/N^\perp\to\PU(\HH)$ 
such that
\begin{eqnarray*}
\lambda(a,z)\kappa(z,\hat z)=\kappa(z+(K \chi)(a),\hat z+(L\hat\chi)(a))\hat\lambda(a,\hat z)
\end{eqnarray*}
is satisfied.
This implies that
\begin{eqnarray*}
(K^*\lambda)(a,z)((K,L)^*\kappa)(z,\hat z)=((K,L)^*\kappa)(z+\chi(a),\hat z+\hat\chi(a))
(L^*\hat\lambda)(a,\hat z)
\end{eqnarray*}
is also satisfied, where 
$(K^*\lambda)(a,z):=\lambda(a,Kz)$,
$(L^*\hat\lambda)(a,\hat z):=\hat\lambda(a,L\hat z)$
and
$((K,L)^*\kappa)(z,\hat z):=\kappa(Kz,L\hat z)$.
It is clear that $(\chi,K^*\lambda)$ is a pair and that
$(\hat\chi,L^*\hat\lambda)$ is a dual pair, 
and Lemma \ref{LemKLAPPTOPI}
implies that $(K,L)^*\kappa$ satisfies the Poincaré condition
of order $KLM$. Hence,
$( (K,L)^*\kappa, (\chi,K^*\lambda),(\hat\chi,L^*\hat\lambda))$ 
is a topological triple of order $KLM$.
This construction defines (\ref{EqThereAreNMF!TOL}).

In section \ref{SecRegET} we introduce the notion of
regular topological triples, and we will see that the restriction 
of (\ref{EqThereAreNMF!TOL}) to the set of regular topological 
triples is an isomorphism. 
In particular, we will consider  the special cases $M=K=1$ and $M=L=1$,
then (\ref{EqThereAreNMF!TOL}) becomes
\begin{eqnarray}\label{EqFromTOp1ToTopL}
{\rm Top}_1(\Gamma,\chi,L\hat\chi)\to{\rm Top}_{L}(\Gamma,\chi,\hat\chi),\\
{\rm Top}_1(\Gamma,K\chi,\hat\chi)\to{\rm Top}_{K}(\Gamma,\chi,\hat\chi).\nonumber
\end{eqnarray}
This means 
that at least partial information of 
higher order triples is encoded 
in  triples of order 1.

It is clear that the maps  of (\ref{EqFromTOp1ToTopL})
respect the structure of the underlying pairs
and dual pairs, i.e. we have commutative squares 
\begin{eqnarray*}
\xymatrix{
{\rm Top}_1(\Gamma,\chi,L\hat\chi)\ar[d]\ar[r]&{\rm Top}_{L}(\Gamma,\chi,\hat\chi),\ar[d]&
{\rm Top}_1(\Gamma,K\chi,\hat\chi)\ar[r]\ar[d]&{\rm Top}_{K}(\Gamma,\chi,\hat\chi).\ar[d]\\
\widehat{\rm Par}(\Gamma,L\hat\chi)\ar[r]^{L^*}&\widehat{\rm Par}(\Gamma,\hat\chi)&
{\rm Par}(\Gamma,K\chi)\ar[r]^{K^*}&{\rm Par}(\Gamma,\chi)\\
}
\end{eqnarray*}
Here $K^*$ and $L^*$ are from (\ref{EqTheLMapOfPairs}), 
and the vertical maps are the natural forgetful maps.

\section{The Classification of Topological Triples}
\label{SecTheClassifOfET}

\noindent
In this section we give a description of 
${\rm Top}_\bullet(\Gamma,\chi,\hat\chi)$ in terms 
of group cohomology.

Fix a pair of homomorphisms
$$
\chi:\Gamma\to G/N,\qquad 
\hat\chi:\Gamma\to\widehat G/N^\perp,
$$
and consider the $\Gamma$-module
$M_{\chi\times\hat\chi}:=C(G/N\times \widehat G/N^\perp,\U(1))$ 
and its submoduls $M_\chi:=C(G/N,\U(1))$ and
$M_{\hat\chi} :=C(\widehat G/N^\perp,\U(1))$,
where $\Gamma$ acts by shift with (the negative of)
$\chi\times \hat\chi$ in the arguments of  
functions. So we have a diagram of 
$\Gamma$-moduls
$$
\resizebox*{!}{2,5cm}{
\xymatrix{
&&M_\chi\ar@{^{(}->}[dr]&&\\
\U(1)\ar@{=}[r]&M_\chi\cap M_{\hat\chi} \ \ \ \ \ \ \ \ \ar@{^{(}->}[ur]\ar@{^{(}->}[dr]&&M_\chi M_{\hat\chi} \ar@{^{(}->}[r]
&M_{\chi\times\hat\chi},\\
&&M_{\hat\chi} \ar@{^{(}->}[ur]&&
}}
$$
which yield to two long 
exact sequences in
group cohomology
\begin{equation}\label{DiagTwoLongExSeq}
\resizebox*{!}{6cm}{
\xymatrix{
&&&\vdots\ar[d]&&\\
&&&H^2(\Gamma,M_\chi\cap M_{\hat\chi} )\ar[d]&&\\
&&&H^2(\Gamma,M_\chi)\oplus H^2(\Gamma,M_{\hat\chi} )\ar[d]^\ominus&&\\
\cdots\ar[r]&
H^1(\Gamma,M_{\chi\times\hat\chi})\ar[r]&
H^1(\Gamma,M_{\chi\times\hat\chi}/M_\chi M_{\hat\chi} )\ar[r]^B&
H^2(\Gamma,M_\chi M_{\hat\chi} )\ar[r]\ar[d]^C&
H^2(\Gamma,M_{\chi\times\hat\chi})\ar[r]&
\cdots,\\
&&&H^3(\Gamma,M_\chi\cap M_{\hat\chi} )\ar[d]&&\\
&&&\vdots&&\\
}}
\end{equation}
where the vertical sequence is the
Mayer-Vietoris sequence of Appendix \ref{SecTheMVSofGC},
and the horizontal sequence is the usual Bockstein
long exact sequence.

Let us introduce
two classes
\begin{eqnarray*}
\chi\cup\hat\chi&\in& H^3(\Gamma, \U(1)),\\
\chi\sqcup\hat\chi&\in& H^2(\Gamma,M_\chi M_{\hat\chi} ),
\end{eqnarray*}
namely, we define a map $\cup$ by the composition
of either two consecutive maps in the following 
commutative diagram
\begin{eqnarray}\label{DiagTheCupProdOfHomos}
\xymatrix{
H^1(\Gamma,G/N)\times H^1(\Gamma,\widehat G/N^\perp)\ar[r]\ar[d]&
H^2(\Gamma,N)\times H^1(\Gamma,\widehat G/N^\perp)\ar[d]\\
H^1(\Gamma,G/N)\times H^2(\Gamma,N^\perp)\ar[r]&
H^3(\Gamma,\U(1)).
}
\end{eqnarray}
Therein the upper horizontal and the left vertical maps are given
by the corresponding  Bockstein homomorphisms in one factor
and the identity in the other. The two maps to the right lower corner 
are given by the canonical  pairings of dual groups, e.g.
a class $([\omega],\hat \chi)\in H^2(\Gamma,N)\times H^1(\Gamma,\widehat G/N^\perp)$ 
is mapped to the class of the 3-cocycle 
$(a,b,c)\mapsto \langle\hat \chi(a),\omega(b,c)\rangle$.
The commutativity of diagram (\ref{DiagTheCupProdOfHomos})
can be explicitly calculated; 
we give the argument in Corollary \ref{CorOnCommOfCup} below.

For the  definition of $\chi\sqcup\hat\chi$ 
we first define a Borel function
\begin{eqnarray}\label{EqDieDickeBerta}
\beta:\Gamma\times G/N\times\widehat G/N^\perp\to \U(1)
\end{eqnarray}
by the explicit formula
$$
\beta(a,z,\hat z):=\langle \hat\sigma(\hat z+\hat\chi(a))-\hat\sigma(\hat z)-\hat\sigma(\hat\chi(a)),
\sigma(z+\chi(a))\rangle \langle\hat\sigma(\hat z),\sigma(\chi(a))\rangle,
$$
where $\sigma:G/N\to G$ and $\hat\sigma:\widehat G/N^\perp\to \widehat G$ are
chosen Borel sections of the quotient maps.
Let us also define 
$$
\chi \underset{\sigma,\hat\sigma}{\sqcup}\hat\chi:=
\delta_{\chi\times \hat\chi}\beta.
$$
Although $\beta$ fails to be continuous, the following statement for
$\chi\underset{\sigma,\hat\sigma}{\sqcup}\hat\chi$ is true.

\begin{lem}\label{LemThePropertiesOfDeltaBeta}
The map 
$\chi \underset{\sigma,\hat\sigma}{\sqcup}\hat\chi:\Gamma^2\times G/N\times \widehat G/N^\perp\to\U(1)$
is continuous,  and we have a decomposition
$$
(\chi \underset{\sigma,\hat\sigma}{\sqcup}\hat\chi)(a,b,z,\hat z)=
\gamma(a,b)\ \langle \hat\eta(a,b),z\rangle^{-1} \
\langle \hat z, \eta(a,b)\rangle,
$$
wherein
\begin{eqnarray*}
\gamma(a,b):=
\langle\hat\sigma(\hat\chi(a)),\sigma(\chi(b))\rangle\ 
\langle \hat\eta(a,b) ,\chi(ba)
\rangle^{-1},
\end{eqnarray*}
and 
$$
\eta(a,b):=  \sigma(\chi(b))-\sigma(\chi(ba))+\sigma(\chi(a))
\in N\cong\widehat{\widehat G/N^\perp}, 
$$
and
$$
\hat\eta(a,b):=  \hat\sigma(\hat\chi(b))-\hat \sigma(\hat\chi
(ba))+\hat\sigma(\hat\chi(a))
\in N^\perp\cong \widehat{G/N}.
$$
\end{lem}

\begin{pf}
One computes  
$$
(\delta_{\chi\times \hat\chi}\beta)(a,b,z,\hat z)=
\beta(b,z+\chi(a),\hat z+\hat\chi(a))\
\beta(ba,z,\hat z)^{-1}\ 
\beta(a,z,\hat z)
$$
right from the definition.  The resulting formula is
\begin{eqnarray*}
(\delta_{\chi\times \hat\chi}\beta)(a,b,z,\hat z)&=&
\langle\hat\sigma(\hat\chi(a)),\sigma(\chi(b))\rangle\\
&&
\langle\hat\sigma(\hat\chi(ba))-\hat\sigma(\hat\chi(a))-\hat\sigma(\hat\chi(b)),\sigma(\chi(ba))
\rangle\\
&&
\langle\hat\sigma(\hat\chi(ba))-\hat\sigma(\hat\chi(a))-\hat\sigma(\hat\chi(b)),\sigma(z)\rangle\\
&&\langle\hat\sigma(\hat z),\sigma(\chi(b))+\sigma(\chi(a))-\sigma(\chi(ba))\rangle
\end{eqnarray*}
which directly implies the statements above.
\end{pf}

By its definition $\chi \underset{\sigma,\hat\sigma}{\sqcup}\hat\chi$
is a Borel boundary, but by the above lemma, 
it is even a continuous cocycle.
So there is a well-defined class
\begin{eqnarray}
\chi\sqcup\hat\chi:=[\chi \underset{\sigma,\hat\sigma}{\sqcup}\hat\chi]\in H^2(\Gamma,M_\chi M_{\hat\chi} ),
\end{eqnarray}
in fact, it is not difficult to check that the class $\chi\sqcup\hat\chi$ 
is independent of the choice of the sections $\sigma$ and $\hat\sigma$.

\begin{cor}\label{CorOnCommOfCup}

\begin{enumerate}
\item Diagram (\ref{DiagTheCupProdOfHomos}) commutes, 
so $\chi\cup\hat\chi$ is well-defined, and
\item the connecting homomorphism $C$ in the vertical long exact sequence of 
(\ref{DiagTwoLongExSeq}) maps the class $\chi\sqcup\hat\chi$ to $\chi\cup\hat\chi$.
\end{enumerate}

\end{cor}

\begin{pf}
1. 
By the lemma above we have
\begin{eqnarray*}
1=(\delta_{\chi\times\hat\chi} \chi \underset{\sigma,\hat\sigma}{\sqcup}\hat\chi)(a,b,c,z,\hat z)
=(\delta_{\chi\times\hat\chi}\gamma)(a,b,c)\ 
\langle \hat\eta(b,c),\chi(a)\rangle^{-1} \
\langle \hat \chi(a), \eta(b,a)\rangle.
\end{eqnarray*}
Thus, we see that the 3-cocycles 
$$
(a,b,c)\mapsto\langle \hat \chi(a), \eta(b,a)\rangle
\quad {\rm and}\quad
(a,b,c)\mapsto\langle \hat\eta(b,c),\chi(a)\rangle$$
are cohomologous, but they represent precisely 
the two compositions in diagram (\ref{DiagTheCupProdOfHomos}).

\item 
2. The second statement follows 
from the definition of the connecting homomorphism (see Appendix \ref{SecTheMVSofGC})
in the vertical long exact sequence of (\ref{DiagTwoLongExSeq}).
\end{pf}

For any integer $L\in\Z$, let
\begin{eqnarray*}
T_L(\Gamma,\chi,\hat\chi)&:=&
\Big\{ ([w],[\psi],[\hat \psi])\Big| B([w])+L \cdot(\chi\sqcup\hat\chi)=[\psi]\ominus[\hat\psi]\Big\}\\
&\subset&
H^1(\Gamma,M_{\chi\times\hat\chi}/M_\chi M_{\hat\chi} )\times H^2(\Gamma,M_\chi)\oplus H^2(\Gamma,M_{\hat\chi} ).
\end{eqnarray*}
The set $T_L(\Gamma,\chi,\hat\chi)$ fits into diagram (\ref{DiagTwoLongExSeq}),
for $L=0$ only it is a pullback:
\begin{equation}\label{DiagTheClassiGroupSitsInAha}
\resizebox*{!}{6cm}{
\xymatrix{
&&&\vdots\ar[d]&&\\
&&&H^2(\Gamma,\U(1))\ar[d]&&\\
&&T_L(\Gamma,\chi,\hat\chi)\ar[r]\ar[d]&H^2(\Gamma,M_\chi)\oplus H^2(\Gamma,M_{\hat\chi} )\ar[d]^
\ominus&&\\
\cdots\ar[r]&
H^1(\Gamma,M_{\chi\times\hat\chi})\ar[r]&
H^1(\Gamma,M_{\chi\times\hat\chi}/M_\chi M_{\hat\chi} )\ar[r]^B&
H^2(\Gamma,M_\chi M_{\hat\chi} )\ar[r]\ar[d]^C&
H^2(\Gamma,M_{\chi\times\hat\chi})\ar[r]&
\cdots\\
&&&H^3(\Gamma,\U(1))\ar[d]&&\\
&&&\vdots&&\\
}}
\end{equation}
The obvious maps
\begin{eqnarray*}
T_L(\Gamma,\chi,\hat\chi)\times T_K(\Gamma,\chi,\hat\chi)
&\to& T_{L+K}(\Gamma,\chi,\hat\chi),\\
T_L(\Gamma,\chi,\hat\chi)&\to& T_{-L}(\Gamma,\chi,\hat\chi)
\end{eqnarray*}  
turn 
$
T_0(\Gamma,\chi,\hat\chi)$ into an abelian group,
and  each 
$T_L(\Gamma,\chi,\hat\chi) $ into a  
$T_0(\Gamma,\chi,\hat\chi)$-torsor.
Moreover, $T_\bullet(\Gamma,\chi,\hat\chi):=\coprod_{L\in\Z} T_L(\Gamma,\chi,\hat\chi)$
becomes an abelian group which fits 
into the short exact sequence
\begin{eqnarray*}
\xymatrix{
0\ar[r]& T_0(\Gamma,\chi,\hat\chi)\ar[r]&T_\bullet(\Gamma,\chi,\hat\chi)\ar[r]&\Z\ar[r]& 0.
}
\end{eqnarray*}
There is another short exact sequence induced by the previous,
\begin{eqnarray*}
\xymatrix{
0\ar[r]& T_0(\Gamma,\chi,\hat\chi)/H^1(\Gamma,M_{\chi\times\hat\chi})\ar[r]&
T_\bullet(\Gamma,\chi,\hat\chi)/H^1(\Gamma,M_{\chi\times\hat\chi})\ar[r]&\Z\ar[r]& 0,
}
\end{eqnarray*}
where $H^1(\Gamma,M_{\chi\times\hat\chi})$ acts in the first factor of 
each $T_L(\Gamma,\chi,\hat\chi)$ via the map
$H^1(\Gamma,M_{\chi\times\hat\chi})\to H^1(\Gamma,M_{\chi\times\hat\chi}/M_\chi M_{\hat\chi} )$.

Now, we state the classification theorem of  
topological  triples.

\begin{thm}\label{TheBigClThmForEquTrip}

There is a natural isomorphism of short exact sequences
\begin{eqnarray*}
\xymatrix{
0\ar[r]& {\rm Top}_0(\Gamma,\chi,\hat\chi)\ar[r]\ar[d]^\cong &
{\rm Top}_\bullet(\Gamma,\chi,\hat\chi)\ar[r]\ar[d]^\cong &\Z\ar[r]\ar[d]^=& 0\\
0\ar[r]& T_0(\Gamma,\chi,\hat\chi)/H^1(\Gamma,M_{\chi\times\hat\chi}) \ar[r]&
T_\bullet(\Gamma,\chi,\hat\chi)/H^1(\Gamma,M_{\chi\times\hat\chi})\ar[r]&\Z\ar[r]& 0.
}
\end{eqnarray*}

\end{thm}

\begin{rem}
As we already noted in Remark \ref{RemOnTorsorStories},
 the empty set is a torsor for any group, 
and the quotient of the empty set by any group
is the one point space.
Thus, Theorem \ref{TheBigClThmForEquTrip}
makes perfectly sense even if 
$T_L(\Gamma,\chi,\hat\chi)=\emptyset$
for all $L\not= 0$.
\end{rem}

\noindent
The proof of Theorem \ref{TheBigClThmForEquTrip}
is the content of the remainder of this section.
\\

We first show how one constructs the
classes 
\begin{eqnarray*}
{[w]}&\in& H^1(\Gamma,M_{\chi\times\hat\chi}/M_\chi M_{\hat\chi} )/H^1(\Gamma,M_{\chi\times\hat\chi}),\\
([\psi],[\hat\psi])&\in& H^2(\Gamma,M_\chi)\oplus H^2(\Gamma,M_{\hat\chi} )
\end{eqnarray*}
of a topological triple.

Let $(\kappa,(P,E),(\widehat P,\widehat E))$ be a topological triple
of order $L$.
We fix  a chart  which defines the corresponding trivialised 
data $(\chi,\lambda')$ and $(\hat\chi,\hat\lambda')$ 
for the pairs $(P,E)$ and $(\widehat P,\widehat E)$.
Without restriction we can assume that  
$\widehat \HH=L^2(G/N)\otimes\HH$ is 
the underlying Hilbert space of the triple. 
The isomorphism $\kappa$ is given by a continuous map which we also denote by    
$\kappa: G/N\times\widehat G/N^\perp\to\PU(\widehat\HH)$. It satisfies 
$$ 
\lambda'(a,z)\ \kappa(z,\hat z)= \kappa(\chi(a)+z,\hat\chi(a)+\hat z)\  \hat\lambda'(a,\hat z).
$$
Due to the Poincaré condition 
the map $\kappa$ is of the  form
$$\kappa(z,\hat z)=
\kappa^x(z)\ \kappa^{\sigma}(z,\hat z)^L\ {\rm Ad}(u(z,\hat z))\ \kappa^y(\hat z),
$$
wherein $\kappa^\sigma(z,\hat z)={\rm Ad}(\langle\hat\sigma(\hat z),\sigma(\_-z)-\sigma(\_)
\rangle)$ 
is from Lemma \ref{LemTheDefiOfThePoincareClass}, 
$\kappa^x, \kappa^y$ are some continuous projective unitary functions, and 
$u: G/N\times \widehat G/N^\perp\to\U(\widehat\HH)$ is a continuous unitary function. 
We introduce some short hands to get rid of $\kappa^x$ and $\kappa^y$.
We let $$\lambda(a,z):=\kappa^x(\chi(a)+z)^{-1}
\lambda'(a,z)\kappa^x(z),$$ 
and analogously  \label{PageOfPrimedTransis}
$$\hat\lambda(a,\hat z):=\kappa^y(\chi(a)+z)
\hat\lambda'(a,z)\kappa^y(z)^{-1}.$$
Thus we have 
\begin{eqnarray}\label{EqDasTrMitDemBegonnen}
\lambda(a,z)&=&\kappa^\sigma(\chi(a)+z,\hat\chi(a)+\hat z)^L
{\rm Ad}(u(\chi(a)+z,\hat \chi(a)+\hat z))\nonumber\\
&&\hat\lambda(a,\hat z){\rm Ad}(u(z,\hat z))^{-1}
\kappa^\sigma(z,\hat z)^{-L}.\label{EqATopTripIsATopTripIsATopTrip}
\end{eqnarray}

\begin{lem}
There exist continuous lifts  
$$
\xymatrix{
&\U(\widehat\HH)\ar[d]\\
\Gamma\times G/N\ar[r]^{\lambda}\ar[ru]^{\overline{\lambda}}&\PU(\widehat\HH)
}
\quad and
\quad
\xymatrix{
&\U(\widehat\HH)\ar[d]\\
\Gamma\times \widehat G/N^\perp\ar[r]^{\hat\lambda}\ar[ru]^{\overline{\hat\lambda}}&\PU
(\widehat\HH)
}.
$$
\end{lem}
\begin{pf}
Take $\hat z=0$ in (\ref{EqATopTripIsATopTripIsATopTrip}) which yields 
$$
\lambda(a,z)=\kappa^\sigma(\chi(a)+z,\hat\chi(a))^L
{\rm Ad}(u(\chi(a)+z,\hat \chi(a)))
\hat\lambda(a,0){\rm Ad}(u(z,0))^{-1}.
$$
Each single factor of the right hand side has a continuous lift, so
$\lambda$ has.
By taking $z=0$ we conclude  that $\hat\lambda$ has a continuous lift.
\end{pf}

By use of this lemma, we define
\begin{eqnarray}\label{EqTheClassesOTC}
[\psi]\in H^2(\Gamma,M_\chi),\quad [\hat\psi]\in H^2(\Gamma,M_{\hat\chi} )
\end{eqnarray}
by 
$$
\overline\lambda(b,\hat\chi(a)+z)\overline\lambda(a,z)=\overline\lambda(ba,z)\ \psi(a,b,z)$$
and 
$$\overline{\hat\lambda}(b,\hat\chi(a)+z) \overline{\hat\lambda}(a,\hat z)=\overline{\hat\lambda}
(ba,\hat z)\ \psi(a,b,z).
$$

Of course, $\kappa^\sigma$ does not have a continuous unitary lift, but 
we take the Borel lift 
$$
(z,\hat z)\mapsto \overline\kappa^\sigma(z,\hat z)
:=\langle\hat\sigma(\hat z),\sigma(\_-z)-\sigma(\_)\rangle
$$
to define 
$$\alpha:\Gamma\to
L^\infty(G/N\times \widehat G/N^\perp,\U(1))$$
to be the 
$\U(1)$-valued perturbation in
\begin{eqnarray}
\overline{\lambda}(a,z)&=&
{\overline\kappa^\sigma}(z+\chi(a),\hat z+\hat\chi(a))^L\
u(z+\chi(a),\hat z+\hat\chi(a))\nonumber\\
&&\overline{\hat\lambda}(a,\hat z)\ u(z,\hat z)^{-1}\
{\overline\kappa^\sigma(z,\hat z)}^{-L}\ \alpha(a,z,\hat z).
\label{EqTheDefiOfAlpha}
\end{eqnarray}
By direct computation, it follows that 
\begin{equation}\label{EqAlphaIstWieIso}
(\delta_{\chi\times\hat \chi}\alpha)(b,a,z,\hat z)=
\psi(b,a,z)\hat\psi(b,a,\hat z)^{-1},
\end{equation}
for $\psi$ and $\hat\psi$ as above.

The next lemma is the defining lemma of $w$.
 
\begin{lem}
There is a continuous function 
$
w:\Gamma\times G/N\times \widehat G/N^\perp\to \U(1)
$
such that $$\alpha = \beta^L\ w,$$
where $\beta$ is from (\ref{EqDieDickeBerta}).
\end{lem}

\begin{pf} 
Equation (\ref{EqTheDefiOfAlpha}) is equivalent to
\begin{eqnarray*}
\alpha(a,z,\hat z){\overline\kappa^\sigma(z,\hat z)}^{-L}{\overline\kappa^\sigma}(z+\chi(a),\hat z+
\hat\chi(a))^L
&=&
\kappa^\sigma(z,\hat z)^{-L}\big[\overline{\lambda}(a,z)\big] u(z,\hat z) \\
&& \overline{\hat\lambda}(a,\hat z)^{-1}\ u(z+\chi(a),\hat z+\hat\chi(a))^{-1} ,
\end{eqnarray*}
but we observe that the right hand side consists of continuous 
terms only. To explore the continuity properties of $\alpha$ 
we compute straightforwardly 
\begin{eqnarray*}
&&{\overline\kappa^\sigma(z,\hat z)}^{-1}{\overline\kappa^\sigma}(z+\chi(a),\hat z+\hat\chi(a))\\
&=&
\langle\hat\sigma(\hat z),-\sigma(\_-z)+\sigma(\_)\rangle
 \langle\hat\sigma(\hat z+\hat\chi(a)),\sigma(\_-z-\chi(a))-\sigma(\_)\rangle\\
&=&
 \langle\hat\sigma(\hat z+\hat\chi(a))-\hat\sigma(\hat z)-\hat\sigma(\hat\chi(a)),-\sigma(z+\chi(a))
\rangle
 \langle\hat\sigma(\hat z),-\sigma(\chi(a))\rangle\\
 &&
  \langle\hat\sigma(\hat\chi(a)),\sigma(\_-z-\chi(a))-\sigma(\_)\rangle
   \langle\hat\sigma(\hat z),\sigma(\chi(a)-\sigma(\_+\chi(a))-\sigma(\_)\rangle\\
&& \langle\hat\sigma(\hat z),-\sigma(\_-z)+\sigma(\_)+\sigma(\_-z-\chi(a))-\sigma(\_-\chi(a))
\rangle,
\end{eqnarray*}
and we see that the factor in the first line is precisely $\beta(a,z,\hat z)^{-1}\in\U(1)$,
and the factors in the second and third line are all continuous functions
$\Gamma\times G/N\times \widehat G/N^\perp\to L^\infty(G/N,\U(1))$.
Therefore it follows that the product map 
$(a,z,\hat z)\mapsto \alpha(a,z,\hat z)\beta(a,z,\hat z)^{-L}$
is continous.
\end{pf}

Note that $w$ of the previous lemma defines
a class $[w]\in H^1(\Gamma,M_{\chi\times\hat\chi}/M_\chi M_{\hat\chi} )$.
In fact, we have $\delta_{\chi\times \hat\chi}w(a,b,z,\hat z)\in M_\chi M_{\hat\chi} $ as
\begin{eqnarray*}
\delta_{\chi\times \hat\chi}w(a,b,z,\hat z)&=&\delta_{\chi\times \hat\chi}\alpha(a,b,z,\hat z)\ 
\delta_{\chi\times \hat\chi}\beta^{-L}(a,b,z,\hat z)\\
&=&\psi(a,b,z)\hat\psi(a,b,\hat z)^{-1}\
(\chi \underset{\sigma,\hat\sigma}{\sqcup}\hat\chi)(a,b,z,\hat z)^{-L}
\end{eqnarray*}
This equation also states that
$B([w])+L\cdot\chi\sqcup\hat\chi=[\psi]\ominus[\hat\psi]$.
Thus we have defined a class
\begin{eqnarray*}
([w],[\psi],[\hat\psi])\in T_L(\Gamma,\chi,\hat\chi).
\end{eqnarray*}
It is clear that this class is not changed if the topological triple 
is replaced by a stable isomorphic one.
However, if we replace the triple by an equivalent one,
we see that the cocycles $\psi,\hat\psi$ remain unchanged
whereas $\alpha$, hence $w$, may be changed by 1-cocycle
$\Gamma\to M_{\chi\times\hat\chi}$. 
This means that the construction of the three classes leads us to a 
well-defined map
\begin{eqnarray}\label{EqThisIsTheCMForJingleBells}
{\rm Top_L}(\Gamma,\chi,\hat\chi)\to T_L(\Gamma,\chi,\hat\chi)/H^1(\Gamma,M_{\chi\times\hat\chi}),
\end{eqnarray}
and we show how to construct an inverse.

Let $([w],[\psi],[\hat\psi])\in T_L(\Gamma,\chi,\hat\chi)/H^1(\Gamma,M_{\chi\times\hat\chi})$
be represented by the $\U(1)$-valued functions
\begin{eqnarray*}
w:\Gamma\to M_{\chi\times\hat\chi},\quad
\psi:\Gamma^2\to M_\chi,\quad
\hat\psi:\Gamma^2\to M_{\hat\chi} .
\end{eqnarray*}
By re-normalising $\psi$ and $\hat\psi$ 
we can assume without restriction that 
$\delta_{\chi\times \hat\chi}w\ \delta_{\chi\times \hat\chi}\beta^L=\psi\ \hat\psi^{-1}$.
We make the following definitions
\begin{eqnarray*}
\lambda_\psi(a,z)&:=&{\rm Ad}\big( \psi(a, : ,z)\rho(a) \big)\quad \in \PU(L^2(\Gamma)),\\
\hat\lambda_{\hat\psi}(a,\hat z)&:=&{\rm Ad}\big( \hat\psi(a, : ,\hat z)\rho(a) \big) 
\quad\in \PU(L^2(\Gamma)),\\
\kappa_{w}(z,\hat z)&:=&{\rm Ad}\Big
(\overline\kappa^\sigma(z+\chi( : ),\hat z+\hat\chi( : ))^L\
\beta( : ,z,\hat z)^L \ w( : ,z,\hat z)\Big)\\
&&\in \PU(L^2(\Gamma)\otimes L^2(G/N)),
\end{eqnarray*}
where $\rho:\Gamma\to \U(L^2(\Gamma))$ is the right regular representation.
We claim that $(\kappa_w,(\chi,\lambda_{\psi}),(\hat\chi,\hat\lambda_{\hat\psi}))$
is a topological triple which is mapped to $([w],[\psi],[\hat\psi])$ under 
(\ref{EqThisIsTheCMForJingleBells}).

The first thing to note is that $\kappa_w$ is continuous and satisfies 
the Poincaré condition. In fact, we have
\begin{eqnarray*}
\kappa_{w}(z,\hat z)&:=&{\rm Ad}\Big
(\overline\kappa^\sigma(z+\chi( : ),\hat z+\hat\chi( : ))^L\
\beta( : ,z,\hat z)^L \ w( : ,z,\hat z)\Big)\\
&=&
{\rm Ad}\Big(\langle\hat\sigma(\hat z+\hat\chi(:)),\sigma(\_-z+\chi(:))-\sigma(\_)\rangle^L\\
&&
\langle \hat\sigma(\hat z+\hat\chi(:))-\hat\sigma(\hat\chi(:)),
\sigma(z)+\sigma(\chi(:))\rangle ^L\\
&&w(:,z,\hat z)\Big)\\
&=&
{\rm Ad}\Big(\langle\hat\sigma(\hat z),\sigma(\_-z)-\sigma(\_)\rangle^L\\
&&
\langle \hat z, \sigma(\_-z+\chi(:))-\sigma(\chi(:))-\sigma(\_-z)\rangle ^L\\
&&\langle \hat\sigma(\hat\chi(:)),
\sigma(\_-z+\chi(:))-\sigma(\_)\rangle ^L\\
&&w(:,z,\hat z)\Big),\\
\end{eqnarray*}
where the first factor guaranties the Poincaré condition 
(cp. Lemma \ref{LemTheDefiOfThePoincareClass}), 
and the remaining factors are already  continuous as unitary functions.

The cocycle conditions for $\lambda_{\psi},\hat\lambda_{\hat\psi}$
follow from the cocycle conditions for $\psi,\hat\psi$ which also 
imply that the classes defined in (\ref{EqTheClassesOTC})
are given by $[\psi],[\hat\psi]$. In fact, 
$(\delta_{\chi}\psi)(a,b,:,z)=1$ is equivalent to
\begin{eqnarray*}
\underbrace{\psi(b,:,z+\chi(a))\rho(b)}\ 
\underbrace{\psi(a,:,z)\rho(a)}& =&
\underbrace{\psi(ba,:,z)\rho(ba)}\ 
\underbrace{\psi(a,b,z)}\in \U(L^2(\Gamma)),\\ 
\overline\lambda_{\psi}(b,z+\chi(a))
\quad \ \overline\lambda_{\psi}(a,z)\quad
&&\qquad\overline\lambda_{\psi}(ba,z)\qquad
\U(1)
\end{eqnarray*}
and the application of ${\rm Ad}:\U(L^2(\Gamma))\to\PU(L^2(\Gamma))$
on both sides yields the cocycle condition.

By similar means, the equality 
\begin{eqnarray}\label{EqThisWillGiveTheDeckerForRET}
(\delta_{\chi\times\hat\chi}w)(a,:,z,\hat z)\ 
(\delta_{\chi\times\hat\chi}\beta^L)(a,:,z,\hat z)
=\psi(a,:,z)\ \hat\psi(a,:,\hat z)^{-1}
\end{eqnarray}
leads directly to the equality 
\begin{eqnarray*}
\lambda_{\psi}(a,z)\ \kappa_w(z,\hat z)= 
\kappa_w(z+\chi(a),\hat z+\hat\chi(a))\ \hat\lambda_{\hat\psi}(a,\hat z)\\
\in\PU(L^2(\Gamma)\otimes L^2(G/N))
\end{eqnarray*}
and also reproduces 
that $\beta^L\ w$ 
is the perturbation $\alpha$ in (\ref{EqTheDefiOfAlpha}).

This proves that
$(\kappa_w,(\chi,\lambda_{\psi}),(\hat\chi,\hat\lambda_{\hat\psi}))$
is a topological triple and is mapped to $([w],[\psi],[\hat\psi])$ by
 (\ref{EqThisIsTheCMForJingleBells}).
If we choose another representative of
$([w],[\psi],[\hat\psi])$, the equivalence  class of the constructed 
triple is not changed. Thus, our construction 
defines a  map
\begin{eqnarray}\label{EqThisIsTheCMForHellsBells}
T_L(\Gamma,\chi,\hat\chi)/H^1(\Gamma,M_{\chi\times\hat\chi})\to {\rm Top}_L(\Gamma,\chi,\hat\chi)
\end{eqnarray}
such that the composition with (\ref{EqThisIsTheCMForJingleBells}) 
is the identity on $T_L(\Gamma,\chi,\hat\chi)$.

It remains to check that the composition of 
(\ref{EqThisIsTheCMForJingleBells}) with 
(\ref{EqThisIsTheCMForHellsBells})
is the identity on ${\rm Top}_L(\Gamma,\chi,\hat\chi)$.
So start with a triple of order $L$
$({\kappa^\sigma}^L{\rm Ad}(u),(\chi,\lambda),(\hat\chi,\hat\lambda))$,
then take its class $([w],[\psi],[\hat\psi])$
and construct $(\kappa_w,(\chi,\lambda_{\psi}),(\hat\chi,\hat\lambda_{\hat\psi}))$
as above. We claim that the two triples are stably equivalent.
In fact, 
$({\kappa^\sigma}^L{\rm Ad}(u),(\chi,\lambda),(\hat\chi,\hat\lambda))$,
is stably equivalent to 
$({\kappa^\sigma}^L{\rm Ad}(u)\otimes \Eins,(\chi,\lambda\otimes {\rm Ad}(\rho)),
(\hat\chi,\hat\lambda\otimes{\rm Ad}(\rho)))$, but
the definitions of $\psi$ and  $\hat\psi$ lead to
\begin{eqnarray*}
\lambda(a,z){\rm Ad}(\rho(a))&=& \underbrace{{\rm Ad}(\overline\lambda(:,z+\chi(a)))^{-1}}\ 
\underbrace{{\rm Ad}( \psi(a,:,z)\rho(a))}\
\underbrace{{\rm Ad}(\overline\lambda(:,z))},\\
&&\qquad\qquad =:\theta(z+\chi(a))^{-1}\ \quad =\lambda_{\psi}(a,z)\qquad=:\theta(z)
\end{eqnarray*}
and
\begin{eqnarray*}
\hat\lambda(a,z){\rm Ad}(\rho(a))&=& 
\underbrace{{\rm Ad}(\overline{\hat\lambda}(:,\hat z+\hat\chi(a)))^{-1}}\ 
\underbrace{{\rm Ad}( \hat\psi(a,:,\hat z)\rho(a))}\
\underbrace{{\rm Ad}(\overline{\hat\lambda}(:,\hat z))}.\\
&&\qquad\qquad =:\hat\theta(\hat z+\hat\chi(a))^{-1}\ \quad 
=\hat\lambda_{\hat\psi}(a,z)\qquad=:\hat\theta(z)
\end{eqnarray*}
Thus the latter triple 
is isomorphic to 
$(\theta ({\kappa^\sigma}^L{\rm Ad}(u))\hat\theta^{-1},(\chi,\lambda_{\psi}),(\hat\chi,\hat\lambda_{\hat\psi}))$.
Finally, the isomorphism
$\theta ({\kappa^\sigma}^L{\rm Ad}(u))\hat\theta^{-1}$
can be computed using equation (\ref{EqTheDefiOfAlpha})
which leads to
\begin{eqnarray*}
\theta(z) {\kappa^\sigma}(z,\hat z)^L{\rm Ad}(u(z,\hat z)) \hat\theta^{-1}(\hat z)
&=&
\kappa_w(z,\hat z){\rm Ad}\Big(
u(z+\chi(:),\hat z+\chi(:))\Big).
\end{eqnarray*}
This shows that the triple
$(\theta ({\kappa^\sigma}^L{\rm Ad}(u))\hat\theta^{-1},(\chi,\lambda_{\psi}),(\hat\chi,\hat\lambda_{\hat\psi}))$
is equivalent to 
$(\kappa_w,(\chi,\lambda_{\psi}),(\hat\chi,\hat\lambda_{\hat\psi}))$.
\\

This proves Theorem \ref{TheBigClThmForEquTrip}.

\section{Regular Topological Triples}
\label{SecRegET}

\noindent
Having the classification of the previous section at hand 
we introduce the notion of regularity of a topological triple.
The motivation for considering regular topological triples 
comes from the fact that the underlying pair 
and  the underlying dual pair of a regular triple (of order 1)
admit natural extensions to dynamical triples in duality.

We will see in section \ref{SecAAppliEquKT} that 
in the important case of $G=\R^n$, $N=\Z^n$
all topological triples are regular.

\begin{defi}\label{DefiRegularity}
A topological triple of order $L$ is called {\bf regular} 
if the first component of its class 
$([w],[\psi],[\hat\psi])\in T_L(\Gamma,\chi,\hat\chi)/H^1(\Gamma,M_{\chi\times\hat\chi})$ 
is zero:
\begin{eqnarray*}
[w]=0\in H^1(\Gamma,M_{\chi\times\hat\chi}/M_\chi M_{\hat\chi} )/H^1(\Gamma,M_{\chi\times\hat\chi}).
\end{eqnarray*}
\end{defi}

For a finite group $\Gamma$ 
we denote by 
\begin{eqnarray*}
{\rm Top}^\dag_L(\Gamma)\subset
{\rm Top}_L(\Gamma)
\end{eqnarray*}
 the set of all (equivalence classes of) 
regular topological triples, and 
similarly, 
\begin{eqnarray*}
{\rm Top}^\dag_L(\Gamma,\chi,\hat\chi)
\subset {\rm Top}_L(\Gamma,\chi,\hat\chi)
\end{eqnarray*}
is the set of all regular topological triples 
with fixed homomorphisms  $\chi, \hat\chi$. 

\begin{cor}
1.
Let $(\kappa,(\chi,\lambda),(\hat\chi,\hat\lambda))$ 
be a regular topological triple of order $L\in \Z$, then 
\begin{eqnarray*}
L\cdot (\chi\cup\hat\chi)=0\in H^3(\Gamma,\U(1)).
\end{eqnarray*}
\item
2. 
If, conversely, for any two homomorphisms 
$\chi$ and $\hat\chi$ we have that 
$L\cdot(\chi\cup\hat\chi)=0$, then there exists 
a regular topological triple of order $L$ extending 
these two homomorphisms.
\item 
3.
Moreover, the set of regular topological triples 
of order $L$, ${\rm Top}^\dag_L(\Gamma,\chi,\hat\chi)$, 
is  a $H^2(\Gamma,\U(1))$-torsor.
\end{cor}
\begin{pf}
As $C:\chi\sqcup\hat\chi\mapsto\chi\cup\hat\chi$ in (\ref{DiagTheClassiGroupSitsInAha}),
all three statements of the corollary follow 
immediately from the exactness of the vertical 
sequence in (\ref{DiagTheClassiGroupSitsInAha})
and the classification of topological triples.
\end{pf}

Now, take a regular topological triple in form of  
a class $([w],[\psi],[\hat\psi])\in T_L(\Gamma,\chi,\hat\chi)$.
Choose representatives $(w,\psi, \hat\psi)$ satisfying
equation (\ref{EqThisWillGiveTheDeckerForRET}). 
The condition $[w]=0\in H^1(\Gamma,M_{\chi\times\hat\chi}/M_\chi M_{\hat\chi} )/H^1(\Gamma,M_{\chi\times\hat\chi})$
means for the function $w:\Gamma\to M_{\chi\times\hat\chi}$ that
it has a decomposition 
\begin{eqnarray*}
w(a,z,\hat z)= v(a,z,\hat z)\ 
v_1(z)\ v_2(\hat z),
\end{eqnarray*}
where $\delta_{\chi\times\hat\chi}v=1.$
Thus, if we re-normalise $\psi$ to  $\psi(a,b,z)v_1(z+\chi(a))^{-1}$
and $\hat\psi$ to  $\hat \psi(a,b,\hat z)v_2(\hat z+\hat\chi(a))$, then
equation (\ref{EqThisWillGiveTheDeckerForRET})
 turns into
\begin{eqnarray}\label{EqDiHirObe}
(\chi \underset{\sigma,\hat\sigma}{\sqcup}\hat\chi)(a,b,z,\hat z)^{L}=\psi(a,b,z)\hat\psi(a,b,\hat z)^{-1}.
\end{eqnarray}
This equality together with  the structure of $\chi \underset{\sigma,\hat\sigma}{\sqcup}\hat\chi$
as stated  in  Lemma \ref{LemThePropertiesOfDeltaBeta}
give rise to a natural map
\begin{eqnarray}\label{EqFTopLtoTop1}
\hat L_*: {\rm Top}^\dag_{L}(\Gamma, \chi,\hat\chi)
&\to& {\rm Top}^\dag_1(\Gamma,\chi,L\hat\chi)
\end{eqnarray}
which we describe next.
Putting $z=0$ in (\ref{EqDiHirObe})
we see that
\begin{eqnarray*}
\hat\psi(a,b,\hat z)&=&\gamma(a,b)^{-L}\langle L\hat z,\eta(a,b)\rangle^{-1}\psi(a,b,0)
\end{eqnarray*}
Thus we can  define 
\begin{eqnarray*}
(\hat L_*\hat\psi)(a,b,\hat z)&:=&\gamma(a,b)^{-L}\langle \hat z,\eta(a,b)\rangle^{-1}\psi(a,b,0).
\end{eqnarray*}
It follows from $\delta_{\hat\chi}\hat\psi=1$ that
$\delta_{L\hat\chi}(\hat L_*\psi)=1$.
We leave it to the reader to verify the identity
\begin{eqnarray}
\label{EqOhoTabako}
\\
\psi(a,b,z) (\hat L_*\hat\psi)(a,b,\hat z)^{-1}
(\delta_{\chi}\nu)(a,b,z)
=(\chi \underset{\sigma,\hat\sigma}{\sqcup}L\hat\chi)(a,b,z,\hat z), 
\nonumber
\end{eqnarray}
wherein $\nu(a,z):=\langle L\hat\sigma(\hat\chi(a))-\hat\sigma(L\hat\chi(a)),z+\chi(a)\rangle$
and $\chi \underset{\sigma,\hat\sigma}{\sqcup}L\hat\chi$ is as in 
Lemma \ref{LemThePropertiesOfDeltaBeta}
but with $\hat\chi$ replaced by $L\hat\chi$.
Taking the cohomology on both sides of (\ref{EqOhoTabako})
we find
$$
[\psi]\ominus [\hat L_*\hat\psi]=\chi\sqcup(L\hat\chi)\in H^2(\Gamma,M_\chi M_{L\hat\chi}).
$$
In other words $(1,[\psi],[\hat L_*\hat\psi])$ is the class
of a topological triple of order 1.
This defines (\ref{EqFTopLtoTop1}).

Similarly, one defines a natural map
\begin{eqnarray}\label{EqFTopLtoTop1222}
 L_*: {\rm Top}^\dag_{L}(\Gamma, \chi,\hat\chi)
&\to& {\rm Top}^\dag_1(\Gamma,L\chi,\hat\chi),
\end{eqnarray}
and it is straight forward to check that
the maps defined are bijections whose inverses 
are the restrictions of the two maps in (\ref{EqFromTOp1ToTopL}) 
to regular topological triples.
In fact, one can easily generalise this construction to
give an inverse of the restriction of (\ref{EqThereAreNMF!TOL})
to regular triples.

We give  another application of (\ref{EqDiHirObe}).
Namely, applying therein the boundary operator $d$ we get
\begin{eqnarray}
(d\psi)(a,b,g,z)&=&\psi(a,b,gN+z)\psi(a,b,z)^{-1}\nonumber\\
&=&\langle \hat\eta(a,b),gN\rangle^{-L}\nonumber\\
&=&\phi(b,g,z+\chi(a))\phi(ba,g,z)^{-1}\phi(a,g,z)\nonumber\\
&=&(\delta_{\chi}\phi)(a,b,g,z),
	\label{EqDPsiEqDelatZeta}
\end{eqnarray}
where $\phi(a,g,z):=\langle\hat\sigma(\hat\chi(a)),g\rangle^{-L}$, and 
if we define $\hat\phi(a,\hat g,\hat z):=\langle \hat g ,\sigma(\chi(a))\rangle^{-L}$,
we obtain similarly
\begin{eqnarray}\label{EqGesternEinhandmeister}
(\hat d\hat \psi)(a,b,\hat g,\hat z)
&=&(\delta_{\hat\chi}\hat\phi)(a,b,\hat g,\hat z).\label{EqDPsiEqDelatZetaDual}
\end{eqnarray}
These two equations determine two classes
\begin{eqnarray*}
[\psi,\phi,1]\in H^{2,\dag}(\Gamma,G,M_\chi), \quad
[\hat\psi,\hat\phi,1]\in H^{2,\dag}(\Gamma,\widehat G,M_{\hat\chi}),
\end{eqnarray*}
and hence, by the classification of dynamical triples,
we obtain two maps
\begin{eqnarray}\label{EqFromTopLToDynnn}
\xymatrix{
&{\rm Top}^\dag_L(\Gamma)\ar[ld]_{\epsilon_L(\Gamma)}\ar[dr]^{\hat\epsilon_L(\Gamma)}&\\
{\rm Dyn}^\dag(\Gamma)&&\widehat{\rm Dyn}^\dag(\Gamma)
}
\end{eqnarray}
which are  natural in $\Gamma$.

The maps of (\ref{EqFromTopLToDynnn}),
(\ref{EqFTopLtoTop1}) and (\ref{EqFTopLtoTop1222})
behave well to each other:
\begin{lem}
The diagrams 
\begin{eqnarray*}
\xymatrix{
&{\rm Top}^\dag_L(\Gamma,\chi,\hat\chi)\ar[ld]\ar[d]&
{\rm Top}^\dag_L(\Gamma,\chi,\hat\chi)\ar[rd]\ar[d]&\\
{\rm Dyn}^\dag(\Gamma)&{\rm Top}_1^\dag(\Gamma,\chi,L\hat\chi),\ar[l]&
{\rm Top}_1^\dag(\Gamma,L\chi,\hat\chi)\ar[r]&\widehat{\rm Dyn}^\dag(\Gamma)
}
\end{eqnarray*}
commute.
\end{lem}

\begin{pf}
To prove the commutativity of the first diagram note that
the boundary operator $d$ applied to  
(\ref{EqOhoTabako}) leads to
\begin{eqnarray*}
d(\psi\delta_\chi\nu)(a,b,g,z)&=&\langle\hat\eta_L(a,b),gN\rangle^{-1}\\
&=&\langle\hat\eta(a,b),gN\rangle^{-L} (d\delta_\chi\nu)(a,g,z)\\
&=&\delta_\chi(\phi d\nu)(a,b,g,z),
\end{eqnarray*}
with $\phi$ as in (\ref{EqDPsiEqDelatZeta}).
So the cocycle $(\psi\delta_\chi\nu,\phi d\nu,1)=(\psi,\phi,1)\cdot\partial_\chi(\nu,1)$
which represents the composition in the 
first diagram is cohomologous to
$(\psi,\phi,1)$ which represents the diagonal map.

The commutativity of the second 
diagram is proven by an analogous argument.
\end{pf}

If we stick the constructed maps for different homomorphisms
$\chi,\hat\chi$ together, then we obtain
commutative diagrams
\begin{eqnarray*}
\xymatrix{
&{\rm Top}^\dag_L(\Gamma)\ar[ld]_{\epsilon_L(\Gamma)}\ar[d]^{\hat L_*}\ar@/^1,5cm/[dd]&&
{\rm Top}^\dag_L(\Gamma)\ar[rd]^{\hat\epsilon_L(\Gamma)}\ar[d]_{ L_*}\ar@/_1,5cm/[dd]&\\
{\rm Dyn}^\dag(\Gamma)\ar@/_0,5cm/[dr]&
{\rm Top}_1^\dag(\Gamma)\ar[d]\ar[l]^{\ \ \epsilon_1(\Gamma)}&&
{\rm Top}_1^\dag(\Gamma)\ar[d]\ar[r]_{\hat\epsilon_1(\Gamma)\ \ }&
\widehat{\rm Dyn}^\dag(\Gamma),\ar@/^0,5cm/[dl]\\
&{\rm Par}(\Gamma)&&\widehat{\rm Par}(\Gamma)&
}
\end{eqnarray*}
wherein the maps ending at ${\rm Par}(\Gamma)$ or at 
$\widehat{\rm Par}(\Gamma)$
are the natural forgetful maps.
Note that in contrast to (\ref{EqFTopLtoTop1}) and (\ref{EqFTopLtoTop1222})
the assembled maps $L_*, \hat L_*:{\rm Top}_L^\dag(\Gamma)\to {\rm Top}_1^\dag(\Gamma)$  
are no isomorphisms in general, as they fail to be surjective.

In the next section we show that 
$\epsilon_1(\Gamma)$ and $\hat\epsilon_1(\Gamma)$
are bijections and that the composition  
$\epsilon_1(\Gamma)^{-1}\circ\hat\epsilon_1(\Gamma)$
is the duality map $\ \widehat{ }\ $
introduced in (\ref{EqTheDualityMap}).

\section{From Dynamical to Topological Triples}
\label{SecFromDToTT}

\noindent
Consider a dualisable dynamical triple which is represented by
a class $[\psi,\phi,1] \in H^{2,\dag}_{\rm cont}(\Gamma,G,M_\chi)$.
Let $[\hat\psi,\hat\phi,1]\in H^{2,\dag}_{\rm cont}(\Gamma, \widehat G,M_{\hat\chi})$
be its dual according to the duality map (\ref{EqTheDualityMap}).
By forgetting $\phi$ and $\hat\phi$ we obtain a class
$([\psi],[\hat\psi])\in H^2(\Gamma,M_\chi)\oplus H^2(\Gamma,M_{\hat\chi})$.

\begin{prop}\label{PropNichtPop}
In the Mayer-Vietoris sequence (\ref{DiagTwoLongExSeq}) the class
$([\psi],[\hat\psi])$ is mapped  to $\chi\sqcup\hat\chi$ by 
the subtraction map $\ominus$.
\end{prop}

\begin{pf}
The class $[\psi]\ominus[\hat\psi]$ is represented 
by $\psi\hat\psi^{-1}$ which we can compute using 
 the definition of $\hat\psi$ in (\ref{EqTheDeOTDPsi}).
Let us denote by $\sigma:G/N\to G$ and $\hat\sigma:\widehat G/N^\perp\to \widehat G$ Borel sections 
of the quotient maps $G\to G/N,\widehat G\to\widehat G/N^\perp$.
We also let $\varphi(a,z):=\phi(a,\sigma(z),0)\langle\hat\sigma(\hat\chi(a)),\sigma(z)\rangle$
which is a well-defined continuous function due to the definition of $\hat\chi$.
We have
\begin{eqnarray*}
\psi(a,b,z)\hat\psi(a,b,\hat z)^{-1}&=&
\psi(a,b,z)\psi(a,b,0)^{-1} \phi(b,\sigma(\chi(a)),0)\\
&&
\langle\hat\chi(ba)+\hat z,\hat\sigma(\chi(b))+\hat\sigma(\chi(a))-\hat\sigma(\chi(ba))\rangle.
\end{eqnarray*}
Using the cocycle conditions $d\psi=\delta_\chi\phi, d\phi=1$ and the definition of $\varphi$
this can finally be transformed to
\begin{eqnarray*}
\psi(a,b,z)\hat\psi(a,b,\hat z)^{-1}&=&
(\delta_\chi\varphi)(a,b,z)\ (\delta_\chi s)(a,b)\\
&&
\langle\hat\sigma(\hat\chi(a)),\sigma(\chi(b))\rangle\\
&&
\langle\hat\sigma(\hat\chi(ba))-\hat\sigma(\hat\chi(a))
-\hat\sigma(\hat\chi(b)),\sigma(\chi(ba))\rangle\\
&&
\langle\hat\sigma(\hat\chi(ba))-\hat\sigma(\hat\chi(a))
-\hat\sigma(\hat\chi(b)),\sigma(z)\rangle\\
&&
\langle\hat z,\sigma(\chi(b))+\sigma(\chi(a))-\sigma(\chi(ba))\rangle,
\end{eqnarray*}
where $s(a):=\langle\hat\sigma(\hat\chi(a)),\sigma(\chi(a))\rangle$.
If we compare this with the results of Lemma 
\ref{LemThePropertiesOfDeltaBeta},
the proposition follows.
\end{pf}

By this proposition, the class 
$([\psi],[\hat\psi])\in H^2(\Gamma,M_\chi)\oplus H^2(\Gamma,M_{\hat\chi})$
may also be considered as a class
$(1,[\psi],[\hat\psi])\in T_1(\Gamma,\chi,\hat\chi)/ H^1(\Gamma,M_{\chi\times\hat\chi})$,
and therefore we obtain (the class of) a regular toplological triple by the classification 
of topological triples.
Hence, the assignment 
$(\chi,[\psi,\phi,1])\mapsto (1,[\psi],[\hat\psi])$
defines  natural maps 
\begin{eqnarray}
\delta(\Gamma,\chi):{\rm Dyn}^\dag(\Gamma,\chi)&\to&{\rm Top}^ \dag_1(\Gamma,\chi),\\
\delta(\Gamma):{\rm Dyn}^\dag(\Gamma)&\to&{\rm Top}^ \dag_1(\Gamma).\nonumber
\end{eqnarray}

In the next three lemmata we point out the 
relation between $\delta(\Gamma,\chi)$ and the maps 
$\epsilon_1(\Gamma,\chi), \hat\epsilon_1(\Gamma,\chi)$
 constructed in the previous section.

\begin{lem}
The composition $\epsilon_1(\Gamma,\chi)\circ\delta(\Gamma,\chi)$
is the identity on ${\rm Dyn}^\dag(\Gamma,\chi)$.
\end{lem}

\begin{pf}
Start with a cocycle $(\chi, (\psi,\phi,1))$ representing a dynamical triple. 
In the course of the proof of Proposition \ref{PropNichtPop} 
we calculated $\psi(a,b,z)\hat\psi(a,b,\hat z)^{-1}.$
If we apply here the boundary operator $d$, we get
after some steps
\begin{eqnarray*}
d(\psi\delta_\chi\varphi^{-1})(a,b,g,z)&=&\langle\hat\sigma(\hat\chi(b))-\hat\sigma(\hat\chi(ba))
-\hat\sigma(\hat\chi(a)),gN\rangle^{-1}\\
&=&\delta_\chi(\phi d\varphi^{-1})(a,b,g,z).
\end{eqnarray*}
Thus   the cocycle $(\psi\delta_\chi\varphi^{-1},\phi d\varphi^{-1},1)=
(\psi,\phi,1)\cdot\partial_\chi(\varphi^{-1},1)$
which represents the composition is cohomologous to 
$(\psi,\phi,1)$.
\end{pf}

\begin{lem}
The composition $\delta(\Gamma,\chi)\circ\epsilon_1(\Gamma,\chi)$
is the identity on ${\rm Top}_1^\dag(\Gamma,\chi)$.
\end{lem}

\begin{pf}
Let $(\psi,\hat\psi)$ satisfy (\ref{EqDiHirObe}) for $L=1$, taking 
therein $ z=0$ we find
\begin{eqnarray*}
\hat\psi(a,b,\hat z)=\gamma(a,b)^{-1}\langle \hat z,\eta(a,b)\rangle^{-1}\psi(a,b,0).
\end{eqnarray*}
Now, let $\phi$ be defined according to (\ref {EqDPsiEqDelatZeta}),
and then let $\hat\psi'$ be defined by (\ref{EqTheDeOTDPsi}),
so
\begin{eqnarray*}
\hat\psi'(a,b,\hat z)&=&
\psi(a,b,0) \langle\hat\sigma(\hat\chi(b)),\sigma(\chi(a))\rangle\\
&&\langle\hat\chi(ba)+\hat z,\sigma(\chi(ba))-\sigma(\chi(b))-\sigma(\chi(a))\rangle.
\end{eqnarray*}
This may be transformed to 
\begin{eqnarray*}
\hat\psi'(a,b,\hat z)
=\hat\psi(a,b,\hat z) (\delta_\chi\nu)(a,b,z)
\end{eqnarray*}
where $\nu(a,z):=\langle\hat\sigma(\hat\chi(a)),\sigma(\chi(a))\rangle^{-1}$ (independent of $z$).
Thus, we have the equality  $[\psi,\hat\psi]=[\psi,\hat\psi']\in H^2(\Gamma,M_\chi M_{\hat\chi})$
which proves the lemma. 
\end{pf}

We have just proven that $\delta(\Gamma,\chi)$ is the 
inverse of $\epsilon_1(\Gamma,\chi)$.
Even more is true:

\begin{lem}
The composition $\hat\epsilon_1(\Gamma)\circ\delta(\Gamma)$
is the duality map $\ \widehat{ }\ $
introduced in (\ref{EqTheDualityMap}).
\end{lem}

\begin{pf}
Let $(\chi,(\psi,\phi,1))$ represent a dualisable dynamical triple.
Its dual is represented by $(\hat\chi,(\hat\psi,\hat\phi,1)$ according 
to (\ref{EqTheDualTorusCoycle123}), (\ref{EqTheDeOTDPsi}) and 
(\ref{EqTDOTDPhi}).

Now, forget $\hat\phi$, and  compute $\psi\hat\psi^{-1}$
as in the course of the proof of Proposition \ref{PropNichtPop}.
Then apply  the boundary operator $\hat d$
which defines  $\hat\phi'$  according to (\ref{EqGesternEinhandmeister}).
We even have $\hat\phi'=\hat\phi$, so
the cocycles $(\hat\psi,\hat\phi',1)$ and 
$(\hat\psi,\hat\phi,1)$ which represent 
the composition $\hat\epsilon_1(\Gamma)\circ\delta(\Gamma)$
and the duality map $\ \widehat{ }\ $ coincide.
\end{pf}

Dually to the construction of $\delta(\Gamma)$, 
we have a natural map
\begin{eqnarray*}
\hat\delta(\Gamma):\widehat{\rm Dyn}^\dag(\Gamma)\to {\rm Top}^\dag(\Gamma)
\end{eqnarray*}
with analogous properties. 

To summarise this section, 
we have a commutative 
diagram of natural  isomorphisms
\begin{eqnarray*}
\xymatrix{
&{\rm Top}^\dag_1(\Gamma)&\\
{\rm Dyn}^\dag(\Gamma)\ar[rr]^{\widehat{\ }}\ar[ru]^{\delta(\Gamma)=\epsilon(\Gamma)^{-1}}&
&\ar[ul]_{\hat\delta(\Gamma)=\hat\epsilon(\Gamma)^{-1}}\widehat{\rm Dyn}^\dag(\Gamma).
}
\end{eqnarray*}

\section{Topological Triples and  K-Theory}
\label{SecAAppliEquKT}

\noindent
In this section we fix a natural number
$n\in \{1,2,3,\dots\}$ and we let
\begin{eqnarray*}
G:=\R^n,\quad N:=\Z^n.
\end{eqnarray*}
We identify the dual group of $\R$ with $2\pi\R$ 
containing the dual lattice $2\pi\Z$.
The $n$-torus is $\T^n:=\R^n/\Z^n$, and 
the dual $n$-torus is $\hat\T^n:=2\pi\R/2\pi\Z$
which should not be mixed up with the dual 
group $\widehat\T^n= 2\pi\Z$ of $\T^n$.

Recall the classification of topological 
triples from section \ref{SecTheClassifOfET}.
There we introduced the $\Gamma$-module
$M:= C(\T^n\times\hat\T^n,\U(1))$ 
with its submodules
$M_\chi:= C(\T^n,\U(1))$ 
and
$M_{\hat\chi} := C(\hat\T^n,\U(1))$ 
on which $\Gamma$ acts 
with (the negative of) homomorphisms 
$\chi:\Gamma\to\T^n$
and $\hat\chi:\Gamma\to\hat\T^n$.
In the classification of topological triples 
the group $H^n(\Gamma,M_{\chi\times\hat\chi}/M_\chi M_{\hat\chi} )/H^1(\Gamma,M_{\chi\times\hat\chi})$
occured, and as we saw in section \ref{SecRegET},
this group has an interpretation as an obstruction group.

However, in the case of $G=\R^n, N=\Z^n$, 
the next lemma simplifies the discussion essentially.

\begin{lem}
If $G=\R^n$ and $N=\Z^n$, then
$H^k(\Gamma,M_{\chi\times\hat\chi}/M_\chi M_{\hat\chi} )=0$  for all $k>0$.
\end{lem}

\begin{pf}
Let $V:=C_{\rm null}(\T^n\times \hat\T^n,\R)$ denote the 
vector space  of real valued, null-homotopic and base point preserving functions 
on the product of the torus and the dual torus.
Let $V_1:=C_{\rm null}(\T^n,\R)$ and
$V_2:=C_{\rm null}(\hat\T^n,\R)$
be the subspaces of those functions of $V$ which 
only depend on one variable.
The group $V\times \U(1)$, respectively $V_1\times V_1\times\U(1)$, 
can be identified with the subgroup of $M_{\chi\times\hat\chi}$, respectively $M_\chi M_{\hat\chi} $, 
consisting of all null-homotopic functions.
The crucial point is that the respective quotients, 
the homotpy classes of maps,  are isomorphic, 
namely  
\begin{eqnarray*}
[\T^n,\U(1)]\times[\hat\T^n,\U(1)]&\stackrel{\cong}{\longrightarrow}&
[\T^n\times\hat\T^n,\U(1)].\\
([f],[g])&\longmapsto& [f(\_)g(..)]
\end{eqnarray*}
This leads us to a diagram of short exact sequences
of $\Gamma$-modules
\begin{eqnarray*}
\resizebox*{!}{5cm}{
\xymatrix{
&0\ar[d]&0\ar[d]&0\ar[d]&\\
0\ar[r]&V_1\times V_2\times\U(1)\ar[r]\ar[d]&
 M_\chi M_{\hat\chi} \ar[r]\ar[d]&[\T^n,\U(1)]\times[\hat\T^n,\U(1)]\ar[r]\ar[d]^\cong&0\\
0\ar[r]&V\times\U(1)\ar[r]\ar[d]& M_{\chi\times\hat\chi}\ar[r]\ar[d]&[\T^n\times\hat\T^n,\U(1)]\ar[r]\ar[d]&0\\
0\ar[r]&V/V_1\times V_2\ar[r]\ar[d]&M_{\chi\times\hat\chi}/ M_\chi M_{\hat\chi} \ar[r]\ar[d]& 0\ar[r]\ar[d]&0.\\
&0&0&0&
}
}
\end{eqnarray*}
It follows that $M_{\chi\times\hat\chi}/M_\chi M_{\hat\chi} $ is $\Gamma$-isomorphic to a vector space,
and this proves the lemma, because the cohomology of a finite group with 
values in a vector space always vanishes.  
\end{pf}

By definition of regularity 
(Def. \ref{DefiRegularity}),
we have:

\begin{cor} If $G=\R^n$ and $N=\Z^n$, then
every topological triple is regular.
\end{cor}

In particular, if we focus our attention to topological triples 
of order 1, then we have isomorphisms
\begin{eqnarray}\label{EqSchonWichtigAuchWennsAm}
{\rm Dyn}^\dag(\Gamma)\cong{\rm Top}_1(\Gamma)\cong\widehat{\rm Dyn}^\dag(\Gamma).
\end{eqnarray}

Now, let   $(\kappa,(E,P),(\widehat E,\widehat P))$  be a topological triple
of order $L$, and let us denote by $(L_*\kappa, (L_*E,L_*P),(\widehat E,\widehat P))$
(a representative of) the image of the triple  $(\kappa,(E,P),(\widehat E,\widehat P))$
under  $L_*$  (s. (\ref{EqFTopLtoTop1222})).
Let $C(L^*)$  denote  the mapping cone of  
the pullback $L^*:C^*(L_*E,L_*P)\to C^*(E,P)$  
(cf. (\ref{EqDochNochNeNumma}), p. \pageref{EqDochNochNeNumma}). 
As the topological triple 
$(L_*\kappa, L_*(E,P),(\widehat E,\widehat P))$ has order 1, 
we have by (\ref{EqSchonWichtigAuchWennsAm}) an isomorphism
$K^{*}_{\Gamma,L_*P}(L_*E)\cong K^{*+n}_{\Gamma,\widehat P}(\widehat E)$
as explained in (\ref{EqAlterDiePapersKlebenNicht}).
So if we replace 
$K^{*}_{\Gamma,L_*P}(L_*E)$ by $K^{*+n}_{\Gamma,\widehat P}(\widehat E)$
in the six term exact sequence (\ref{DiagMyFirstSixPack}), we find an 
exact sequence 
\begin{eqnarray*}
\xymatrix{
 &
K^{n}_{\Gamma,\widehat P}(\widehat E)\ar[r]&
K^{0}_{\Gamma, P}( E)\ar[r]& 
K_{1}(C(L^*))
\ar[d]&\\
&
K_{0}(C( L^*))\ar[u]&
K^{1}_{\Gamma,P}(E)\ar[l]&
K^{1+n}_{\Gamma,\widehat P}(\widehat E) .\ar[l]& 
}
\end{eqnarray*}
Dually, we also have an exact sequence 
\begin{eqnarray*}
\xymatrix{
 &
K^{n}_{\Gamma, P}( E)\ar[r]&
K^{0}_{\Gamma,\widehat P}(\widehat E)\ar[r]& 
K_{1}(C(\hat L^*))
\ar[d]&\\
&
K_{0}(C(\hat L^*))\ar[u]&
K^{1}_{\Gamma,\widehat P}(\widehat E)\ar[l]&
K^{1+n}_{\Gamma, P}( E),\ar[l]& 
}
\end{eqnarray*}
where $C(\hat L^*)$ is the mapping cone of the pullback
$\hat L^*:C^*(\hat L_*\widehat E,\hat L_*\widehat P))\to C^*(\widehat E,\widehat P)$.
Here $(L_*\kappa, (E,P),(\hat L_*\widehat E,\hat L_*\widehat P))$ is 
(a representative of) the image of   $(\kappa,(E,P),(\widehat E,\widehat P))$
under $\hat L_*$ (s. (\ref{EqFTopLtoTop1})).

In the special case of $L=1$ these sequences 
reduce to the isomorphism (\ref{EqAlterDiePapersKlebenNicht}), 
namely:

\begin{cor}\label{CorDickUndAuch}
Let $G=\R^n$ and $N=\Z^n$, and let $(\kappa,(E,P),(\widehat E,\widehat P))$
be a topological triple of order $1$. Then there is
an isomorphism in twisted equivariant K-theory 
\begin{eqnarray*}
K^*_{\Gamma,P}(E)\cong K^{*+n}_{\Gamma,\widehat P}(\widehat E).
\end{eqnarray*}
\end{cor}

\begin{pf}
For $L=1$ we have $1^*={\rm id}_{C^*(E,P)}$ and 
$\hat 1^*={\rm id}_{C^*(\widehat E,\widehat P)}$. 
Thus, the mapping cones $C(1^*)$ and $ C(\hat 1^*)$ 
above are contractible and their K-theories vanish.
\end{pf}

\section{Conclusion and Outlook}

\noindent
In this work we considered topological (T-duality) triples 
over the singleton space $\ast$ with an action
of a finite group $\Gamma$. In more learned terms, 
this is a topological triple over the proper étale groupoid 
$\Gamma \rightrightarrows \ast$.
Topological (T-duality) triples  over a space $B$ 
as as considered in \cite{BS,BSST,Sch} may be 
understood as topological triples over the proper étale groupoid 
$B\rightrightarrows B$.
In a future work we are going to combine the methods 
of \cite{Sch} and this work to deal with topological 
(T-duality) triples over any proper étale  groupoid 
$\Gamma_1\rightrightarrows \Gamma_0$, that is a 
topological (T-duality) triple  over $\Gamma_0$
in the sense of \cite{BS,BSST,Sch} equipped with
an action of $\Gamma_1$.
We then will  also derive the corresponding isomorphism in
K-theory (cf. Corollary \ref{CorDickUndAuch}), and,
after proving Morita invariance of the set-up, we also obtain the corresponding 
statement for topological (T-duality) triples over orbi\-folds, 
i.e. for stacks which have  presentations by proper étale groupoids.

\begin{appendix}

\section{A Mayer-Vietoris Sequence in Group Cohomology}

\label{SecTheMVSofGC}

\noindent
Let $(M,+)$ be an abelian group and $\Gamma$-module.
Assume that there are sub-modules $M_1,M_2\subset M$. 
Thus, there is a diagram of $\Gamma$-modules
\begin{eqnarray*}
\resizebox*{!}{2,5cm}{
\xymatrix{
&&M_1\ar@{^{(}->}[dr]&&\\
&M_1\cap M_2 \ \ \ \ \ \ \ \ \ar@{^{(}->}[ur]\ar@{^{(}->}[dr]&&M_1+ M_{2},& \\
&&M_{2} \ar@{^{(}->}[ur]&&
}}
\end{eqnarray*}
where $M_1+M_2$ is the sub-module of $M$ consisting of
all $m\in M$ such that there exist $x_i\in M_i$ satisfying $ m=x_1+x_2$.
These modules give rise to a long 
exact sequence in group cohomology
\begin{eqnarray*}
\xymatrix{
\cdots\ar[r] &
H^n(\Gamma,M_1\cap M_2)\ar[r]^{\hspace{-0.7cm}{\rm diag}}&
H^n(\Gamma,M_1)\oplus H^n(\Gamma, M_2)
\ar 
`[r] 
`[l]^\ominus 
`[lld]
`[r]
[dl]
&\\
&H^n(\Gamma,M_1+ M_2)\ar[r]^{\!\!\!\! C}&
H^{n+1}(\Gamma,M_1\cap M_2)\ar[r]^{\hspace{1cm}{\rm diag}}&\cdots .
}
\end{eqnarray*}
Namely, the diagonal map ${\rm diag}$ is defined by
${\rm diag}([\omega]):=([\omega],[\omega])$, and
the subtraction map $\ominus$ is defined 
by $[\omega_1]\ominus[\omega_2]:=[\omega_1-\omega_2]$.
The connecting homomorphism $C$ is given 
by $C([\Omega]):=[\delta \omega_1]=-[\delta\omega_2]$, where $\delta$
is the boundary operator of the cochain complex $C^\bullet(\Gamma,M)$, and 
$\Omega=\omega_1+\omega_2$ is any decompositon
such that $\omega_1:\Gamma^n\to M_1$ and 
$\omega_2:\Gamma^n\to M_2$.

It is left to the reader to  to check that these maps are well-defined
and  that the resulting sequence is exact.

\section{The Equivariant $\check{\rm C}$ech Class
$ [v]$}
\label{SecTECeCOBI12}

\noindent
For  a  homomorphisms $\chi:\Gamma\to G/N$ 
consider a $\Gamma$-equivariant covering $V_{\boldsymbol\cdot}:=\{V_i|i\in I\}$ of the space 
$G/N$.
This yields a double complex
\begin{eqnarray}\label{DigGrischssbljab12}
\xymatrix{
\vdots&\vdots&\vdots& \\
 C^{2,0}(\Gamma,V_{\cdot})\ar[r]^{\delta_{\chi}}\ar[u]^{\check\delta}&
 C^{2,1}(\Gamma,V_{\cdot})\ar[r]^{\delta_{\chi}}\ar[u]^{\check\delta}&
 C^{2,2}(\Gamma,V_{\cdot})\ar[r]^{\quad \delta_{\chi}}\ar[u]^{\check\delta}&\cdots\\
 C^{1,0}(\Gamma,V_{\cdot})\ar[r]^{\delta_{\chi}}\ar[u]^{\check\delta}&
 C^{1,1}(\Gamma,V_{{\cdot}})\ar[r]^{\delta_{\chi}}\ar[u]^{\check\delta}&
 C^{1,2}(\Gamma,V_{{\cdot}})\ar[r]^{\quad \delta_{\chi}}\ar[u]^{\check\delta}&\cdots\\
 C^{0,0}(\Gamma,V_{{\cdot}})\ar[r]^{\delta_{\chi}}\ar[u]^{\check\delta}&
 C^{0,1}(\Gamma,V_{{\cdot}})\ar[r]^{\delta_{\chi}}\ar[u]^{\check\delta}&
 C^{0,2}(\Gamma,V_{{\cdot}})\ar[r]^{\quad \delta_{\chi}}\ar[u]^{\check\delta}&\cdots,
}
\end{eqnarray}
where $C^{k,l}(\Gamma,V_{\boldsymbol\cdot})$ is the set of all continuous functions
$\Gamma^l\times \coprod_{\alpha\in I^{k+1}} V_\alpha\to \U(1)$, for
$V_\alpha:= V_{i_0}\cap\dots\cap V_{i_k}$ if $\alpha=(i_0,\dots,i_k)$.
The boundary operators $\delta_{\chi}$ and $\check\delta$ are the usual
boundary operators of group and $\check{\rm C}$ech cohomology.
The corresponding total cohomology of this complex
is denoted by $\check H^\bullet_\Gamma(V_{\boldsymbol\cdot},\underline{\U(1)})$,
and we define the colimit over all $\Gamma$-equivariant coverings
\begin{eqnarray*}
\check H^\bullet_\Gamma(G/N,\underline{\U(1)}):=
{\rm colim}_{V_{\boldsymbol\cdot}}\check H^\bullet_\Gamma(V_{\boldsymbol\cdot},\underline{\U(1)}).
\end{eqnarray*}
For an equivariant bundle isomorphism
\begin{eqnarray*}
\xymatrix{
\lambda {\rm Ad}(l)\ \rotatebox{90}{\rotatebox{90}{\rotatebox{90}{$\circlearrowright$}}}
\hspace{-0.9cm}& 
G/N\times\PU(\HH)\ar[d]\ar[r]^v
&G/N\times\PU(\HH)\ar[d]&
\hspace{-1.2cm} \rotatebox {90}{$\circlearrowleft$}\ \lambda\\
\chi\ \rotatebox{90}{\rotatebox{90}{\rotatebox{90}{$\circlearrowright$}}}
\hspace{-3.4cm}
&
G/N\ar[r]^=\ar[d]&G/N\ar[d]&\hspace{-2.8cm} \rotatebox {90}{$\circlearrowleft$}\ \chi\\
&\ast\ar[r]^=&\ast&
}
\end{eqnarray*}
we are going to define its class
$
[v]\in \check H^1_\Gamma(G/N,\underline{\U(1)}).
$

Let $\{ V_i | i\in I\}$ be a $\Gamma$-equivariant covering 
of $G/N$ such that
the isomorphism $v$ is given 
by functions $v_i: V_i\to\PU(\HH)$
which permit unitary lifts:
\begin{eqnarray*}
\xymatrix{
&&\U(\HH)\ar[d]\\
V_i\ar[rru]^{\overline v_i}\ar[rr]_{v_i}&&\PU(\HH).
}
\end{eqnarray*}
We define  $c_{ij}:V_i\cap V_j\to \U(1)$ by
\begin{eqnarray*}
c_{ij}(z,\hat z):=  v_i(z) v_j(z)^{-1}.
\end{eqnarray*}

Now, for the cocycles  $\lambda,\lambda{\rm Ad}(l):\Gamma\times G/N \to\PU(\HH)$ 
we have
\begin{eqnarray*}
\lambda(a,z) v_i(z)=v_i(z+\chi(a))\hat\lambda(a, z){\rm Ad}(l(a,z)),
\end{eqnarray*}
for all $i\in I$. 
This equality can be used to define functions 
$\gamma_i:\Gamma\times V_i\to\U(1)$
by
\begin{eqnarray*}
v_i(z )&=&\lambda(a,z)^{-1}[v_i(z+\chi(a))]
l(a,z) 
 \gamma_i(a,z,\hat z).
\end{eqnarray*}

A short calculation shows that $(c_{\boldsymbol{..}}, \gamma_{\boldsymbol.})$ is a 1-cocycle in 
(\ref{DigGrischssbljab12}), i.e.
\begin{eqnarray*}
\check\delta c_{\boldsymbol{\cdot\cdot}}=1, \quad \delta_{\chi\times\hat\chi} c_{\boldsymbol{\cdot\cdot}}=\check\delta\gamma_{\boldsymbol\cdot},
\quad \delta_{\chi\times\hat\chi}\gamma_{\boldsymbol\cdot}=1.
\end{eqnarray*}
The class of $(c_{\boldsymbol{..}},\gamma_{\boldsymbol .})$
 in $\check H^1_\Gamma(G/N,\underline{\U(1)})$ 
is denoted by $[v]$.
We leave it to the reader to verify that the class 
is independent of the choices of $v_i, \overline v_i$.

\section{The Equivariant $\check{\rm C}$ech Class
$ [\kappa, \kappa']$}
\label{SecTECeCOBI}

\noindent
The definition of  $[\kappa, \kappa']$ is  essentially the 
same as the definition of $[v]$, it just concerns more free parameters.

For  two homomorphisms $\chi:\Gamma\to G/N$ and 
$\hat\chi:\Gamma\to\widehat G/N^\perp$ we 
consider a $\Gamma$-equivariant covering $W_{\boldsymbol\cdot}:=\{W_i|i\in I\}$ of the space 
$G/N\times \widehat G/N^\perp$.
This yields a double complex
\begin{eqnarray}\label{DigGrischssbljab}
\xymatrix{
\vdots&\vdots&\vdots& \\
 C^{2,0}(\Gamma,W_{\cdot})\ar[r]^{\delta_{\chi\times\hat\chi}}\ar[u]^{\check\delta}&
 C^{2,1}(\Gamma,W_{\cdot})\ar[r]^{\delta_{\chi\times\hat\chi}}\ar[u]^{\check\delta}&
 C^{2,2}(\Gamma,W_{\cdot})\ar[r]^{\quad \delta_{\chi\times\hat\chi}}\ar[u]^{\check\delta}&\cdots\\
 C^{1,0}(\Gamma,W_{\cdot})\ar[r]^{\delta_{\chi\times\hat\chi}}\ar[u]^{\check\delta}&
 C^{1,1}(\Gamma,W_{{\cdot}})\ar[r]^{\delta_{\chi\times\hat\chi}}\ar[u]^{\check\delta}&
 C^{1,2}(\Gamma,W_{{\cdot}})\ar[r]^{\quad \delta_{\chi\times\hat\chi}}\ar[u]^{\check\delta}&\cdots\\
 C^{0,0}(\Gamma,W_{{\cdot}})\ar[r]^{\delta_{\chi\times\hat\chi}}\ar[u]^{\check\delta}&
 C^{0,1}(\Gamma,W_{{\cdot}})\ar[r]^{\delta_{\chi\times\hat\chi}}\ar[u]^{\check\delta}&
 C^{0,2}(\Gamma,W_{{\cdot}})\ar[r]^{\quad \delta_{\chi\times\hat\chi}}\ar[u]^{\check\delta}&\cdots,
}
\end{eqnarray}
where $C^{k,l}(\Gamma,W_{\boldsymbol\cdot})$ is the set of all continuous functions
$\Gamma^l\times \coprod_{\alpha\in I^{k+1}} W_\alpha\to \U(1)$, for
$W_\alpha:= W_{i_0}\cap\dots\cap W_{i_k}$ if $\alpha=(i_0,\dots,i_k)$.
The boundary operators $\delta_{\chi\times\hat\chi}$ and $\check\delta$ are the usual
boundary operators of group and $\check{\rm C}$ech cohomology.
The corresponding total cohomology of this complex
is denoted by $\check H^\bullet_\Gamma(W_{\boldsymbol\cdot},\underline{\U(1)})$,
and we define the colimit over all $\Gamma$-equivariant coverings
\begin{eqnarray*}
\check H^\bullet_\Gamma(G/N\times \widehat G/N^\perp,\underline{\U(1)}):=
{\rm colim}_{W_{\boldsymbol\cdot}}\check H^\bullet_\Gamma(W_{\boldsymbol\cdot},\underline{\U(1)}).
\end{eqnarray*}

Consider two topological triples 
$\big(\kappa,(P,E),(\widehat P,\widehat E)\big)$ 
and
$\big(\kappa',(P',E'),(\widehat P',\widehat E')\big)$
(with underlying Hilbert spaces $\HH,\HH'$ respectively)
such that the underlying 
pairs $(P,E)$ and  $(P',E')$ and the dual pairs 
$(\widehat P,\widehat E)$
and $(\widehat P',\widehat E')$ are outer conjugate.
We are going to define the class
$$
[\kappa,\kappa']\in \check H^1_\Gamma(G/N\times \widehat G/N^\perp,\underline{\U(1)}).
$$

Let $\{ W_i | i\in I\}$ be a $\Gamma$-equivariant covering 
of $G/N\times\widehat G/N^\perp$ such that
the isomorphisms $\kappa, \kappa'$ are given 
by functions $\kappa_i,\kappa_i': W_i\to\PU(\HH)$
which permit unitary lifts:
\begin{eqnarray*}
\xymatrix{
&&\U(\HH)\ar[d]\\
W_i\ar[rru]^{\overline\kappa_i,\ \overline\kappa_i'}\ar[rr]_{\kappa_i,\ \kappa_i'}&&\PU(\HH).
}
\end{eqnarray*}
We define $k_i(z,\hat z):=\overline\kappa_i(z,\hat z)\overline\kappa_i'(z,\hat z)^{-1}$
and $c_{ij}:W_i\cap W_j\to \U(1)$ by
\begin{eqnarray*}
c_{ij}(z,\hat z):=  k_i(z,\hat z) k_j(z,\hat z)^{-1}.
\end{eqnarray*}

Now, if $\lambda,\lambda':\Gamma\times G/N \to\PU(\HH)$ and
$\hat\lambda,\hat\lambda':\Gamma\times \widehat G/N^\perp\to\PU(\HH)$
are the cocycles implementing the $\Gamma$-actions,
we have
\begin{eqnarray*}
\lambda(a,z)\kappa_i(z,\hat z)&=&\kappa_i(z+\chi(a),\hat z+\hat\chi(a))\hat\lambda(a,\hat z),\\
\lambda'(a,z)\kappa'_i(z,\hat z)&=&\kappa'_i(z+\chi(a),\hat z+\hat\chi(a))\hat\lambda'(a,\hat z),
\end{eqnarray*}
 for all $i\in I$. By use of the outer equivalence of the actions, 
 $\lambda'=\lambda {\rm Ad}(l), \hat\lambda'=\hat\lambda{\rm Ad}(\hat l)$,
 we eliminate $\hat\lambda$ and obtain
\begin{eqnarray*}
\kappa_i(z,\hat z)\kappa_i'(z,\hat z)^{-1}&=&\lambda(a,z)^{-1}
\kappa_i(z+\chi(a),\hat z+\hat\chi(a))\kappa_i'(z+\chi(a),\hat z+\hat\chi(a))^{-1}\lambda(a,z)
\\
&&{\rm Ad}(l(a,z)) 
\kappa_i'(z,\hat z){\rm Ad}(\hat l(a,\hat z))^{-1}\kappa_i'(z,\hat z)^{-1}.
\end{eqnarray*}
This equality can be used to define functions 
$\gamma_i:\Gamma\times W_i\to\U(1)$
by
\begin{eqnarray*}
k_i(z,\hat z)&=&\lambda(a,z)^{-1}[k_i(z+\chi(a),\hat z+\hat\chi(a))]
l(a,z) 
\kappa_i'(z,\hat z)[\hat l(a,\hat z))^{-1}]\ \gamma_i(a,z,\hat z).
\end{eqnarray*}

A short calculation shows that $(c_{\boldsymbol{..}}, \gamma_{\boldsymbol.})$ is a 1-cocycle in 
(\ref{DigGrischssbljab}), i.e.
\begin{eqnarray*}
\check\delta c_{\boldsymbol{\cdot\cdot}}=1, \quad \delta_{\chi\times\hat\chi} c_{\boldsymbol{\cdot\cdot}}=\check\delta\gamma_{\boldsymbol\cdot},
\quad \delta_{\chi\times\hat\chi}\gamma_{\boldsymbol\cdot}=1.
\end{eqnarray*}
The class of $(c_{\boldsymbol{..}},\gamma_{\boldsymbol .})$ in $\check H^1_\Gamma(G/N\times\widehat G/N^\perp,\underline{\U(1)})$ 
is denoted by $[\kappa,\kappa']$, and we leave 
it to the reader to check that it is independent of the 
chosen charts and the chosen lifts.
\end{appendix}

\end{document}